\documentclass{amsart}%
\usepackage{amsmath,amsthm}
\usepackage{amsfonts}
\usepackage{amssymb}
\usepackage{graphicx}
\usepackage[usenames,dvipsnames]{color}

\usepackage[english]{babel}

\usepackage[utf8]{inputenc}

\usepackage{verbatim}
\usepackage[usenames,dvipsnames]{color}

\usepackage{amsmath, amssymb,latexsym,amsthm}

\usepackage{enumitem}
\usepackage{moreenum}

\theoremstyle{definition}

\numberwithin{example}{section}

\newcommand{\ul}{\ulcorner}
\newcommand{\ur}{\urcorner}
\newcommand{\U}{\mathcal{U}}
\newcommand{\V}{\mathcal{V}}

\newcommand{\modRr}{\mod \llangle \Rel_0\cup\ldots\cup \Rel_r\rrangle}

\newcommand{\fracn}{\frac{n}{2}}
\newcommand{\Z}{\mathcal{Z}}

\theoremstyle{plain}
\newtheorem{theorem}{Theorem}
\numberwithin{theorem}{section}

\newtheorem{lemma}{Lemma}
\numberwithin{lemma}{section}

\newtheorem{proposition}[lemma]{Proposition}
\newtheorem{corollary}[lemma]{Corollary}
\newtheorem{definition}[lemma]{Definition}
\newtheorem{examples}[lemma]{Examples}

\theoremstyle{definition}

\newtheorem{remark}[lemma]{Remark}
\newtheorem{convention}[lemma]{Convention}

\usepackage{tikz}
\usetikzlibrary{decorations.markings}
\usetikzlibrary{decorations.pathreplacing}
\usetikzlibrary{shapes.geometric,positioning,calc}
\usetikzlibrary{intersections}
\usetikzlibrary{positioning}

\newcommand{\Fr}{\mathrm{F}}
\newcommand{\Rel}{\mathrm{Rel}}
\newcommand{\Can}{\mathrm{Can}}
\newcommand{\Canc}{\mathrm{Cycl}}
\newcommand{\Scan}{\mathrm{SCan}}
\newcommand{\FracM}{\Lambda}
\newcommand{\can}{\mathrm{can}}

\renewcommand{\prod}{\mathrm{prod}}

\newcommand{\llangle}{\langle \langle }
\newcommand{\rrangle}{\rangle \rangle }
\newcommand{\half}{{\lceil\frac{t}{2}\rceil}}
\newcommand{\Relr}{\llangle \Rel_0\cup\ldots\cup \Rel_r \rrangle}

\newcommand{\deglex}{\mathrm{deglex}}
\renewcommand{\leq}{\leqslant}
\renewcommand{\geq}{\geqslant}
\renewcommand{\epsilon}{\varepsilon}
\renewcommand{\emptyset}{\varnothing}
\renewcommand{\phi}{\varphi}
\newcommand{\abs}[1]{|{#1}|}
\newcommand{\inv}{^{-1}}
\newcommand{\X}{\mathcal{X}}
\newcommand{\Y}{\mathcal{Y}}
\renewcommand{\tilde}{\widetilde}
\renewcommand{\hat}{\widehat}

\tikzset{arrow/.style={
        decoration={markings,
            mark= at position #1 with {\arrow{stealth}},
        },
        postaction={decorate}
    }
}

\tikzset{reversearrow/.style={
        decoration={markings,
            mark= at position #1 with {\arrow{stealth reversed}},
        },
        postaction={decorate}
    }
}

\begin{document}
\title{The Burnside Problem for Odd Exponents}
\date{\today}

\author{Agatha Atkarskaya}
\thanks{AA:  atkarskaya.agatha@gmail.com}
\author{Eliyahu Rips}
\thanks{ER: Hebrew University, Jerusalem,  eliyahu.rips@mail.huji.ac.il}

\author{Katrin Tent}
\thanks{KT: Universit\"at M\"unster,  tent@wwu.de}
 \thanks{The first and third author were partially supported by the Deutsche Forschungsgemeinschaft (DFG, German Research Foundation) under Germany's Excellence Strategy EXC 2044-390685587, Mathematics M\"unster: Dynamics-Geometry-Structure and SFB 1442, and by their stay at the Oberwolfach Research Institute. The first author was supported by an ISF fellowship and wishes to thank E. Plotkin for constant support and encouragement.}

\maketitle

\begin{abstract}We show that the free Burnside groups $B(n,m)$ are infinite for $n\geq 557$ and $m\geq 2$. The proof uses iterated small cancellation theory where the induction based on the nesting depth of relators. The main instrument at every step is a new concept of a certification sequence. This decreases the best known lower bound in the Burnside problem for odd exponents from~$665$ to~$557$.
\end{abstract}

\section{Introduction}
\label{intro_section}

In 1902 Burnside asked whether any finitely generated group of finite exponent is necessarily finite. This
question was first answered in the negative in 1964 by Golod and Shafarevitch who constructed an infinite
finitely generated torsion group. However, their example has unbounded exponent raising the question
whether the so-called free Burnside group
\[B(m, n) = F_m/\llangle  w^n\colon w \in F_m\rrangle \]
of exponent $n$ is finite where $F_m$ is the free group in $m$ generators.
For exponents $n = 2, 3, 4$ and $6$ it is known by work of Burnside \cite{B}, Sanov \cite{Sa}, and M. Hall \cite{Ha} that the free Burnside
group is indeed finite for any finite number $m$ of generators.
On the other hand, in 1968 Adian and Novikov gave the first proof that the free Burnside group $B(m, n)$
is infinite for odd $n \geq 4381$ \cite{AN}. Later on, Adian improved the bound to odd $n \geq 665$ \cite{Ad}. The case of even exponent $n$ turned out to be much harder. This case was treated by Ivanov in 1992 \cite{Iv}, he established that $B(m, n)$ is infinite for $n > 2^{48}$ \cite{Iv}. Then Lysenok in 1996 improved the exponent for the even case to $n \geq 8000$ \cite{Lys}. Together with the work of Adian~\cite{Ad}, this yields that $B(m , n)$ is infinite for all $m > 1$ and all $n \geq 8000$.
The proofs of Adian and Novikov use a very involved induction process with a  list of 178 assumptions. So Ol’shanskii's  geometric proof based on a deep study of van-Kampen diagrams was an important step. It resulted in the paper \cite{Ols} for exponents $n > 10^{10}$. The proof is much shorter and more transparent than the one by Adian and Novikov, at the expense of a significantly larger exponent.

Another more geometric approach to free Burnside groups of odd exponent was suggested by Gromov and Delzant in~\cite{GrDl}. This has been further developed by Coulon \cite{Cu}. However, their arguments also require a very large exponent~$n$.

Note that the restricted version of the Burnside problem asks whether there exist finitely many different finite groups in $m$ generators of exponent~$n$, up to isomorphism. This question was solved in positive by Zelmanov in 1989, \cite{Zel1}, \cite{Zel2}, for arbitrary exponents.

While arguably the Burnside question has thus long been settled,  the precise lower bound for the infiniteness of $B(m, n)$ remains open. Experience shows that decreasing the exponent requires huge efforts even for small steps. We  hope that our methods pave the way for further reductions and we believe that an  exponent around $300$ might be in reach.

We also believe that it is  important to provide readable and accessible proofs which give useful lower bounds for the infiniteness of $B(m, n)$ and that the methods developed in this paper are applicable for addressing other Burnside type questions, for instance Engel and quasi-Engel problems, which deal with identities equal to Engel and quasi-Engel words.

 Our proof works  inductively by choosing
a canonical representative for every coset in $B(m, n)$. The induction is based on the rank of a word 
$w \in F_m$ where we (roughly speaking) define the rank $rk(w)$ to be greater or equal to $k + 1$ with respect to our \emph{nesting constant} $\tau$ if the word $w$
(cyclically) contains a subword of the form $v^{\tau}$
for some word $v\in F_m$ with $rk(v) \geq k$.

We define
\[N_k=\llangle 
w^n
\colon rk(w) \leq k\rrangle .\]

Thus, we obtain an ascending sequence of normal subgroups
\[N_0\leq N_1\ldots \leq N_i\leq\ldots\bigcup N_i=N=\llangle  w^n\colon w \in F_m\rrangle .\]
We inductively define the canonical form $\can_k(w)$ for a word $w$ as a canonical representative for $wN_k$. In particular, for all
$w_0, w_1 \in F_m$ we  have \[w_0N_k = w_1N_k \mbox{ \ if and only if \  } \can_k(w_0) = \can_k(w_1)\] 
and we can define a group operation on the set of canonical forms of rank~$k$ making this group isomorphic to $F_m/N_k$.

In order to define the canonical form $\can_k(w)$ on the basis of $\can_{k-1}(w)$ we use the concept of a \emph{certification sequence}. 
 We think of it as carefully choosing the sides of the relators in a given word. 
The important point is that for any $w \in F_m$, the
canonical form stabilizes, i.e. for any $w\in F_m$ there is some $k$ such that $\can_k(w) =\can_l(w)$ for all $l > k$ and thus $\can_k(w)$ will be the canonical representative for $wN\in B(m,n)$.

In this way we obtain a section
$\can: F_m/N \longrightarrow F_m $ i.e. we have \[ \can(w)=\can(w')  \mbox{\  if and only if \ }
wN=w'N\in B(m,n).\]

The set of canonical forms $\can(F_m)$ with the appropriate multiplication then forms a group isomorphic to $B(m,n)$.

Thus, the main thrust of the paper lies in inductively defining $\can_k(w)$ for any $k$ based on 13 induction hypotheses.
We will see that any cube-free element of $F_m$ is
already in canonical form and so the infinity of the Burnside group  follows immediately from the
fact that there are infinitely many cube-free words on two letters.

For our method to give a relatively short and accessible proof, we currently need the exponent $n$ to be
at least $n > 36\cdot 15+16 =556$. However, we expect that this can still be much improved.
The proof also yields (the previously known result) that the infinite free Burnside groups are not finitely
presented.

\section{The set-up}\label{sec:set-up}

Let $\Fr = \langle  x_1, \ldots, x_m \rangle $ be the free group with free generators $x_1, \ldots, x_m$, $m\geq~2$.  Then
\begin{equation*}
B(m, n) = \Fr /\llangle  x_1, \ldots, x_m \mid w^n, \ w\in \Fr\ \rrangle 
\end{equation*}
is called \emph{the free Burnside group of rank $m$ and exponent $n$}.

In this paper we prove the following
\begin{theorem}
\label{main thm}
The free Burnside group $B(m, n)$ is infinite for $m \geqslant 2$ and odd exponents $n \geq 557$.
\end{theorem}
Throughout the paper $n$ is an odd natural number $\geqslant 557$. Section~\ref{sec:IH} and Section~\ref{ind_step_section} describe an inductive process for the definition of a canonical form. We apply the results of this induction  in Section~\ref{proof_completion_section} and show that $B(m, n)$ is infinite for $m \geqslant 2$ and odd exponents $n \geq 557$.
\bigskip

The free generators $\lbrace x_1, \ldots, x_m\rbrace$  of $\Fr$ and their inverses $\lbrace x_1\inv,\ldots , x_m\inv  \rbrace $ are called~\emph{letters}, sequences of letters are called~\emph{words}. A word without cancellations is called a~\emph{reduced word}.

We say that a word $w$ \emph{cyclically contains a word} $A$ if $A$ is a subword of a cyclic shift of~$w$.

 A \emph{prefix} of a reduced word is any (not necessarily proper) initial segment of this word. Similarly, a \emph{suffix} of a reduced word is any (not necessarily proper) final segment of it.

If $N$ is a normal subgroup of $G$ and  $w_1, w_2\in G$ represent the same element in~$G/ N$, we say that $w_1$ and $w_2$ \emph{are equivalent} $\mod N$ and we write \[w_1\equiv w_2 \mod  N.\] We write $w_1=w_2$  to   denote equality of (reduced) words in the free group.

Let $A, B\in \Fr$. We denote their product by $A\cdot B$. If we just write $AB$ this implies that $A\cdot B$ has no cancellation. In particular, if we write $A^m$ for some exponent $m\in\mathbb{Z}$,  this indicates that $A$ is cyclically reduced. 

We will frequently use the following easy observation:
\begin{remark}\label{rem:cancellation}
Suppose that $A$ and $B$ are reduced words. Then
the product $A\cdot B\inv $ has cancellation if and only if $A$ and $B$ have a non-trivial common suffix. Similarly, $A\inv \cdot B$ has cancellation if and only if $A$ and $B$ have a non-trivial common prefix.
\end{remark}

For any word $w$ we denote the number of letters in~$w$ by $| w|$ and call it the \emph{length of $w$}.

\begin{remark}\label{rem:centralizers in free groups}
Note that if $w\neq 1$ is a reduced word in the free group, then $\mathrm{Cen}(w^n)=\langle  w\rangle $ if and only if $w$ is not a proper power. In this case we say that $w$ is \emph{primitive}.
\end{remark}

\section{The list of induction hypotheses}
\label{sec:IH}

The purpose of the induction is to define
 \emph{the canonical form of rank~$i$}, $\can_i(A)$, of $A$,   for all words $A$  in the alphabet $\lbrace x_1, \ldots, x_m \rbrace \cup \lbrace x_1\inv , \ldots, x_m\inv  \rbrace$ and for all $i\geq 0$. Then $\Can_i$ denotes the set of canonical forms of rank $i$.
To start the induction, $\Can_{-1} $ is the set of all words (not necessarily reduced) and the canonical form of rank~$0$ of a word in $\Can_{-1}$ is its reduced form, i.e. $\Can_0$ is the set of all reduced words. Then we inductively define the  canonical form of rank $i$
for all $A \in \Can_{i-1}$ and extend the definition to all words in $\Can_{-1}$ via 
\[\can_i(A)=\can_i(\can_{i - 1}(\ldots \can_0(A)\ldots)).\]
 The elements of $\Can_i$ are called \emph{canonical words of rank~$i$}.

Furthermore, we will specify pairwise disjoint sets $\Rel_i\subset\{w^n\colon w\in\Fr$ primitive$\}$ of relators which are invariant under inverses and cyclic shifts.
Note that relators from $\Rel_i$  may not belong to $\Can_{i - 1}$.

Throughout the paper we fix our \emph{nesting constant} $\tau = 15$.

\begin{definition}[Fractional powers and $\FracM_{i}$-measure]
\label{def:frac power} 
If $u$ is a subword of $a^k$ for some $k\in\mathbb{Z}$, we call $u$ a \emph{fractional power of $a$} and put
\[ \FracM_a(u) = \frac{| u|}{| a|}.\]
If $a^n \in \Rel_i$, we call $u$ a \emph{fractional power of rank} $i$ and if $k\geq \tau+1$ we  put
$\FracM_i(u) =  \FracM_a(u)$.
If $k<\tau+1$  we only define its $\FracM_{i}$-measure if it is clear from the context with respect to which relator from $\Rel_{i}$ the measure is taken.

We say that $u$ has $\FracM_{i}$-measure \emph{at most} $m$ for $m\geq\tau$ if either $\FracM_{i}(u)\leq m$  or the $\FracM_i$-measure of $u$ is not defined.
\end{definition}
 We show inductively for $i \geqslant 0$  that $\Can_i$ is a group with respect to an appropriately defined  multiplication.

\textbf{The induction hypothesis at stage $r$:} At stage $r$ we assume inductively that the following statements hold for $i = 0, \ldots, r - 1$.
Here and in what follows we will refer to Induction Hypothesis~1 as IH~1 etc.
\begin{enumerate}[label=IH\,\arabic*., ref=\arabic*]
\item
\label{IH can unique}
The canonical form of rank~$i$ of every word of $\Can_{i - 1}$ is uniquely defined and \[\Can_i=\{ \can_i(w)\mid\ w\in\Can_{i-1}\}.\]
\item
\label{IH  can decreasing sequence}
$\Can_{i} \subseteq \Can_{i - 1}$.

\item
\label{IH Rel}
The sets $\Rel_i$, $0\leq i \leq r-1$, are closed under cyclic shifts and inverses and pairwise disjoint.
We have $\Rel_0 =\{1\}$,  and $\Rel_i \subseteq\{ w^n\mid w\in \Fr \textit{ primitive}\}$  for $1\leq i\leq r-1$.

\item
\label{IH immediate}
If $A \in \Can_{i - 1}$ does not contain fractional powers of rank~$i$ of $\Lambda_i$-measure $> \fracn - 5\tau - 2$, then $A \in \Can_i$.
\end{enumerate}
\begin{remark}

Note that by IH~\ref{IH immediate} we also have
$\can_i(1) = 1\in\Can_i$ where $1$ denotes the empty word, and
 $\can_i(x)=x \in \Can_i$ for every single letter $x$.
\end{remark}

The small cancellation condition  is contained in the following induction hypothesis (see Lemma~\ref{common_part_of_powers}):
\begin{enumerate}[label=IH\,\arabic*., ref=\arabic*]
\setcounter{enumi}{4}
\item
\label{rel_common_parts_in_diff_ranks_hyp}\label{IH small cancellation}
Let $x^n \in \Rel_i,\ y^n\in \Rel_j$, $1\leq i\leq j\leq r-1 $,  and  let $c$ be their common prefix. If $i<j$, then $|c|<2|y|$ and if $i=j$ and $|x|\leq |y|$, then $|c|< \min\{(\tau + 1)|x|, 2|y|\}$.

\item
\label{IH idempotent}
If $A \in \Can_i$, then $A = \can_i(A)$.

\item
\label{IH can inverse}
$\can_i(A\inv ) = (\can_i(A))\inv $.
\end{enumerate}

\begin{remark}\label{rem:group IH}
By IH~\ref{IH idempotent} we have $\can_i(B) = \can_i(\can_i(B))$ for every $B \in \Can_{-1}$. In other words, $\can_i$ is an idempotent operation equivariant with respect to taking inverses by IH~\ref{IH can inverse}.
\end{remark}

The following axiom states that the canonical form picks unique coset representatives:
\begin{enumerate}[label=IH\,\arabic*., ref=\arabic*]
\setcounter{enumi}{7}
\item
\label{IH can of equivalent words}
Let $A, B \in \Can_{-1}$. Then $A\equiv B\mod \llangle \Rel_0, \ldots, \Rel_i\rrangle$  if and only if $\can_i(A) = \can_i(B)$.
\end{enumerate}

\begin{remark}
\label{can_equiv_remark}\label{can_preserves_coset_hyp} 

Note that for $i \geqslant 0$ the set $\Can_i$ is a group with respect to the multiplication defined by
\begin{equation*}
A\cdot_i B = \can_i(A\cdot B),\ A, B \in \Can_i,
\end{equation*}
with identity element $1=\can_i(1)$ and inverses given by inverses in the free group.

In particular, $(\Can_0, \cdot_0)$ is precisely the free group~$\Fr$.

Notice that if $A \equiv B \mod \llangle  \Rel_0, \ldots, \Rel_i\rrangle$ and $A \in \Can_i$, then $\can_i(B) =\can_i(A)= A$ by IH~\ref{IH can of equivalent words} and IH~\ref{IH idempotent}. 

For $A \in \Can_{-1}$ we thus have
 $A\equiv\can_i(A)\mod \llangle\ \Rel_0, \ldots, \Rel_i\ \rrangle$.
Furthermore, since $\can_i(1)=1$,
 we have $\can_i(v) = 1$ for $v \in \Rel_i$, $1\leq i\leq r-1$.
\end{remark}

These previous remarks can be rephrased as:

\begin{corollary}
\label{can_of_can_product_hyp}
Let $A, B \in \Can_{-1}$. Then for $i \geqslant 0$ we have
\begin{align*}
\can_i(A) \cdot_i \can_i(B) &= \can_i(\can_i(A)\cdot \can_i(B)) = \can_i(A\cdot \can_i(B)) \\
&=\can_i(\can_i(A)\cdot B) = \can_i(A\cdot B).
\end{align*}
\end{corollary}

\begin{enumerate}[label=IH\,\arabic*., ref=\arabic*]
\setcounter{enumi}{8}

\item
\label{can_non_trivial_subwords_hyp}
Any non-empty subword of a word from $\Can_i$, $i \geqslant 0$, is not equal to $1$ in the group $\Fr / \llangle\ \Rel_0, \ldots, \Rel_i\ \rrangle$.
\end{enumerate}

\begin{definition}
\label{alpha_free_def}
A reduced word $A$ is $\alpha$-free modulo rank~$i$ if $A$ it does not contain subwords of the form $a^{\alpha}$ where $a$ is primitive and $a^n \notin \Rel_1\cup\ldots \cup \Rel_i$.

A reduced word $A$ is $\alpha$-free of rank~$i$ if it does not contain subwords of the form $a^{\alpha}$ with $a^n \in \Rel_i$.

\label{def:canonical triangle} We call a triple of words $(D_1, D_2, D_3)$ a \emph{canonical triangle} of rank~$i$ if they are $\tau$-free modulo rank $i+1$ and
$D_1\cdot D_2\cdot D_3\inv\equiv 1 \mod \llangle \Rel_0, \ldots, \Rel_i\rrangle$.
\end{definition}

The following axiom is crucial:
\begin{enumerate}[label=IH\,\arabic*., ref=\arabic*]
\setcounter{enumi}{9}
\item
\label{can_triangle_hyp1} {\bf (Canonical triangle hypothesis)}
For $A, B \in \Can_i$  there is a canonical triangle $(D_1, D_2, D_3)$ of rank~$i$ such that  $A = A'D_1X,  \ B = X\inv D_2B'$ (where $X\cdot X\inv $ is the maximal cancellation in $A\cdot B$) such that \[\can_i(A\cdot B)= A'D_3B'.  \] 
Furthermore, if $(D_1^{(i)}, D_2^{(i)}, D_3^{(i)})$ is a canonical triangle of rank~$i-1$ such that  $A = A''D_1^{(i)}X,  \ B = X\inv D_2^{(i)}B''$ and $\can_{i-1}(A\cdot B)= A''D^{(i)}_3B''$, then $A'$ is a prefix of $A''$ and $B'$ is a suffix of $B''$ and if $D_1=D_1^{(i)}, D_2=D_2^{(i)}$, then $D_3=D_3^{(i)}$. 
\end{enumerate}
Note that if $D_1=D_1^{(i)}, D_2=D_2^{(i)}$, then $A'=A''$ and $B'=B''$ since the maximal cancellation is independent of $i$.

The multiplication $A\cdot_i B = \can_i(A\cdot B)$ in the group $(\Can_i, \cdot_i)$ can be graphically expressed as follows:
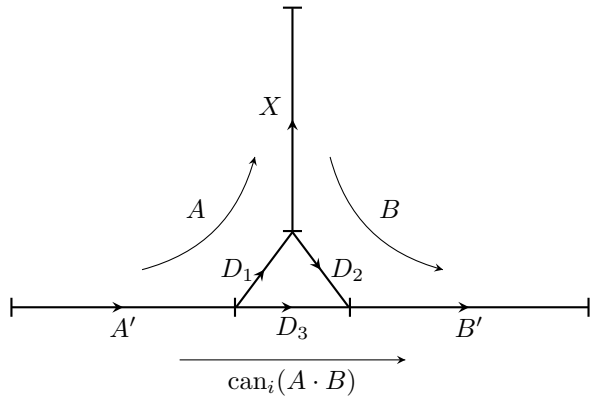
\begin{center}
\begin{tikzpicture}
\draw[|-|, thick, arrow=0.5] (0, 0) to node[midway, below] {$A'$} (3, 0);
\draw[thick, arrow=0.5] (3, 0) to node[midway, below] {$D_3$} (4.5, 0);
\draw[|-|, thick, arrow=0.5] (4.5, 0) to node[midway, below] {$B'$} (7.7, 0);
\draw[thick, arrow=0.5] (3, 0) to node[midway, left] {$D_1$} (3.75, 1);
\draw[thick, arrow=0.5] (3.75, 1) to  node[midway, right] {$D_2$}(4.5, 0);
\draw[|-|, thick, arrow=0.5] (3.75, 1) to node[midway, left, yshift=5] {$X$} (3.75, 4);

\draw[black, arrow = 1] (-1.5 + 3.75, -0.7)--(1.5 + 3.75, -0.7) node[midway, below] {$\can_i(A\cdot B)$};
\draw[black, arrow = 1] (-2 + 3.75, 0.5) to [bend right] node[midway, left, yshift=8] {$A$} (-0.5 + 3.75, 2);
\draw[black, arrow = 1] (0.5 + 3.75, 2) to [bend right] node[midway, right, yshift=8] {$B$} (2 + 3.75, 0.5);

\end{tikzpicture}
\end{center}
Note that $A\cdot B$ and $\can_i(A\cdot B)$ represent the same element  in $\Fr / \llangle \Rel_0, \ldots, \Rel_i\rrangle $ by IH~\ref{IH can of equivalent words}. Hence after cancelling  $A'$ from the left and $B'$ from the right it  follows that $D_1\cdot D_2$ and $D_3$ represent the same element in $\Fr / \llangle \Rel_0, \ldots, \Rel_i\rrangle $.
In particular, if two of $D_1, D_2, D_3$ are equal to $1$, then so is the remaining one by IH~\ref{can_non_trivial_subwords_hyp}.

The triangles constitute the '\emph{smoothing process}' in the multiplication of canonical words. So IH~\ref{can_triangle_hyp1} states that in this smoothing process the perturbation on both sides of the multiplication seam is very limited and, furthermore, in order to obtain higher canonical forms the \emph{smoothing area} given by the canonical triangles  may need to increase (but will never shrink).
\noindent

\begin{enumerate}[label=IH\,\arabic*., ref=\arabic*]
\setcounter{enumi}{10}

\item
\label{can_periodic_extension_constriction_hyp}
If $L_1A^{\tau}R_1,\  L_2 A^{\tau} R_2 \in \Can_i$ for $A$ primitive, $A^n\notin \Rel_0\cup\ldots\cup\Rel_i$ then $L_1A^NR_2 \in \Can_i$ for any $N\geq\tau$.

\item
\label{can_of_prefix_and_suffix_hyp}

If $A_1$ is a prefix of $A \in \Can_i$, there is a canonical triangle $(D_1, 1, D_3)$ such that $A_1=A_1'D_1$ and $ \can_{i}(A_1) = A_1'D_3$.
\end{enumerate}
By taking inverses IH~\ref{can_of_prefix_and_suffix_hyp} implies also
that for a suffix $A_2$ of $A\in\Can_{r-1}$ there is a canonical triangle $(E_1, E_2, 1)$ such that $\can_{i}(A_2) = E_3A_2'$ and $A_2=E_2A_2'$.

\begin{enumerate}[label=IH\,\arabic*., ref=\arabic*]
\setcounter{enumi}{12}
\item
\label{can_power_hyp}
If  $A \in \Can_{-1}$ and $A^n\notin\llangle \Rel_0, \ldots, \Rel_i\rrangle$, then there are natural numbers $K, M_0$ and words $W, Z$  depending only on $A$ and $i$ such that 
\begin{equation*}
\can_i(\underbrace{A\cdot \ldots \cdot A}_{M \textit{ times}}) = W\widetilde{A}^{M - K}Z \ \textit{ for all } M \geqslant M_0,
\end{equation*}
 and $A$ and $\widetilde{A}$ are conjugate in the group $\Fr / \llangle \Rel_0, \ldots, \Rel_i\rrangle$.
\end{enumerate}

We now collect a few immediate consequences of the induction hypotheses which will be widely used throughout:

\begin{corollary}
\label{can_power_inverse_hyp_cor}
Let $La^{N_1}Aa^{N_2}R\in\Can_i$ where $A$ may be empty, $a$ is primitive, $a^n \notin \Rel_0\cup\ldots \cup \Rel_i$ and $N_1, N_2\geq2\tau$. Then
\begin{align*}
\can_i(La^{N_1}) &= La^{N_1 - \tau}X,\\
\can_i(a^{N_2}R) &= Ya^{N_2 - \tau}R,\\ 
\can_i(a^{N_1}Aa^{N_2}) &= Ya^{N_1 - \tau}Aa^{N_2 - \tau}X,
\end{align*}
where $X \equiv Y \equiv a^{\tau} \mod \llangle  \Rel_0, \ldots, \Rel_i\rrangle $ and $X, Y$ only depend on $a$ and $i$. 
\end{corollary}
\begin{proof}
 By IH~\ref{can_of_prefix_and_suffix_hyp}, there is a canonical triangle $(D_1, 1, D_3)$ of rank~$i$ such that 
\[La^N=La^{N-\gamma}a_1D_1 \mbox{\  and \ } \can_i(La^N)= La^{N-\gamma}a_1D_3\]
for some $\gamma\leqslant\tau$ and a prefix $a_1$ of $a$. Write $X=a^{\tau-\gamma}a_1D_3$, so $La^{N-\gamma}a_1D_3=La^{N-\tau}X$. 
Since $D_1 \equiv D_3  \mod \llangle \Rel_0, \ldots, \Rel_i\rrangle$, we have
\[a^{\tau}= a^{\tau-\gamma}a_1D_1\equiv a^{\tau-\gamma}a_1D_3=X \mod \llangle \Rel_0, \ldots, \Rel_i\rrangle.\]
Since $N\geq 2\tau$, by IH~\ref{can_periodic_extension_constriction_hyp} we have $L_1a^KX\in\Can_i$ for any $K\geq\tau$ and any $L_1$ such that $L_1a^{\tau}$ is a prefix of a word from $\Can_i$. Now $L_1a^{K+\tau}\equiv L_1a^KX\mod \llangle  \Rel_0, \ldots, \Rel_i\rrangle$, hence by Remark~\ref{can_equiv_remark} we obtain that $\can_i(L_1a^{K+\tau}) = L_1a^KX$. So, $X$ depends only on $a$ and $i$. By taking inverses and applying the previous case on both sides the remaining claims follow. 
\end{proof}

For convenience we also note the following:

\begin{corollary}
\label{can_powers_merging_hyp_cor}
Let $La^{N_1 + N_2}R\in\Can_i$, $N_1 + N_2\geq \tau$, where $a$ is primitive and $a^n \notin \Rel_0\cup\ldots \cup \Rel_i$. Let $M \in \Can_{-1}$ be such that $M \equiv a^{\alpha} \mod \llangle  \Rel_0, \ldots, \Rel_i\rrangle $. Then
\begin{equation*}
\can_i(La^{N_1}\cdot M\cdot a^{N_2}R) = La^{N_1 + N_2 + \alpha}R.
\end{equation*}
\end{corollary}
\begin{proof}
Since $M \equiv a^{\alpha} \mod \llangle  \Rel_0, \ldots, \Rel_i\rrangle $, we see that
\[
La^{N_1}\cdot M\cdot a^{N_2}R \equiv La^{N_1 + N_2 + \alpha}R \mod \llangle  \Rel_0, \ldots, \Rel_i\rrangle .
\]
IH~\ref{can_periodic_extension_constriction_hyp} implies that $La^{N_1 + N_2 + \alpha}R \in \Can_i$. Therefore Remark~\ref{can_equiv_remark} implies the result.
\end{proof}

Since canonical triangles are $\tau$-free of rank $i$,  fractional powers of rank~$i$ and $\Lambda_i$-measure $\geq\tau$ block the influence of the smoothing process obtained from the canonical triangles in the computation of the canonical form for subwords and products:

\begin{corollary}
\label{can_power_context_hyp_cor1}
Let $A = A'D_1X,\ B_1 = X\inv D_2Ma^{\tau}R  \in \Can_i$ and \ $\can_i(A\cdot B_1) = A'D_3Ma^{\tau}R$ for some canonical triangle  $(D_1, D_2, D_3)$  of rank~$i$ and primitive $a$ with $a^n \notin \Rel_1 \cup\ldots \cup \Rel_i$ (where $M$ may be empty). \\
 If $B_2 = X\inv D_2Ma^{\tau}R_1  \in \Can_i$, then $\can_i(A\cdot B_2) = A'D_3Ma^{\tau}R_1$.
 \end{corollary}
  \begin{proof}
   By IH~\ref{can_periodic_extension_constriction_hyp} applied to $A'D_3Ma^{\tau}$ and $a^{\tau}R_1$ we have  $A'D_3Ma^{\tau}R_1 \in \Can_i$. Since $D_1\cdot D_2 \equiv D_3 \mod \llangle  \Rel_1, \ldots, \Rel_i\rrangle $, we see that $A'D_3Ma^{\tau}R_1 \equiv A\cdot B_2 \mod \llangle  \Rel_1, \ldots, \Rel_i\rrangle $. Thus, Remark~\ref{can_equiv_remark} implies the claim.
\end{proof}
Clearly the corresponding statement for $A_1\cdot B, A_2\cdot B$ follows from this by considering inverses.
Similarly we have

\begin{corollary}
\label{can_power_context_hyp_cor2}
Let $A=La^{\tau}Mb^{\tau}WR,\  A_1=L_1a^{\tau}Mb^{\tau}WR_1\in\Can_i$ where $a, b$ are primitive and $a^n, b^n \notin \Rel_1 \cup\ldots \cup \Rel_i$ (where $M, W$ may be empty).  Then
$\can_i(L_1a^{\tau}Mb^{\tau}W)$ is obtained from $\can_i(La^{\tau}Mb^{\tau}W)$ by replacing $L$ by~$L_1$.
\end{corollary}
\begin{proof}
By IH~\ref{can_of_prefix_and_suffix_hyp}  there is a word $D$ $\tau$-free of rank~$i + 1$ such that 
\[\can_i(La^{\tau}Mb^{\tau}W)=La^{\tau}MXD\]
where $X$ is non-empty and $b^{\tau}W\equiv XD\mod\llangle  \Rel_1, \ldots, \Rel_i\rrangle$.
By IH~\ref{can_periodic_extension_constriction_hyp} applied to $L_1a^{\tau}$ and $a^{\tau}MXD$ we have  $L_1a^{\tau}MX D\in \Can_i$. Since
\[L_1a^{\tau}MXD\equiv L_1a^{\tau}Mb^{\tau}W \mod \llangle  \Rel_1, \ldots, \Rel_i\rrangle,\] Remark~\ref{can_equiv_remark} implies the claim.  
\end{proof}

\section{The induction}
\label{ind_step_section}

In this section we start showing that the induction step works. We first establish the 
\textbf{induction basis for $i=0$}. Note that although we have defined $\Can_{-1} $ with index~$-1$, ranks of the canonical form and canonical triangles start from~$0$.

\subsection{Induction basis}

\begin{proposition}
\label{basis_of_induction}
The sets $\Rel_0$, $\Can_0$,  and $\Can_{-1} $ satisfy IH~\ref{IH can unique}--\ref{can_power_hyp}.
\end{proposition}
\begin{proof}
Since $\Rel_0=\{1\}$ and the canonical form of rank~$0$ of a word from $\Can_{-1}$ is its reduced form,   all the induction hypotheses are easily verified. In particular, all sides of canonical triangles of rank $0$ are equal to $1$.
\end{proof}
 Note that IH~\ref{IH immediate},  IH~\ref{IH small cancellation} and IH\ref{rel_common_parts_in_diff_ranks_hyp} are not defined for $i=0$, but will be verified for $i\geq 1$ inside the proofs.

Now assume that IH~\ref{IH can unique}--\ref{can_power_hyp} hold for $\Can_{-1}, \ldots, \Can_{r - 1}$, $\Rel_0, \ldots, \Rel_{r - 1}$. In order to prove \textbf{the  induction step}, we now construct $\Rel_r$ and $ \Can_r $  such that $\Can_{-1}, \ldots, \Can_r$, and $\Rel_0, \ldots, \Rel_r$ also satisfy IH~\ref{IH can unique}--\ref{can_power_hyp}.

\subsection{Cyclically canonical words}
The multiplication of canonical words requires the smoothing process given by canonical triangles at the seam between the words (see IH~\ref{can_triangle_hyp1}). Hence in general $\can_i(A)\cdot\can_i(A)\neq \can_i(A)\cdot_i\can_i(A)=\can_i(A\cdot A)$.
We now define \emph{cyclically canonical words} of rank~$i$, $i = 0, \ldots, r - 1$, for which equality holds at least approximately: 
\begin{definition}[cyclically canonical words]
\label{canc_def}
We say that a word $A$ is \emph{cyclically canonical of rank~$0$} if it is cyclically reduced. 
We call a cyclically reduced word $A$ \emph{cyclically canonical of rank~$i$ for $i \geqslant 1$} if $A^{\tau}$ is a subword of a word in $\Can_i$ and $A=A_1^K$ for a primitive word $A_1^n\notin\Rel_0\cup\ldots\cup\Rel_i$.

The set of all cyclically canonical words of rank~$i$, $i \geqslant 0$, is denoted by~$\Canc_i$.
\end{definition}
Clearly $\can_0(A^K) = A^K$ for every $A\in\Canc_0$ and $K \geqslant 0$ and if $A=A_1^K$ is cyclically canonicalof rank $i$, then so is $A_1$.

The following are immediate consequences of the induction hypotheses:

\begin{lemma}
\label{canc_properties}

\begin{enumerate}
\item
$\Canc_i \subseteq \Canc_{i - 1}$.
\item
$\Canc_i$ is closed under taking cyclic shifts and inverses.
\item
If $A \in \Canc_0$ and $A^{N_1} \in \Canc_i$, then $A^{N_2} \in \Canc_i$ for all $N_1, N_2 \geqslant 1$.

\item
$\Canc_i \cap (\Rel_0\cup \ldots \cup \Rel_i) = \emptyset$.
\item If $A\in\Canc_i$, $i\geq 0$ and $K\geq 4\tau$, then $\can_i(A^K)=T_1 A^{K-2\tau} T_2$ where $T_1, T_2$ only depend on $A$ and $i$.
\end{enumerate}
\end{lemma}

\begin{proof}  (1) and (2) are clear and
(3)  follows from IH~\ref{can_periodic_extension_constriction_hyp}.

(4) This follows directly from the definition.

(5) is a consequence of IH~\ref{can_periodic_extension_constriction_hyp} and Corollary~\ref{can_power_inverse_hyp_cor}.
\end{proof}

\subsection{Sets of relators~$\Rel_r$ and their common parts}
\label{slices_of_relators_section}

Recall that $\Rel_0 = \{1\}$, $\Canc_0$ is the set of all cyclically reduced words and  that throughout the paper we fix the nesting constant $\tau =15$. We put
\begin{align*}
\Rel_1 = \lbrace x^n \in \Canc_0 \mid \ & |x| = 1\rbrace,\\
\Rel_2 = \lbrace x^n \in \Canc_1 \mid \ &  Cen(x)=\llangle  x\rrangle , | x| > 1 \textit{ and }
x \textit{ does not}\\ & \textit{cyclically contain }  a^{\tau}
\textit{for } a \in \Canc_0\setminus\{1\} \big\rbrace.
\end{align*}
For $r\geq 3$ we define:
\begin{align*}
\Rel_r = \big\lbrace x^n \in \Canc_{r - 1} \mid \ &  Cen(x)=\langle  x\rangle    
\textit{ and if  } x \textit{ cyclically contains } a^{\tau}\
\textit{for}\\
& a \in \Canc_0,\ Cen(a)=\langle  a\rangle ,\
\textit{then } a^n \in \Rel_{1} \cup \ldots \cup \Rel_{r - 1} \big\rbrace.
\end{align*}

\begin{remark}
\label{nesting_remark}
Note that by definition for $r\geq 3$, if $x^n\in\Canc_{r-1}$ and $x$ does \emph{not} cylically contain a subword $a^\tau$ with $a^n\in\Rel_{r-1}$, then $x^n\in\Rel_1\cup\ldots\cup \Rel_{r-1}$.
In this way the sets of relators $\Rel_i$ for $i \geqslant 2$ are defined by the nesting depth of power words that contain at least $\tau$ periods (see Corollary~\ref{subwords_r-1_of_rel_r}).

After completing the induction process we prove in Corollary~\ref{rel_covers_all_powers} that by organizing the relators according to their nesting depth,  we obtain 
\[ B(m,n)\cong  \Fr / \left\langle\left\langle \bigcup_i \Rel_i\right\rangle\right\rangle.\]
\end{remark}

Since $\Canc_{r-1}$ is closed under inverses and cyclic shifts, if $x^n \in \Canc_{r - 1}$ is such that $x$ cyclically contains some $a^{\tau}$ with $a^n\in \Rel_{1} \cup \ldots \cup \Rel_{r - 1}$, then so does  $x^{-n}$ and any cyclic shift of  $x^n$. So with Lemma~\ref{canc_properties}(4) and IH~\ref{IH Rel} for ranks $<r$ we obtain:
\begin{lemma}\label{rel_slices_intersection}
{\rm IH~\ref{IH Rel}} holds for $\Rel_r$.
\end{lemma}

\begin{corollary}\label{rem:no nesting in same rank}
If $x^n,y^n\in\Rel_i, i\geq 1$, then $x^{\tau}$ is not cyclically contained in $y$.
\end{corollary}
\begin{proof}
This follows directly from $\Rel_i\cap\Rel_j=\emptyset$ for $i\neq j$ and the definition of $\Rel_i$.
\end{proof}

We also note the following:

\begin{lemma}
\label{lem:nesting_occurrences_rel_r}
If $x\in\Canc_{r-1}$ is primitive, then either $x^n \in \Rel_1\cup\ldots\cup\Rel_r$, or $x$ cyclically contains $a^{\tau}$ for some $a^n \in \Rel_r$ if $r \geqslant 3$, (or  $a^n \in \Rel_1\cup \Rel_2$ if $r = 2$). 
\end{lemma}
\begin{proof}
If
$x$ does not cyclically contain any subwords of the form $a^{\tau}$, then, by definition, $x^n \in \Rel_2$. 
If $x$ cyclically contains only subwords of the form $a^{\tau}$ with $a^n \in \Rel_1\cup\ldots\cup\Rel_{r - 1}$, then, again by definition, $x^n \in \Rel_1\cup\ldots\cup\Rel_r$. So assume $x$ cyclically contains a subword $a^{\tau}$ where $a$ is  primitive and $a^n \notin \Rel_1\cup\ldots\cup\Rel_{r - 1}$. Then $|a|<|x|$, and by Lemma~\ref{canc_properties}(4)  and induction on $|x|$, we have $a^n\in\Rel_r$ or $a$ (and hence $x$) cyclically contains $b^{\tau}$ for some $b^n \in \Rel_r$ if $r \geqslant 3$, (or  $b^n \in \Rel_1\cup \Rel_2$ if $r = 2$). By Corollary~\ref{rem:no nesting in same rank} the cases are mutually exclusive.
\end{proof}
Lemma~\ref{lem:nesting_occurrences_rel_r} and Lemma~\ref{canc_properties}(3) with $r-1$ in place of $r$ now imply:
\begin{corollary}\label{subwords_r-1_of_rel_r}
 If $x^n \in \Rel_r$,
then $x$ cyclically contains $a^{\tau}$ for some $a^n \in \Rel_{r-1}$ if $r \geqslant 3$, (or  $a^n \in \Rel_1\cup \Rel_2$ if $r = 2$).
\end{corollary}

The following important statement is proved in~\cite{FW}.
\begin{lemma}
\label{common_part_of_powers}
Let $x^n, y^n$ be two reduced words such that $x$ and $y$ do not centralize each other in $F$. Let $c$ be a common prefix of $x^n$ and $y^n$. Then $| c| < | x| + | y| - \gcd(| x|, | y|)$, where $\gcd(| x|, | y|)$ is a greatest common divisor of $| x|$ and $| y|$.
\end{lemma}

\begin{lemma}
\label{common_part_in_rank_r}
Let $x^n, y^n \in \Rel_i$, $i \geq 1$, $x \neq y$, and $c$ be a common prefix of $x^n$ and $y^n$. Assume $| x| \leq | y |$. Then $|c| <\min\{ 2|y|, (\tau + 1)| x|\}$.
\end{lemma}
\begin{proof}
For $r=1$, by definition of $\Rel_1$, we have $| x| = | y| = 1$, so the claim is obvious.
Now let $i\geq 2$. Since $x^n\neq y^n \in \Rel_i$ we have $Cen(x)=\langle  x\rangle, Cen(y)=\langle  y\rangle $ and $\langle x\rangle\cap\langle y\rangle=\{1\}$. So it follows from Lemma~\ref{common_part_of_powers} that $| c| < | x| + | y|\leq 2|y|$. 

From $|c|<|y|+|x|$ we see that if $|c|>(\tau+1)|x|$, then we must have $|y|>\tau|x|$. Since $c$ is a common prefix of $x^n$ and $y^n$, this implies that $y$ contains  $x^{\tau}$ as a prefix, contradicting Corollary~\ref{rem:no nesting in same rank}.
\end{proof}

\begin{corollary}
\label{common_part_in_different_ranks}
Let $x, y$ be primitive, $x\in\Canc_{r-1}, x^n \notin \Rel_0\cup\ldots\cup\Rel_i, y^n\in \Rel_i, 1\leq i< r$,  and let $c$ be a common prefix of $x^n$ and $y^n$. Then $|c| <  2|x|$. 
\end{corollary}
\begin{proof}
If  $i=1$, then $y$ is a single letter and since $x$ is primitive, we have $|c|<|x|$. Thus the claim holds for $i=1$.

If $i \geq 2$, then  $x$ cyclically contains a subword $a^{\tau}$ with $a^n \in \Rel_i$  by Corollary~\ref{subwords_r-1_of_rel_r}.   If $|c|\geq 2|x|$, then $|x|<|y|$ by Lemma~\ref{common_part_of_powers} and since any cyclic shift of $x$ is a subword of $x^2$  and hence of $y^2$, we see that  $y$ also cyclically contains $a^{\tau}$, contradicting Corollary~\ref{rem:no nesting in same rank}. 
\end{proof}

Now Lemma~\ref{common_part_in_rank_r},  Corollary~\ref{common_part_in_different_ranks} and $\Rel_i\cap\Rel_j=\emptyset$ for $i\neq j$ directly imply:
\begin{corollary}
\label{ind_hyp_holds_for_rel}
IH~\ref{IH small cancellation} and IH~\ref{rel_common_parts_in_diff_ranks_hyp} hold for $\Rel_i$, $i = 1, \ldots, r$.
\end{corollary}

\subsection{Turns of rank~$r$}
\label{turns_r_section}

If $u$ is a fractional power of $a$, there exists a cyclic shift  $\widehat{a} = a_2a_1$  ($a_1, a_2$  may be empty) of $a=a_1a_2$ such that $u$ can be written in the form
\begin{equation}
\label{frac_power_forms2}
u=\widehat{a}^ka_2 \quad\mbox{  or \ } u = \widehat{a}^{-k}a_1\inv ,\ k \in \mathbb{N} \cup \lbrace 0 \rbrace.
\end{equation}

\medskip


The set  of \emph{fractional powers} of rank $j, 1\leq j\leq r,$ is defined as
\begin{equation*}
\lbrace u \mid u \textit{ is a subword of } R^N, R \in \Rel_j, N \in \mathbb{Z}\rbrace.
\end{equation*}
Note that since $\Rel_j$ is closed under cyclic shifts and inverses, this coincides with 
\begin{equation*}
\lbrace u \mid u \textit{ is a prefix of } R^N, R \in \Rel_j, N \in \mathbb{N}\rbrace.
\end{equation*}

\begin{remark}
\label{subword_of_rel_power} 
If $u$ is a fractional power of rank $j$, $1 \leq j \leq r$ of $\Lambda_j$-measure  $\geqslant \tau + 1$, then by IH~\ref{IH small cancellation} (for $1 \leq j< r $) and Lemma~\ref{common_part_in_rank_r} (for $j = r$) there exists a unique relator  $a^n\in\Rel_j$ such that $u$ is a prefix of $a^K, K\geq 0$. So $u$ can be written uniquely as
\begin{equation}
\label{rel_subword_form}
u = a^ka_1,\textit{ where } a^n \in \Rel_j,\ a = a_1a_2, \ k \in \mathbb{N} \cup \lbrace 0\rbrace.
\end{equation}

Clearly, any fractional power $u$ of rank $j$  can be represented as in~\eqref{rel_subword_form}. However, without  the condition that $u$ contains $\geqslant \tau + 1$ periods of a relator, the relator $a^n\in\Rel_j$ need not be unique, which is why we require in Definition~\ref{def:frac power} that either $k\geq\tau+1$ or that the corresponding relator is clear from the context.
\end{remark}

The following simple definition is a crucial concept for everything that follows:
\begin{definition}[occurrences of rank~$j$, $1 \leq j\leq r$]
\label{max_occurr_rank_r_def}
Let $U$ be a subword of $A\in\Can_{-1}$. Then  \emph{\bf the occurrence of $U$ in $A$} is determined by its position inside $A$.
We say that an occurrence $U$ is properly contained in $A$ if it is neither prefix nor suffix of $A$.

Let $A=LuR \in \Can_{r - 1}$ where $u$  is  a fractional power of rank $j$. If $u$ is not properly contained inside an occurrence $u_1$ in $A$ which is  also a fractional power of rank $j$, then $u$ is called \emph{a maximal occurrence} (of rank~$j$) in $A$.
\end{definition}
I.e.  if $A=LUR=L'UR'$ with $L'\neq L, R'\neq R$, then $A$ contains two different occurrences of $U$.

Note that for any prefix $u$ of $a^k$ with $k\in \mathbb{N}$ and suffix $w$ of $a$, the word $wu$ is reduced and contained in $a^{k+1}$. This motivates the following definition:
\begin{definition}[Prolongation of occurrences of fractional powers]\label{def:prolongation} 
Let  $a^n\in \Rel_r$ and suppose $u,w$ are occurrences in $a^K$ for some  $K\in \mathbb{N}$. If $u$ is properly contained in $w$, we call an occurrence of $w$ in $A\in\Can_{r-1}$ a \emph{prolongation of the occurrence $u$} in $A$.
\end{definition}

\begin{remark}\label{rem:occ with tau_1+1}
If $u = a^ka_1$ with $a^n \in \Rel_r$, then all prolongations of $u$ with respect to $a$ are fractional powers of $a$. If $k\geq\tau+1$, then for prolongations of $u$ we do not have to mention $a$ by Remark~\ref{subword_of_rel_power} as $a$ is unique (up to cyclic shift).
However, if $u$  contains $< \tau + 1$ periods of $a$, then $u$ may also be a prefix of another relator $b^n \in \Rel_{r}$. In that case it is possible that $u$ has no proper prolongation in~$A$ with respect to $a$, but $u$ does have a proper prolongation in $A$ with respect to $b$.
\end{remark}


For further reference we can now state the following characterization of maximal occurrences:

\begin{remark}\label{rem:max_occurrence}
Suppose  $A = LuR\in \Can_{r - 1}$ with $u = a^ka_1$,  $a^n \in \Rel_j, 1\leq j\leq r$, $a = a_1a_2$ (where $a_2$ can be empty) and $k \in \mathbb{N} \cup \lbrace 0\rbrace$.
Then  by Remark~\ref{rem:cancellation},   $u$ does not have a proper prolongation in $A$ with respect to $a$ if and only if there are no cancellations in the words $L\cdot a\inv $ and $a_1\inv a_2\inv \cdot R$.

In particular we see that if $vu$ and $wu$ are prolongations of $u$ with respect to $a$, then $v$ is a suffix of $w$ or conversely and the word $v\cdot w\inv$ is not reduced.
\end{remark}

\begin{corollary}
\label{max_occ_prolongation} 
Let $A \in \Can_{r - 1}$ and let $u=a^ka_1, k\geq\tau+1, a^n\in \Rel_j, a=a_1a_2,$  for some $1 \leq j \leq r$. Then
there exists a unique maximal occurrence of rank~$j$ containing $u$  and it coincides with the maximal prolongation of $u$ in $A$.
\end{corollary}
\begin{proof}
This follows directly from Lemma \ref{common_part_in_rank_r} and the previous remark.
\end{proof}

For further reference we now record the following version of  Lemma~\ref{common_part_in_rank_r}. Here and in what follows we say that a word $c$ is an  \emph{overlap} of words $v, w$ if $c$ is a suffix of $v$ and prefix of $w$.

\begin{corollary}
\label{overlaps_size}
Let $A \in \Can_{r - 1}$, let $u_1$ be a maximal occurrences of rank~$r$ in $A$, let $u_2$ be an occurrence in $A$ of rank~$r$ not contained in $u_1$. Write $u_1 = a^{k}a_1$, $u_2 = b^{s}b_1$, where $a^{n}, b^n \in \Rel_r$, $a = a_1a_2$, $b = b_1b_2$, and $| a| \leq | b|$. If $c$ is the overlap of $u_1$ and $u_2$, then $| c| < \min\lbrace (\tau + 1)| a|, 2| b| \rbrace$.
\end{corollary}
\begin{proof}
By taking inverses if necessary we assume that $c$ is a suffix of $u_1$ and prefix of $u_2$. Then 
 we can write $c$ as $c = \widehat{a}^{k_1}\widehat{a}_1$, where $\widehat{a}$ is a cyclic shift of $a$ and $\widehat{a} = \widehat{a}_1\widehat{a}_2$ and $k_1\geq 0$. Then $c$ is a common prefix of $\widehat{a}^N$ and $b^N$ for some $N\in \mathbb{N}$ and  a cyclic shift $\widehat{a}$ of $a$. Since $u_1$ is a maximal occurrence and $u_2$ is not contained in $u_1$ we see that $\hat{a}\neq b$ and hence the claim follows from Lemma~\ref{common_part_in_rank_r}.
\end{proof}

If  $u = a^ka_1$ is a proper prefix of the relator $a^n \in \Rel_r, a=a_1a_2$,  we call $v = a^{-n}\cdot u=a^{k - n + 1}a_2\inv $  \emph{the complement of~$u$ with respect to the relator~$a^n$}.
Clearly, $v$ is the complement of $u$ with respect to $a^n$ if and only if $u\cdot v\inv=v\inv\cdot u=a^n$.

\begin{remark}\label{max_occurrence_replacement}
Let $A =LuR\in \Can_{r - 1}$ where $u= a^ka_1, k\geq 0, a^n\in\Rel_r,$ is a maximal occurrence of rank~$r$ in $A$ and put $v = a^{-n}\cdot u$. If $v\neq 1$, then we have $L\cdot v\cdot R=LvR$: if $k< n$ there are no cancellations in $L\cdot v\cdot R$  by Remark~\ref{rem:max_occurrence} as $u$ is a maximal occurrence and if $k> n$ there are no cancellations in $L\cdot v\cdot R$ because there are no cancellations in the initial word $LuR = La^ka_1R$. In particular, if $k> n$, then $v$ has no prolongation with respect to $a$.

Note also that if $u$ contains $\geqslant \tau + 1$ periods of the relator $a^n$ (that is, if $| u| \geqslant (\tau + 1)| a| $), then  $a^n$ is the unique relator in $\Rel_r$ with prefix $u$ by Remark~\ref{subword_of_rel_power}. So, the complement of $u$ is defined without referring to the particular relator and in this case we will simply  call $v$ \emph{the complement of~$u$}. 
\end{remark}
 
The next definition is central to our approach:
\begin{definition}[turns of rank~$r$]\label{def:turns of rank r}
Let $A =LuR\in \Can_{r - 1}$ where $u= a^ka_1, a=a_1a_2, k\geq 0, a^n\in\Rel_r,$ is a maximal occurrence of rank~$r$ in $A$.  Let $v = a^{-n}\cdot u$. 
The  transformation
\begin{equation*}
A = LuR \longmapsto  \can_{r - 1}(L vR)
\end{equation*}
is called \emph{\bf a turn of rank~$r$ in $A$}, or, more specifically, \emph{\bf the turn of $u$ in $A$}.
\end{definition}
Note that we may have $k\geq n$ 
and that the reduced form of $v$ is one of the following:
\begin{equation}
\label{turn_first_step}
v = \begin{cases}
a^{k - n}a_1 &\textit{ if } k \geq n ,\\
a^{k - n + 1}a_2\inv  &\textit{ if } k < n.
\end{cases}
\end{equation}

The following observation will be convenient:

\begin{remark}\label{rem:cyclic shift}
(i) Let $a=a_1a_2$ be a cyclically reduced word and $\widehat{a} = a_2a_1$ a cyclic shift of $a$. For any $t, k\in\mathbb{N}$ we have 
\[a^{-t}\cdot (a^ka_1)=a^{k-t}\cdot a_1=a_1\cdot (a_2a_1)^{k - t}= a_1 (a_2a_1)^k \cdot (a_2a_1)^{-t}\]
\[=a^ka_1 \cdot (a_2a_1)^{-t} =(a^ka_1)\cdot \widehat{a}^{-t}.\]
Hence if $A = LuR\in\Can_{r-1}$ where  $u = a^ka_1$  is a maximal occurrence of rank~$r$ in $A$ with complement $v=a^{-n}\cdot u$, then we see that we can move the multiplication with $a^{-n}$ to any position across $u$ by using the appropriate cyclic shift $\widehat{a}^n$ of $a^n$:

\[LvR=L\cdot a^{-n}\cdot uR= Lu'\cdot \widehat{a}^{-n}\cdot u^{\prime\prime}R \]
where $u= u'u^{\prime\prime}$ and $u^{\prime\prime}$ is a prefix of $\widehat{a}^N$ for some $N\geq 0$.

(ii) Note that this also shows that if $u$ is a not necessarily maximal occurrence of rank~$r$ in $A\in\Can_{r-1}$, $u=a^ka_1$ for some $a^n\in\Rel_r$ with $k\geq\tau+1$, then if we  multiply $u$ from the left by $a^{-n}$ (and take its canonical form),  we automatically turn the maximal prolongation of $u$ with respect to $a$. 
\end{remark}

\begin{remark}\label{rem:turn induction beginning}
Let $A = LuR\in\Can_{r-1}$ where  $u = a^ka_1$  is a maximal occurrence of rank~$r$ in $A$ with complement $v=a^{-n}\cdot u$. Then by Remark~\ref{can_preserves_coset_hyp} and Remark~\ref{rem:cyclic shift} the result $\can_{r - 1}(LvR)$ of turning $u$ in $A$ satisfies
\begin{align*}
\label{turn_characterisation}
\can_{r - 1}(LvR) =\can_{r - 1}(L\cdot a^{-n}\cdot uR)& \equiv L\cdot a^{-n}\cdot uR \\
&\equiv Lu'\cdot \widehat{a}^{-n}\cdot u^{\prime\prime}R \mod \llangle \Rel_0, \ldots, \Rel_{r - 1}\rrangle,
\end{align*}
for any decomposition $u= u'u^{\prime\prime}$, where $\hat{a}^n$ is the appropriate cyclic shift of~$a^n$.
\end{remark}

In order to describe the resulting word after a turn of rank~$r$ we first establish the following general lemma:
\begin{lemma}\label{lem:product preparation}
Suppose $a^n\in\Rel_r$, $a_2,a_3$ are (possibly empty) suffixes of $a$,
$La_2a^M$ and $R\inv a_3a^K$ are prefixes of words in $\Can_{r-1}$ and assume $a_3a^K$ is a maximal occurrence,  $\Lambda_r(a_2a^M)-\tau\geq \Lambda_r(a_3a^K)\geq 2\tau$. Then there exists a canonical triangle $(D_1, D_2, D_3)$ such that
\[\can_{r-1}(La_2a^M\cdot a^{-K}a_3\inv R)=\tilde{L}D_3 R',\tilde{L}D_1 = La_2a^M\cdot a^{-K}a_3^{-1}, R = D_2R'.\]

Furthermore, if $\Lambda_r(a_2a^M)-\tau \geq \Lambda_r(a_3a^K)$, then $\tilde{L} = Lw_0$ for a prefix $w_0$ of $a_2a^M$ with $\Lambda_r(w_0) > \Lambda_r(a_2a^M\cdot a^{-K}a_3\inv)-\tau$.
\end{lemma}

Note that by considering inverses and using the fact that $\Can_{r-1}$ and $\Rel_r$ are closed under inverses, for the case $2\tau \leq \Lambda_r(a_2a^M) \leqslant \Lambda{a_3a^K} - \tau$ we also obtain
\[
\can_{r-1}(La_2a^M\cdot a^{-K}a_3\inv R)=L'D_3 w_0R,
L = L'D_1, D_2w_0R = a_2a^M\cdot a^{-K}a_3^{-1}R\]
for a suffix $w_0$ of $a^{-K}a_3\inv$ and some canonical triangle $(D_1, D_2, D_3)$.

\begin{proof}
By Corollary~\ref{can_of_can_product_hyp} we have 
\[\can_{r-1}(La_2a^M\cdot a^{-K}a_3\inv R)=\can_{r-1}(\can_{r-1}(La_2a^M)\cdot \can_{r-1}(a^{-K}a_3\inv R)).\]

By Corollary~\ref{can_power_inverse_hyp_cor} we have
\[\can_{r-1}(La_2a^M)=La^{M-{\tau}}X \quad \mbox{ and \ }\can_{r-1}(a^{-K}a_3\inv R))=X\inv a^{-K+{\tau}}a_3\inv R.\]
\[\mbox{So}\quad
\can_{r-1}(La_2a^M\cdot a^{-K}a_3\inv R)=\can_{r-1}(La_2a^{M-\tau}X \cdot X\inv a^{-K+\tau}a_3\inv R).\]
Put $W=a_3a^{K-\tau}X$. Since $a_3a^K$ is maximal, $W\cdot W\inv$ is the maximal cancellation in this product.
By IH~\ref{can_triangle_hyp1} there is a canonical triangle $(D_1,D_2, D_3)$
such that  \[La_2a^{M-\tau}X= LwW=\tilde{L}D_1W \mbox{ and } X\inv a^{-K+\tau}a_3\inv R=W\inv R= W\inv D_2R'\] for a prefix $w$ of $a_2a^{M-K}$, prefix $\tilde{L}$ of $Lw$ and suffix $R'$ of $R$ such that
\[\can_{r-1}(La_2a^{M-K}\cdot a_3\inv R)=L'D_3R'.\]

If $\Lambda_r(a_2a^M\cdot a^{-K}a_3\inv)\geq\tau$, then $\tilde{L}=Lw_0$ for a nonempty prefix $w_0$ of $a_2a^M$ with $\Lambda_r(w_0)>\Lambda_r(a_2a^M\cdot a^{-K}a_3\inv)-\tau$.
\end{proof}

The following proposition describes the resulting word after a turn of  a maximal occurrence of rank~$r$ of $\Lambda_r$-measure $\geq \tau$.
Below is an illustration that presents both $A=LuR$ and the result $B$ of the turn of $u$ in $A$ with complement $v$. Note that  the canonical triangles $(D_1, D_2, D_3)$ and $(E_1,E_2,E_3)$ could intersect. In fact, the relative position of these two triangles on the circle corresponds to the Types 2.  and 3. in the following proposition.
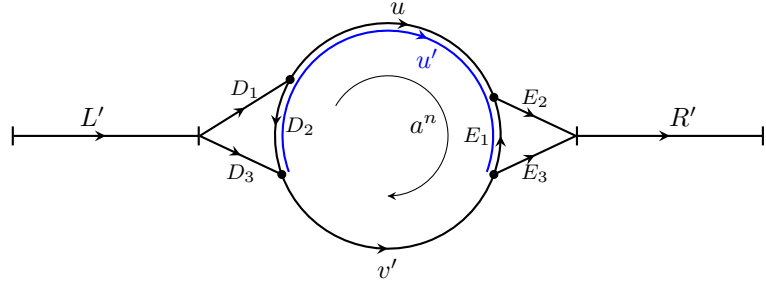
\begin{center}
\begin{tikzpicture}
\draw[thick, arrow=0.5] (-20:1.5) arc (-20:20:1.5) node[midway, left] {\footnotesize $E_1$};
\draw[thick, reversearrow=0.5] (20:1.5) arc (20:150:1.5) node[midway, above] {$u$};
\draw[thick, arrow=0.5] (150:1.5) arc (150:200:1.5) node[midway, right] {\footnotesize $D_2$};
\draw[thick, arrow=0.5] (200:1.5) arc (200:340:1.5) node[midway, below] {$v'$};
\filldraw (-20:1.5) circle (1.5pt);
\filldraw (20:1.5) circle (1.5pt);
\filldraw (150:1.5) circle (1.5pt);
\filldraw (200:1.5) circle (1.5pt);
\draw[thick, arrow=0.6, blue] (200:1.4) arc (200:-20:1.4) node[pos=0.6, below] {$u'$};
\draw[arrow=1] (150:0.8) arc (150:-90:0.8) node[pos=0.5, below, xshift=-6] {$a^n$};

\draw[thick, arrow=0.5] (180:2.5) to node[midway, above] {\footnotesize $D_1$} (150:1.5);
\draw[thick, arrow=0.5] (180:2.5) to node[midway, below] {\footnotesize $D_3$} (200:1.5);
\draw[thick, |-|, arrow=0.5] (180:5) to node[midway, above, xshift=-5] {$L'$} (180:2.5);

\draw[thick, arrow=0.5] (20:1.5) to node[midway, above] {\footnotesize $E_2$} (0:2.5);
\draw[thick, arrow=0.5] (-20:1.5) to node[midway, below] {\footnotesize $E_3$} (0:2.5);
\draw[thick, |-|, arrow=0.5] (0:2.5) to node[midway, above, xshift=5] {$R'$} (0:5);
\end{tikzpicture}
\end{center}

\begin{lemma}
\label{turn_res_complement}
Let $A=LuR \in \Can_{r - 1}$ where $u = a^ka_1$ is a maximal occurrence of rank~$r$ in $A$,  $a^n \in \Rel_r$, $k\geq \tau$, and $a = a_1a_2$ (where $a_1$ can be empty). Let $v = a^{-n}\cdot u$ and consider the turn of rank~$r$ in $A$:
\[
LuR \longmapsto LvR \longmapsto \can_{r - 1}(LvR) \mbox{ if } v\neq 1;\]
\[LuR \longmapsto L\cdot  R \longmapsto \can_{r - 1}(L\cdot R) \mbox{ if } v= 1.
\]
Put $m=\tau-1$ if $k<n$ and $m=\tau+k-n$ if $k\geq n$.

\medskip
\noindent
\rm{(i)} The result of the turn is of the form \[L'Q \tilde{R}\]
where $L'$ is a prefix of $L$, $\tilde{R}$ is a proper  suffix of $uR$ and we have one of the following 
three possibilities:

\begin{enumerate}[label=Type\,\arabic*., ref=\arabic*]
\item
\label{v_long_power}
\begin{equation*}
\can_{r - 1}(LvR) = LvR, \quad L'=L, Q=v \mbox{ is a fractional power of } a \mbox{ and }\tilde{R}=R;
\end{equation*}

\item
\label{non_empty_v_part}
\begin{align*}
&\can_{r - 1}(L v R) = L'D_3v'E_3R',
&L = L'D_1,\ R = E_2R',\ v = D_2v'E_1,
\end{align*}

where the {\bf remainder} $v'$ of $v$ is non empty, and $(D_1, D_2, D_3)$ and $(E_1, E_2, E_3)$ form canonical triangles of rank~$r - 1$. Here $v$ is a fractional power of  $a\inv$ and $ Q=D_3v'E_3$.
\item
\label{empty_v_part1}
\[
\can_{r - 1}(LvR) = L'D_3'E_3\tilde{R},\ L = L'D_1,\ Q=D_3'E_3,
\]
where  $D_3'$ is a not-empty prefix of a side of a canonical triangle of rank~$r - 1$, $D_1$ and $E_3$ are a sides of canonical triangles of rank~$r - 1$, $L'$ is a prefix of $L$ and $\tilde{R}$ is a proper suffix of $a^m a_1R$.

\item
\label{empty_v_part2}
\[
\can_{r - 1}(LvR) = L'E_3\tilde{R}, \ Q=E_3,
\]
where $E_3$ is a side of a canonical triangle of rank~$r - 1$, $L'$ is a prefix of $L$ and $\tilde{R}$ is a proper suffix of $a^m a_1R$.
\end{enumerate}

\medskip
\noindent
\rm{(ii)}  If $k\geq n+\tau$, the result is of   Type~\ref{v_long_power} and if $\tau\leq \Lambda_r(u)\leq n-2\tau$ (or equivalently, $n-\tau\geq \Lambda_r(v) \geqslant 2\tau$), the result is of  Type\,\ref{non_empty_v_part} with $\Lambda_r(v') > \Lambda_r(v) - 2\tau  \geqslant 0$.

\medskip
\noindent
\rm{(iii)} Unless $Q = D_3v'E_3$ with $\abs{v'} \geqslant (\tau + 1)\abs{a}$ or $Q = v$ with $\abs{v} \geqslant (\tau + 1)\abs{a}$, $Q$ is $(3\tau +1)$-free of rank~$r$.
\end{lemma}

\begin{proof}
First assume $k-n \geq \tau$. Since $LuR = La^ka_1R \in \Can_{r - 1}$ and $k, k-n \geqslant~\tau$, it follows from IH~\ref{can_periodic_extension_constriction_hyp} that $LvR = La^{k - n}a_1R \in \Can_{r - 1}$ as well. Hence $\can_{r - 1}(LvR) = LvR$ by IH~\ref{IH idempotent} and so the result is of  Type~\ref{v_long_power}.

Now suppose that $\tau\leq k<n+\tau$. Then $v=a^{k - n + 1}\cdot a_2\inv $ and this product is reduced if and only if $k<n$. While we would like to compute $\can_{r-1}(LvR)$ by applying Lemma~\ref{lem:product preparation} to $La^{-n}\cdot uR$, we do not know that $La^{-n}$ is a prefix of a word in $\Can_{r-1}$. Therefore we appeal to Corollary~\ref{can_of_can_product_hyp}  and first write $LvR$ as a product of suitable subwords from $\Can_{r-1}$. To this end we take~$N\geq 2\tau$ and rewrite $LvR$ as follows:
\begin{align*}
LvR = La^{k - n + 1}\cdot a_2\inv R &= (La^{N+\tau}) \cdot (a^{-N-\tau}\underbrace{\cdot a^{k - n + 1}\cdot a_2\inv }_{v} a_1\inv a^{-N-\tau})\cdot (a^{N+\tau}a_1R) \\
&= (La^{N+\tau})\cdot (a^{-2N-2\tau + k - n})\cdot (a^{N+\tau}a_1R).
\end{align*}

Since $A = LuR = La^ka_1R \in \Can_{r - 1}$ and $k \geqslant \tau$, it follows from  IH~\ref{can_periodic_extension_constriction_hyp} that also $La^{K}a_1R\in\Can_{r - 1}$ for any $K\geq \tau$. Therefore, $La^N$ and $a^{N}a_1R$ are prefix and suffix, respectively,  of a canonical word of rank~$r - 1$. Since  $\Can_{r-1}$ is closed under taking inverses, also
$a^{-2N -2\tau + k - n}$ is a subword of a word from $\Can_{r-1}$. Thus, Corollary~\ref{can_power_inverse_hyp_cor} applies to the words $La^N$, $a^{N}a_1R$, $a^{- 2N -2\tau + k - n}$ yielding
\begin{align}
\begin{split}
\label{can_components}
Z_1 &= \can_{r - 1}(La^{N+\tau}) = La^NX,\\
Z_2 &= \can_{r - 1}(a^{- 2N -2\tau + k - n}) = X\inv a^{-2N+k-n}Y\inv ,\\
Z_3 &= \can_{r - 1}(a^{N+\tau}a_1R) = Ya^{N }a_1R.
\end{split}
\end{align}
\begin{align*}
\mbox{Hence }\quad \can_{r-1}(Z_1\cdot Z_2)&= \can_{r-1}(La^NX\cdot  X\inv a^{-2N+k-n}Y\inv )\\
&= \can_{r-1}(La^N\cdot a^{-2N+k-n}Y\inv ). 
\end{align*}
By Lemma~\ref{lem:product preparation} applied to  $La^N$ and $Ya^{2N-k+n}$ (see the comment after the lemma) we find a canonical triangle $(D_1,D_2,D_3)$ such that 
\[Z=\can_{r-1}(Z_1\cdot Z_2)=L'D_3v_0Y\inv \mbox{ with } L=L'D_1,\ D_2v_0Y\inv=a^{-N-n+k}Y\inv.\]
Clearly if $v$ is a fractional power of $a^{-1}$ (i.e. if $k < n$), then $v$ is a prefix of $D_2v_0Y^{-1}$.

\begin{align*}
\mbox{Now}\quad \can_{r-1}(LvR)=\can_{r-1}(Z\cdot Z_3)&=\can_{r-1}(L'D_3v_0Y\inv\cdot Ya^{N }a_1R)\\
   &= \can_{r-1}(L'D_3v_0\cdot a^{N }a_1R).
\end{align*}
Let $v_0'$ be the prefix of $v_0$ (if any) that is not cancelled in the product $v_0\cdot a^Na_1$. Note that $v_0$ may have a proper prolongation $\tilde{v}= v_1v_0$ in $L'D_3v_0$ with respect to $a\inv$. We again apply Lemma~\ref{lem:product preparation} to $L'D_3v_0$, $a^Na_1R$ or their inverses and obtain another canonical triangle $(E_1, E_2,E_3)$ and the following cases according to the position of this triangle relative to $D_3$:

{\bf Type 2} 
\[\can_{r-1}(Z\cdot Z_3)=L'D_3v'E_3R'\]
where $L=L'D_1$, $R= E_2 R'$, and $v=D_2v'E_1$.  So in particular $\Lambda_r(v')>\Lambda_r(v)-2\tau$, and this happens exactly if $v_0'$ is not contained inside $E_1$:

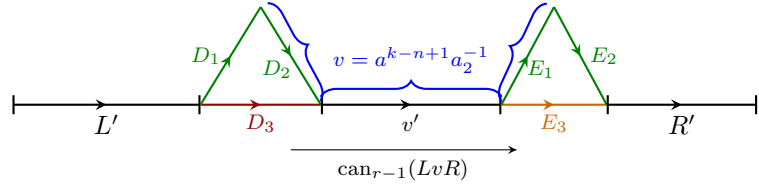
\begin{center}
\label{non_empty_v_part_picture}
\begin{tikzpicture}
\draw[|-|, thick, arrow=0.5] (-2, 0) to node[midway, below] {$L'$} (0.5, 0);

\draw[thick, arrow=0.5, black!40!red] (0.5, 0) to node[midway, below] {\footnotesize $D_3$} (2.1, 0);
\draw[thick, black!50!green, arrow=0.5] (1.3, 1.3) to node[midway, left, yshift=-5, xshift=4] {\footnotesize $D_2$} (2.1, 0);
\draw[thick, black!50!green, arrow=0.5] (0.5, 0) to node[midway, left] {\footnotesize $D_1$} (1.3, 1.3);

\draw[|-|, thick, arrow=0.5] (2.1, 0) to node[midway, below] {\footnotesize $v'$} (4.5, 0);

\draw[thick, black!20!orange, arrow=0.5] (4.5, 0) to node[midway, below] {\footnotesize $E_3$} (5.9, 0);
\draw[thick, black!50!green, arrow=0.5] (5.2, 1.3) to node[midway, right] {\footnotesize $E_2$} (5.9, 0);
\draw[thick, black!50!green, arrow=0.5] (4.5, 0) to node[midway, right, yshift=-5, xshift=-4] {\footnotesize $E_1$} (5.2, 1.3);

\draw[|-|, thick, arrow=0.5] (5.9, 0) to node[midway, below] {$R'$} (7.9, 0);

\draw[blue, thick, decorate, decoration={brace, amplitude=6pt, raise=3pt}] (2.1, 0) to (4.5, 0);
\draw[blue, thick, decorate, decoration={brace, amplitude=6pt, raise=3pt}] (4.6, 0) to (5.2, 1.3);
\draw[blue, thick, decorate, decoration={brace, amplitude=6pt, raise=3pt}] (1.3, 1.3, 0) to (2, 0);
\node[blue] at (3.3, 0.6) {\footnotesize $v = a^{k - n + 1}a_2\inv $};

\draw[black, arrow = 1] (-1.5 + 3.2, -0.6)--(1.5 + 3.2, -0.6) node[midway, below] {\footnotesize $\can_{r-1}(LvR)$};
\end{tikzpicture}
\end{center}

{\bf Types 3 and 4} 
If $v_0'$ is contained in $E_1$ (in particular if  $v_0$ cancels completely), then 
\[\can_{r-1}(Z\cdot Z_3)=L'D_3'E_3\tilde{R} \textit{ or } \can_{r-1}(Z\cdot Z_3)=L''E_3\tilde{R},\]
where $D_3'$ is a non-empty prefix of $D_3$, $L''$ is a prefix of $L'$ and $\tilde{R}$ is a suffix of $a^m a_1R$. Notice that the suffix of $a^Na_1R$ remaining after cancellation with $L'D_3v_0$ is a proper suffix of $a^ma_1R$, so $\widetilde{R}$ is a proper suffix of $a^ma_1R$. If $E_1$ is properly contained in $D_3v_0^{\prime}$, the we obtain the first formula that gives Type~3. Otherwise we obtain the second formula that gives Type~4.
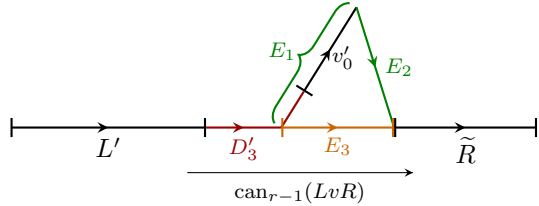
\begin{center}
\begin{tikzpicture}
\draw[|-|, thick, arrow=0.5] (0, 0) to node[midway, below] {$L' $} (2.6, 0);
\draw[thick, arrow=0.5, black!40!red] (2.6, 0) to node[midway, below] {\footnotesize $D_3'$} (3.6, 0);
\path (3.6, 0) to coordinate[pos=0.3] (A1) (4.6, 1.6);
\draw[thick, black!40!red] (3.6, 0) to (A1);
\draw[thick, |-, arrow=0.5] (A1) to node[midway, right, yshift=-3, xshift=-3] {\footnotesize $v_0'$} (4.6, 1.6);
\draw[thick, black!50!green, arrow=0.5] (4.6, 1.6) to node[midway, right] {\footnotesize $E_2$} (5.1, 0);

\draw[black!50!green, thick, decorate, decoration={brace, amplitude=6pt, raise=3pt}] (3.6, 0) to node[left, midway, xshift=-5, yshift=6] {\footnotesize $E_1$} (4.6, 1.6);

\draw[|-|, thick, black!20!orange, arrow=0.5] (3.6, 0) to node[midway, below] {\footnotesize $E_3$} (5.1, 0);
\draw[|-|, thick, arrow=0.5] (5.1, 0) to node[midway, below] {$\tilde{R}$} (7, 0);

\draw[black, arrow = 1] (-1.5 + 3.85, -0.6)--(1.5 + 3.85, -0.6) node[midway, below] {\footnotesize $\can_{r-1}(LvR)$};
\end{tikzpicture}
\end{center}

Now we prove the last part: if $Q = D_3'E_3$ or $Q = E_3$, then $Q$ is $2\tau$-free of rank~$r$, because $D_3'$ and $E_3$ are $\tau$-free of rank~$r$.

Let $Q = D_3v'E_3$ and $\abs{v'} < (\tau + 1)\abs{a}$. Assume that $Q$ contains $b^{3\tau+1}$, where $b^n\in \Rel_r$. Since $D_3$ and $E_3$ are $\tau$-free of rank~$r$, we obtain that $v'$ contains $b^{\tau + 1}$. Hence, it follows from Lemma~\ref{common_part_in_rank_r} that $a\inv $ is a cyclic shift of $b$. Therefore, $\abs{v} \geqslant (\tau + 1)\abs{b} = (\tau + 1)\abs{a}$, a contradiction.

Let $Q = v$ and $\abs{v} < (\tau + 1)\abs{a}$. Assume that $Q$ contains $b^{3\tau+1}$, where $b\in \Rel_r$. Then it follows from Lemma~\ref{common_part_in_rank_r} that $a$ is a cyclic shift of $b$. Therefore, $|v| \geqslant (3\tau+1) |b| > (\tau + 1)|a|$, a contradiction.

\end{proof}

By considering inverses and using IH~\ref{IH can inverse} we also obtain the following ``left'' version of Type~\ref{empty_v_part1}  in Lemma~\ref{turn_res_complement} (instead of the current ``right'' version). It is important to note that while the description of the canonical form may differ, it is in fact uniquely defined and therefore, these two versions agree.

\begin{remark}\label{rem:left version} In the situation of Lemma~\ref{turn_res_complement} (and with the same notation)   we obtain the following ``left'' description of $\can_{r - 1}(LvR)$ for Types~\ref{empty_v_part1} and~\ref{empty_v_part2}:  \\
\noindent
\emph{
Type~\ref{empty_v_part1}' and~\ref{empty_v_part2}'.\\
\[
\can_{r - 1}(LvR) = \tilde{L}F_3G_3'R',
\]
where $R'$ is a suffix of $R$, $F_3$ is a side of a canonical triangle of rank~$r - 1$, $G_3'$ is a (possibly empty) suffix of a side of a canonical triangle of rank~$r - 1$, and $\tilde{L}$ is a proper prefix of $La^ma_1$.}\\

\end{remark}

\begin{convention}\label{def:Type 2}
If $A=LuR\in\Can_{r-1}$ for some maximal occurrence $u$ of rank~$r$ with $\tau\leq\Lambda_r(u)<n$ we say that the turn $A\mapsto \can_{r-1}(LvR)=B$ is of Type 2 provided
$B=L'D_3v'E_3R', L = L'D_1, R = D_2R', v = D_2v'E_1$ as in Type~\ref{non_empty_v_part} of Lemma~\ref{turn_res_complement}.
\end{convention}

\begin{corollary}\label{cor:prefix turn}
Let $A_1=Lu_1R_1,\quad A_2=Lu_2 R_2 \in \Can_{r - 1}$ where $u_1=a^ka_1, u_2=a^ma_2$ are maximal occurrences of rank~$r$, $u_1$ is a prefix of $u_2$ and $\tau\leq \Lambda_r(u_1)\leq\Lambda_r(u_2)\leq n-2\tau$. Let $v_i, i=1,2,$ be the complement of $u_i$.
Then there is a canonical triangle $(D_1, D_2, D_3)$ such that the result $B_i$ of turning $u_i$ in $A_i, i=1, 2$, is of the form
\[B_1=L'D_3v_1'E_3R_1' \mbox{ \quad and } B_2=L'D_3v_2'F_3R_2'\]
where $v_2'$  and $v_1'$ have a common a prefix of $\Lambda_r$-measure $> n-\Lambda_r(u_2)-2\tau$, $R_i'$ is a suffix of $R_i, i= 1, 2,$ and $E_3, F_3$ are sides of respective canonical triangles of rank $r-1$.
\end{corollary}
\begin{proof}
Consider the decomposition of $Lv_1R_1$ and $Lv_2R_2$ into three factors $Z_1, Z_2, Z_3$ as in the proof of Lemma~\ref{turn_res_complement}. Then the factor $Z_1=La^NX$ is identical in both cases, the factor $Z_2$ is of the form $X\inv a^{-2N+k-n}Y\inv$ and $X\inv a^{-2N+m-n}Y\inv$, so differs only in the exponent of $a$ by $m-k$, and the third factor is of the form $Ya^Na_iR_i$, $i=1, 2$.
By Corollary~\ref{can_power_context_hyp_cor1} we see that in either case the product $Z=\can_{r-1}(Z_1\cdot Z_2)$ is of the form $L'D_3v_{0,1}Y\inv$ and $L'D_3v_{0,2}Y\inv$, respectively, where  $\Lambda_r(v_{0,1})=\Lambda_r(v_{0,2})+\Lambda_r(u_2)-\Lambda_r(u_1)$.
Thus, looking at the proof of Lemma \ref{turn_res_complement} we see that after multiplying either of these results with $Z_3$ the product will have a prefix of the form $L'D_3v_i'$ for some prefix $v_i'$ of $v_{0,1}$ of $\Lambda_r$-measure $>n - \Lambda_r(u_2) - 2\tau$.
\end{proof}

\begin{corollary}
\label{cor:left-most after turn}
Let $A=LuR \in \Can_{r - 1}$ where $u$ is a maximal occurrence of rank~$r$  with   
  $\tau\leq \Lambda_r(u)\leq n-2\tau$ and let $B=L'D_3v'E_3R'$ be the result of turning $u$ in $A$.
Let $w$ be an  occurrence of rank~$r$ in $R$ with $\Lambda_r(w)\geq \tau$. Then $R'$ contains a non-empty suffix $w'$ of $w$ with $\Lambda_r(w')> \Lambda_r(w)-\tau$.
\end{corollary}
\begin{proof}
By Lemma~\ref{turn_res_complement} $B$ is of the given form. 
By the description of  Type~\ref{non_empty_v_part}, we have $R=E_2R'$ where $E_2$ is $\tau$-free of rank~$r$. So $w$ cannot be entirely contained in $E_2$, so  $\Lambda_r(w')> \Lambda_r(w)-\tau$.
\end{proof}

\begin{corollary}
\label{cor:turn protection aux}
Let $A_i = LuMb^\tau R_i\in \Can_{r - 1}, i= 1, 2,$ where $u$ is a maximal occurrence of rank~$r$  with   
 $\tau\leq \Lambda_r(u)\leq n-2\tau$, $b^n \in \Rel_r$, and assume that the result $B_1$ of turning $u$ in $A_1$ is of the form
 \[B_1 = L'D_3v'E_3M'b^\tau R_1,\]
where $M^{\prime}$ is a suffix of $M$ and $(D_1, D_2, D_3), (E_1, E_2, E_3)$ are canonical triangles of rank~$r - 1$. Then the result
$B_2$ of turning $u$ in $A_2$ is of the form
\[B_2 = L'D_3v'E_3M'b^\tau R_2\]
with the same canonical triangles.
 \end{corollary}
 \begin{proof}
 This follows directly from  Corollary~\ref{cor:left-most after turn}, IH~\ref{can_periodic_extension_constriction_hyp}, IH~\ref{IH idempotent}, and IH~\ref{IH can of equivalent words} (see also the proofs of Corollaries~\ref{can_power_context_hyp_cor1} and \ref{cor:prefix turn}).
 \end{proof}

The following statement is a  useful particular case of Corollary~\ref{cor:turn protection aux} with $M = M_1a^\tau M_2$:
\begin{corollary}
\label{cor:turn protection}
Let $A_i=LuM_1a^\tau M_2b^\tau R_i\in \Can_{r - 1}, i= 1, 2,$ where $u$ is a maximal occurrence of rank~$r$  with   
  $\tau\leq \Lambda_r(u)\leq n-2\tau$, $a^n, b^n \in \Rel_r$.
  Then there are canonical triangles $(D_1, D_2, D_3), (E_1, E_2, E_3)$ such that the result $B_i, i=1, 2$, of turning $u$ in $A_i$ is of the form 
  \[B_i=L'D_3v'E_3 M_2'b^\tau R_i\]
  where $M_2'$ contains a non-empty suffix of $a^\tau$, i.e. the canonical triangles do not depend on $R_i$.
 \end{corollary}

\subsection{Inverse turns}\label{sec:inverse turns}

 From now on we will fix the following notational conventions:
\begin{convention}
For $A=LuR\in\Can_{r-1}$ where $u$ is a maximal occurrence of rank~$r$ with $\tau\leq \Lambda_r(u)\leq n-2\tau$ (i.e. the turn of $u$ is of Type 2) we use the following conventions:
\begin{enumerate}
\item $v$ denotes the complement of $u$;
\item $v'$ denotes the remainder of $v$ after turning $u$ in $A$; 
\item $B=\can_{r-1}(LvR)=L'D_1v'E_1R'$ is the result of turning $u$ in $A$ for canonical triangles $(D_1, D_2, D_3)$ and $(E_1, E_2, E_3)$ where $L=L'D_1, R=E_2R'$ and $v=D_2v'E_1$;
\item if $\Lambda_r(v')\geq\tau+1$, then $\hat{v}$ is the maximal prolongation of $v'$ in $B$ (and coincides with the maximal occurrence of rank~$r$ in $B$ containing $v'$).
\end{enumerate}
\end{convention}

We next prove that a turn of Type~\ref{non_empty_v_part}  of a maximal occurrence $u$ with  complement $v$ has a natural inverse turn, namely the turn of the maximal  prolongation $\hat{v}$ of $v'$ in $B$ (provided $\hat{v}$ is a maximal occurrence).

\begin{lemma}
\label{inverse_turn_formal_existence}
Let $A=LuR\in \Can_{r - 1}$, where $u$ is a maximal occurrence of rank~$r$ in $A$ with  $\tau\leq \Lambda_r(u) <n-(3\tau+1)$ and let  $B=L'D_3v'E_3R'$ be the result of  turning   $u$ in $A$. Then the result of turning $\hat{v}$ in $B$ is equal to $A$.
\[\can_{r-1}(L'D_3\cdot\tilde{a}^n\cdot v'E_3R')=A\]
where $u$ is a prefix of $a^n\in\Rel_r$ and $v'$ a prefix of some cyclic shift $\tilde{a}^{-n}$ of $a^{-n}$.
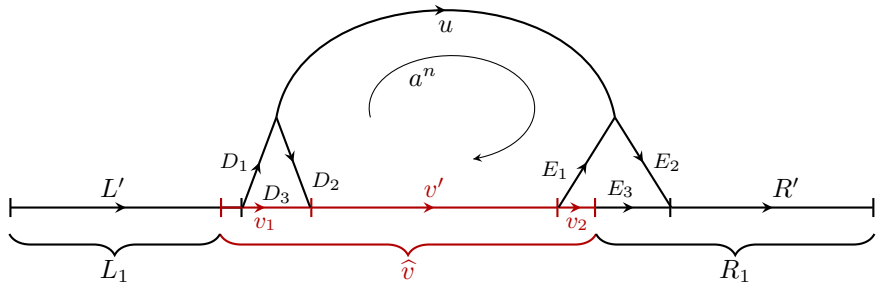
\begin{center}
\begin{tikzpicture}
\draw[thick, |-|, arrow = 0.5] (-0.5, 0) to node[midway, above, xshift=-5] {$L'$} (2.6, 0);
\path (2.6, 0) to  node[midway, above, yshift=-2] {\footnotesize $D_3$}  (3.5, 0);
\draw[thick, |-, arrow = 0.5, black!30!red] (2.3, 0) to node[midway, below] {\small $v_1$}  (3.5, 0);

\draw[thick, arrow = 0.5, black!30!red, |-|] (3.5, 0) to  node[midway, above] {$v'$} (6.3 + 0.5, 0);
\draw[thick, -|, black!30!red, arrow = 0.6] (6.3 + 0.5, 0) to  node[midway, below] {\small $v_2$}  (6.8 + 0.5, 0);
\draw[thick, -|, arrow = 0.5] (6.8 + 0.5, 0) to  node[midway, above, xshift=-5] {\footnotesize $E_3$}  (7.8 + 0.5, 0);

\draw[thick, -|, arrow = 0.5] (7.8 + 0.5, 0) to node[midway, above, , xshift=5] {$R'$} (11, 0);

\draw[thick, arrow = 0.5] (2.6, 0) to node[midway, left] { \footnotesize $D_1$} (3.05, 1.2);
\draw[thick, arrow = 0.5] (3.05, 1.2) to node[midway, right, xshift=3, yshift=-7] { \footnotesize $D_2$} (3.5, 0);

\draw[thick, arrow = 0.5] (6.3 + 0.5, 0) to node[midway, left, xshift=-3, yshift=-3] { \footnotesize $E_1$} (7.05 + 0.5, 1.2);

\draw[thick, arrow = 0.5] (7.05 + 0.5, 1.2) to node[midway, right] { \footnotesize $E_2$} (7.8 + 0.5, 0);

\draw[thick, arrow=0.5] (3.05, 1.2) to[bend left, in=100, out=80, looseness=1.1] node[midway, below] {$u$} (7.05+ 0.5, 1.2);

\draw[arrow=1] (4.3, 1.2) arc (190:-75: 1.1 and 0.7) node[pos=0.3, below] {$a^n$};

\draw[thick, decorate, decoration={brace, amplitude=8pt, raise=10pt, mirror}, black!30!red] (2.3, 0) to node[below, midway, yshift=-16] {$\widehat{v}$} (7.3, 0);
\draw[thick, decorate, decoration={brace, amplitude=8pt, raise=10pt, mirror}] (-0.5, 0) to node[below, midway, yshift=-16] {$L_1$} (2.3, 0);
\draw[thick, decorate, decoration={brace, amplitude=8pt, raise=10pt, mirror}] (7.3, 0) to node[below, midway, yshift=-16] {$R_1$} (11, 0);
\end{tikzpicture}
\end{center}
\end{lemma}
Thus, if the maximal prolongation $\hat{v}$ of $v'$ with respect to $a\inv$ is a maximal occurrence  of rank~$r$ in $B$, then the turn of $\hat{v}$ in $B$ is defined and  by Remark~\ref{rem:cyclic shift}, (ii), the result of turning $\hat{v}$ in $B$ is equal to $A$.
\begin{proof}
Suppose that $u$ is a prefix of $a^n\in\Rel_r$ and so $v'$ a prefix of some cyclic shift $\tilde{a}^{-n}$ of $a^{-n}$. 

We see that $\tilde{a}^n\cdot v'=D_2\inv uE_1\inv$ and by the properties of canonical triangles we have $D_1uE_2\equiv D_3\cdot D_2\inv uE_1\inv\cdot E_3\mod \llangle \Rel_0,\ldots, \Rel_{r-1}\rrangle$. Hence 
\begin{align*}
A&=L'(D_1 uE_2)R'\equiv L'(D_3\cdot D_2\inv u E_1\inv \cdot E_3) R' \\
&=L'D_3\cdot \tilde{a}^n\cdot v' E_3R'\mod  \llangle \Rel_0,\ldots, \Rel_{r-1}\rrangle 
\end{align*}
Since  $A\in\Can_{r-1}$ by assumption, the claim  now follows from IH~\ref{IH can of equivalent words} and IH~\ref{IH idempotent}.
\end{proof}

\begin{remark}\label{rem:not all turns have inverse}
Note that turns of Type 3 do not have inverses. Furthermore, if the remainder $v'$ after a turn of Type 2 has $\Lambda_r$-measure $<\tau+1$, then the maximal prolongation $\hat{v}$ of $v'$ with respect to $a\inv$ need not be a maximal occurrence, so again the inverse turn need not exist.
\end{remark}

In turns of Type 2, the $\Lambda_r$-measure of the maximal prolongation of $v'$   in either direction is bounded by $\tau$:

\begin{lemma}
\label{v_prolongation}
Let $A=LuR \in \Can_{r - 1}$, where $u$ is a maximal occurrence of rank~$r$ in $A$ with $\tau\leq \Lambda_r(u)<n$. Let $u$ be a prefix of $a^n\in\Rel_r$ and put $v = a^{-n}\cdot u$. Assume that the result $B$ of turning $u$ in $A$ is of Type 2 and write $B=\can_{r - 1}(LvR) = L'D_3v'E_3R'$.
 Let $\hat{v}=  v_1v'v_2$ be the maximal prolongation of $v'$ in $L'D_3v'E_3R'$ with respect to $a\inv$. Then $|v_1|\leq\max\{|D_2|, |D_3|\}$ and $|v_2|\leq\max\{|E_1|, |E_3|\}$. Thus, 
 $\Lambda_r(v_i)<\tau, i=1, 2$ and \[ n-\Lambda_r(u)-2\tau<\Lambda_r(v')\leq \Lambda_r(\hat{v})<\Lambda_r(v')+2\tau\leq n-\Lambda_r(u)+2\tau.\]
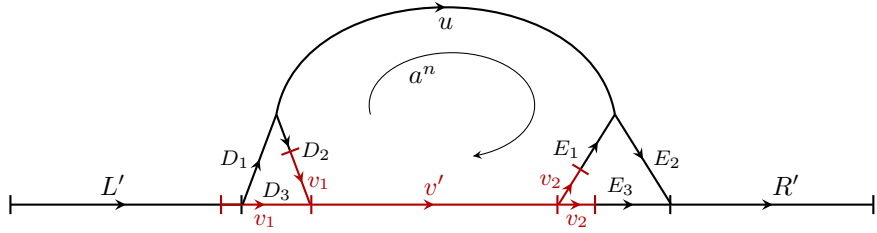
\begin{center}
\begin{tikzpicture}
\draw[thick, |-|, arrow = 0.5] (-0.5, 0) to node[midway, above, xshift=-5] {$L'$} (2.6, 0);
\path (2.6, 0) to  node[midway, above, yshift=-2] {\footnotesize $D_3$}  (3.5, 0);
\draw[thick, |-, arrow = 0.5, black!30!red] (2.3, 0) to node[midway, below] {\small $v_1$}  (3.5, 0);

\draw[thick, arrow = 0.5, black!30!red, |-|] (3.5, 0) to  node[midway, above] {$v'$} (6.3 + 0.5, 0);
\draw[thick, -|, black!30!red, arrow = 0.6] (6.3 + 0.5, 0) to  node[midway, below] {\small $v_2$}  (6.8 + 0.5, 0);
\draw[thick, -|, arrow = 0.5] (6.8 + 0.5, 0) to  node[midway, above, xshift=-5] {\footnotesize $E_3$}  (7.8 + 0.5, 0);

\draw[thick, -|, arrow = 0.5] (7.8 + 0.5, 0) to node[midway, above, , xshift=5] {$R'$} (11, 0);

\draw[thick, arrow = 0.5] (2.6, 0) to node[midway, left] { \footnotesize $D_1$} (3.05, 1.2);
\path (3.05, 1.2) to coordinate[pos=0.4] (A) (3.5, 0);
\draw[thick, arrow = 0.9] (3.05, 1.2) to node[midway, right, xshift=3, yshift=-7] { \footnotesize $D_2$} (A);
\draw[thick, arrow = 0.6, |-, black!30!red] (A) to node[midway, right, xshift=-1, yshift=-2] {\small $v_1$} (3.5, 0);

\path (6.3 + 0.5, 0) to coordinate[pos=0.4] (B) (7.05 + 0.5, 1.2);
\draw[thick, arrow = 0.6, -|, black!30!red] (6.3 + 0.5, 0) to node[midway, left, yshift=2, xshift=1] {\small $v_2$} (B);
\draw[thick, arrow = 0.5] (B) to node[midway, left, xshift=-3, yshift=-3] { \footnotesize $E_1$} (7.05 + 0.5, 1.2);

\draw[thick, arrow = 0.5] (7.05 + 0.5, 1.2) to node[midway, right] { \footnotesize $E_2$} (7.8 + 0.5, 0);

\draw[thick, arrow=0.5] (3.05, 1.2) to[bend left, in=100, out=80, looseness=1.1] node[midway, below] {$u$} (7.05+ 0.5, 1.2);

\draw[arrow=1] (4.3, 1.2) arc (190:-75: 1.1 and 0.7) node[pos=0.3, below] {$a^n$};
\end{tikzpicture}
\end{center}
\end{lemma}
\begin{proof}
Assume towards a contradiction that $|v_1| >|D_2|, |D_3|$.
Let  $zD_2$ denote the maximal common suffix of $u\inv D_2$ and $v_1$. If $z\neq u\inv$, then $zD_2=v_1$. Otherwise $\Lambda_r(v_1)\geq\tau$ since $\Lambda_r(u)\geq\tau$. Since $D_3$ is $\tau$-free of rank~$r$, we see that $|zD_2|>|D_3|$  in either case. Thus, $L'D_3=L''zD_2=L''z'D_3$ and hence $z'D_1z\inv$ is a subword of $A$. However,  since $D_3\cdot D_2\inv\equiv D_1\mod\llangle \Rel_0\cup\ldots\cup \Rel_{r-1}\rrangle$   by the definition of a canonical triangle,  we have $z'D_1z\inv\equiv z'D_3\cdot D_2\inv z\inv \equiv 1\mod\llangle \Rel_0\cup\ldots\cup \Rel_{r-1}\rrangle$, contradicting IH~\ref{can_non_trivial_subwords_hyp}.
\end{proof}

For further reference we note the following immediate consequence of the previous lemma:
\begin{corollary}\label{cor:one side always certified}
Let $A=LuR\in\Can_{r-1}$ where $u=a^ka_1$ is a maximal occurrence of rank~$r$ with $\tau\leq \Lambda_r(u)\leq n-(3\tau+1)$ and let $B=L'D_3v'E_3R'$ be the result of turning $u$. Then the maximal occurrence $\widehat{v}$ of rank~$r$ in $B$ containing $v'$ is the maximal prolongation of $v'$ with respect to $a\inv$ by Corollary~\ref{max_occ_prolongation}  and  
\[\Lambda_r(\widehat{v})< \Lambda_r(v')+2\tau\leq \Lambda_r(v)+2\tau = (n-\Lambda_r(u)) +2\tau.\]
Hence at least one of $u$ and $\hat{v}$ has $\Lambda_r$-measure $< \frac{n}{2}+\tau$, at least one of $u$ and $\hat{v}$ has $\Lambda_r$-measure $> \frac{n}{2}-\tau$, and
the turn $B\mapsto A$ of turning $\hat{v}$ in $B$ is inverse to the turn $A\mapsto B$.
\end{corollary}
\begin{proof}
 Since $v'$ has $\Lambda_r$-measure $\geq\tau+1$ by Lemma~\ref{turn_res_complement} (ii), the maximal prolongation $\hat{v}$ in $B$ with respect to $a^{-1}$ is a maximal occurrence of rank~$r$ in $B$ and its turn is inverse to the turn of $u$ in $A$ by Lemma~\ref{inverse_turn_formal_existence}. The bound on $\Lambda_r(\hat{v})$ is given in Lemma~\ref{v_prolongation}. Furthermore if $\Lambda_r(u)\geq\frac{n}{2}+\tau$, then $\Lambda_r(\widehat{v})<n-(\frac{n}{2}+\tau)+2\tau=\frac{n}{2}+\tau$.
\end{proof}

For convenience we will say that a turn of an occurrence $u$ is a turn of $\Lambda_r$-measure $\Lambda_r(u)$.

\subsection{Influence of turns on other maximal occurrences}
\label{small_interaction_turns_section}

In this subsection we describe the effect that a turn of a maximal occurrence has on other maximal occurrences in the original word. We will use the following conventions:

\begin{convention}
We will say that an occurrence $u$ in $A$ is \emph{to the left} (right) of an occurrence $w$ in $A$ if the starting point of $u$ is left (right, resp.) of the starting point of $w$ and between occurrences $w$ and $w'$ if the starting point is.
\end{convention}

\begin{convention}\label{simpler notation for occurrences}
Let $A\in\Can_{r-1}$ and assume that $u_1,\ldots, u_t$ is a  sequence of maximal occurrences of rank~$r$  in $A$. 
We use the notation $A=L\ul u_1\ldots u_t\ur R$ (thereby, in slight abuse of notation, ignoring overlaps between the occurrences or subwords separating them) where $L, R$ are prefix and suffix of $A$ which may have overlaps with $u_1, u_t$, respectively, of $\Lambda_r$-measure $<\tau +1$.  

If the word is clear from the context we may also ignore the prefix and suffix and simply write $A=\ul u_1\ldots u_t\ur$, especially in power words.
\end{convention}

\begin{convention}\label{rem:turn remainder}
Let $A\in\Can_{r-1}$ and let $u_1, u_2$ be maximal occurrences of rank~$r$. Let $B_1$ be the result of turning $u_1$ in $A$. Clearly, when turning $u_1$ the occurrence $u_2$  might be truncated to a subword $u_2''$ or even be canceled completely. However, if $\Lambda_r(u_2'')\geq\tau+1$, the maximal prolongation $\tilde{u_2}$ of $u_2''$ in $B_1$ is uniquely defined and coincides with the maximal occurrence containing $u_2''$ by Corollary~\ref{max_occ_prolongation}. We call  $\tilde{u_2}$ \emph{the occurrence corresponding to} $u_2$. For ease of notation we may then also write $\Lambda_r(u_2, B_1)$ to refer to the $\Lambda_r$-measure of $\tilde{u_2}$ in $B_1$.
\end{convention}

Turns of occurrences and multiplication of canonical words introduce perturbations on the boundaries of these operations that are captured  by the introduction of canonical triangles. Since the sides of these triangles are $\tau$-free in the corresponding ranks, an occurrences $u$  measure $\geq\tau$ in the corresponding rank absorbs the effect of the canonincal triangle and protect the remaining word from further perturbation. In other words, we will see that if $A=LuR$ for a maximal occurrence $u$ of $\Lambda_r$-measure $\geq\tau$, then a turn of rank~$r$ of Type~2 inside $L$ will have no effect on $R$ and vice versa. Therefore we introduce the following terminology:

\begin{definition}
\label{isolated_occurrences}
Let $A  =Lu_1 W u_2 R$ be a reduced word  and $u_1$ and $u_2$  maximal occurrences of rank~$r$. We say that $u_1, u_2$  are  \emph{isolated} in $A$ if  $W$ contains an  occurrence $a^{\tau}$ and  \emph{strongly isolated from each other} if $W$ contains a subword of the form $a^{\tau}M_1b^{\tau}M_2c^{\tau}$ with $a^n, b^n, c^n \in \Rel_r$ (where $M_1, M_2$ may be empty) and in this case we call $W$ a strong isolation word   (in rank~$r$).
We say that $u_1$ and $u_2$ are \emph{close neighbours} in $A$ if they are not isolated from each other.

Furthermore, we say that $u_1$ and $u_2$ are \emph{essentially non-isolated} if there are $f_1\in\{u_1, v_1\}$ and $f_2\in\{u_2, v_2\}$ such that turning $f_i$ in $W  =L\ul f_1 f_2\ur R'$ does not leave $f_j$ invariant for $\{i, j\}=\{1, 2\}$. 

We say that a word $W$ is a strong separation word  (in rank~$r$) from the right if in any word $A  =Lu WR\in\Can_{r-1}$ the maximal occurrence $u$  of rank~$r$ is   strongly isolated from any maximal occurrence of rank~$r$ in $A$ which has overlap with $R$ (and similarly for the left).

\end{definition}

\begin{examples}
\label{isolation_words_examples}
Words of the following form are strong separation words from the right:
\begin{itemize}

\item  If $A  =Lu_1 W u_2 R$ is a reduced word  such that $u_1$ and $u_2$ are essentially non-isolated maximal occurrences with $\tau \leq \Lambda_r(u_i) \leq n - (3\tau + 1)$, then $W$ does not contain a subword of the form $a^{\tau}M_1b^{\tau}M_2c^{\tau}$ with $a^n, b^n, c^n\in\Rel_r$ by Lemmas~\ref{turn_res_complement}, \ref{inverse_turn_formal_existence} and \ref{v_prolongation}.

\item
$W=a_0^{\tau}M_1a_1^{\tau}M_2a_2^{\tau}M_3a_3^{\tau + 1}M_4$, where $a_0^n, a_1^n, a_2^n, a_3^n \in \Rel_r$, $M_1, M_2, M_3$ can be empty, $M_4$ is not empty and $a_3^{\tau + 1}$ cannot be prolonged to the right.
\item
$W=a_0^{\tau}M_1a_1^{\tau}M_2a_2^{\tau}M_3a_3^{\tau + 1}$, where $a_0^n, a_1^n, a_2^n, a_3^n \in \Rel_r$, $M_1, M_2, M_3$ may be empty,  $a_3^{\tau + 1}$ cannot be prolonged to the left.
\item
$W=a_0^{\tau}M_0a_1^{\tau}M_1a_2^{\tau}M_2a_3^{\tau}M_3a_3^{\tau}$, where $a_0^n, a_1^n, a_2^n, a_3^n \in \Rel_r$, $M_0, M_1, M_2$ can be empty, and $M_3a_3^{\tau}$ is a primitive word (in particular, $M_3$ is not empty).
\end{itemize}
\end{examples}
\begin{proof}
Clearly every $W$ is a strong isolation word.
We have to show that if  $W=W_1y$ for some fractional power $y=b^kb_1$ of rank~$r$, then $W_1$ is still a strong isolation word.
For the first two cases this  follows directly from Lemma~\ref{overlaps_size}.
For the third case, this is immediate if $|b|<|a_3|$ from Corollary~\ref{overlaps_size}. If $|b|\geq|a_3|$, then $|b|<\tau|a_3|$ by Corollary~\ref{rem:no nesting in same rank}. Hence  by Lemma~\ref{common_part_of_powers}  comparing the suffixes of $b^kb_1$ and $a_3^{\tau}M_3a_3^{\tau}$ we see that $|y|<|M_3a_3^{\tau}|+|b|<|a_3^{\tau}M_2a_3^{\tau}|$.
\end{proof}

If $u_1$ and $u_2$ are isolated from each other in $A$, then $u_2$ is not affected from turning $u_1$ and vice versa:

\begin{lemma}
\label{isolated_stays_after_turn}
Let $A =Lu_1 M u_2 R\in \Can_{r - 1}$, where $u_1$, $u_2$ are maximal occurrences of rank~$r$  isolated from each other in $A$ and  $\tau \leq \FracM_r(u_1)  \leq n - 2\tau$. Let $B_1$ denote the result of turning $u_1$. Then 
\[B_1=L'D_3v_1'E_3M'u_2R\quad\mbox{ for some non-empty suffix  } M' \mbox{ of } M.\] 
In particular $\tilde{u_2}=u_2$ (as words occurring in $B_1$ and $A$, respectively, see Convention~\ref{rem:turn remainder}).

If $u_1, u_2$ are strongly isolated, then  $\hat{v_1}$ (if it is defined) is isolated from $\tilde{u_2}$ in $B_1$
\end{lemma}
\begin{proof}
Since $M$ contains an occurrence $w$ of rank~$r$ with $\Lambda_r(w)\geq \tau$, both claims follows from Corollary~\ref{cor:left-most after turn} and Lemma~\ref{v_prolongation}.
\end{proof}

In order to consider the influence of a turn on a close neighbour we first note the following:
\begin{lemma}
\label{overlap_type_occurrence}
Consider $Dv'$, where $D$  is $\tau$-free of rank~$r$ and $v'$ is a fractional power or rank~$r$. If $z$ is a maximal occurrence of rank~$r$ in $Dv'$ not containing $v'$, then $\Lambda_r(z)< 2\tau + 1$.
\end{lemma}
\begin{proof}
Write $z=z_0z_1$ where $z_0$ is a suffix of $D$ and $z_1$ a prefix of $v'$.
 Since $D$ is $\tau$-free  of rank~$r$, we have $\Lambda_r(z_0)<\tau$ and by Lemma~\ref{common_part_in_rank_r} we have $\Lambda_r(z_1)<\tau+1$. 
\end{proof}

\begin{lemma}
\label{influence_on_neighbours_complement}
Let $A  = LuR\in \Can_{r - 1}$ where $u, z$ are distinct maximal occurrences of rank~$r$ in $A$ with $\Lambda_r(u)\geq \tau+1$ and  $\Lambda_r(z)\geq 3\tau+2$. Let $B$ be the result of turning $u$ in $A$. If $B=L'D_3v'E_3R'$ is of Type 2 with $\Lambda_r(v')\geq\tau+1$, the occurrence $\tilde{z}$ corresponding to $z$ in $B$ is well-defined and 
\[\Lambda_r(z)-(2\tau+1)<\Lambda_r(\tilde{z})<\Lambda_r(z)+(2\tau+1).\]
\end{lemma}
\begin{proof} 
Since $z$ does not contain $u$, by symmetry we may assume that $z$ is contained in $Lu=L'D_1u$. Thus we may write $z=z'X$ where $z'$ is contained in $L'$ and $X$ is a prefix of $D_1u$. By Lemma~\ref{overlap_type_occurrence} we have $\Lambda_r(X)<2\tau+1$ and so $\Lambda_r(z')>\tau+1$. Hence the the corresponding occurrence $\tilde{z}=z'Y$ in $B$ is well-defined and cannot contain $v'$  by Lemma~\ref{v_prolongation} since $\Lambda_r(v') \geq\tau + 1$ and $\Lambda_r(z') > \tau$. So 
$\tilde{z}=z'Y$ where $Y$  is a  proper prefix of $D_3v'$. 
Again by Lemma~\ref{overlap_type_occurrence} we have $\Lambda_r(Y)<2\tau+1$ and the result follows.
\end{proof}
We call $z'$ the {\bf remainder} of $z$ after turning $u$ (in analogy to $v'$ in Lemma~\ref{turn_res_complement}).

\begin{remark}
\label{influence_on_neighbours_complement_increasing_decreasing}
In the previous lemma, both $z$ and $\tilde{z}$ are prolongations of $z'$.  Since $\Lambda_r(z')\geq \tau+1$, either $z=\tilde{z}$ or one is a proper prefix of the other  by Remark~\ref{rem:occ with tau_1+1}. 
\end{remark}

\begin{remark}
\label{influence on neighbour small Q} If $A=LuR\in\Can_{r-1}$ where $ u, z$ are as in Lemma~\ref{influence_on_neighbours_complement} and the result $B=L'QR'$  of turning $u$  is  of Type 3 or of Type 1 or 2 with $\Lambda_r(v')<\tau+1$, then $\tilde{z}$ may contain $Q$  and we may have $\Lambda_r(\tilde{z})\geq\Lambda_r(z)+(2\tau+1)$. On the other hand, if the turn is of Type 3 in Lemma~\ref{turn_res_complement}, it is also possible that the occurrence $z$ is completely cancelled and has no trace in the result of the turn.
\end{remark}

\begin{corollary}
\label{close neighbours_stays_after_turn}
Let $A \in \Can_{r - 1}$ and let $u_1$, $u_2$  be  maximal occurrences of rank~$r$ in $A$  and $2\tau+1\leq \FracM_r(u_1) ,\Lambda_r(u_2) \leq n - 2\tau$. Write $A=Lu_1R$ and assume that $u_2$ is contained in $u_1R$. Let $v_1$ be the complement of $u_1$. Then  the result of turning $u_1$ in $A$ is of the form 
\[\can_{r-1}(Lv_1R)=L'D_3v_1'E_3M'u_2'R' \mbox{\quad for a suffix \ } M'u_2'R' \mbox{\ of  } R\] 
where $u_2'$ is a non-empty suffix of $u_2$ with
$\Lambda_r(u_2')>\Lambda_r(u_2)-(2\tau+1)$ and  $u_2'=u_2$ if $M'\neq 1$.
\end{corollary}
\begin{proof}
This follows immediately from (the symmetric version of) Lemma~\ref{influence_on_neighbours_complement} applied with $u=u_1$ and $z=u_2$.
\end{proof}

\begin{remark}\label{lem: two close neighbours}
Let $A\in\Can_{r-1}$ and let $u_1, u_2, u_3$ be maximal occurrences of rank~$r$  in $A$ enumerated from left to right. By Lemma~\ref{common_part_in_rank_r} the overlap of $u_2$ with  $u_1$ and $u_3$  has $\Lambda_r$-measure  $<\tau+1$. So if $\Lambda_r(u_2)\geq 2\tau+2$, then there is a subword of $u_2$ of  $\Lambda_r$-measure $> \Lambda_r(u_2)-(2\tau+2)$ not contained in either $u_1$ or $u_2$. 
In particular, if   $\Lambda_r(u_2)\geq 3\tau+2$ (or  $\geq 5\tau+2$, respectively), then $u_1, u_3$ are isolated (strongly isolated,  respectively) in $A$ (by a subword of $u_2$).
\end{remark}

The previous remark implies the following:
\begin{corollary}\label{cor: two close neighbours}
Let $A\in\Can_{r-1}$ and let $\X$ be a set of maximal occurrences of rank~$r$ in $A$ of $\Lambda_r$-measure $\geq 3\tau+2$. 
Then  any maximal occurrence in $\X$ to the left of $u\in\X$ is isolated from any maximal occurrence in $\X$ to the right of $u$, and so any $u\in X$ has at most one close neighbour in $\X$ on either side.
Similarly, if all occurrences in $\X$ have $\Lambda_r$-measure $\geq 5\tau+2$, then on either side of $u$ there is at most one maximal occurrence in $\X$ which is not strongly isolated  from $u$. 
\end{corollary}

\begin{lemma}\label{lem:isolated complements}
Let  $A \in \Can_{r - 1}$ and let 
  $u_1, u_2, u_3$ be  maximal occurrences of rank~$r$ in $A$ enumerated from left to right and of $\Lambda_r$-measure $\geq 3\tau+2$. Suppose $\FracM_r(u_2) \leq n - k(\tau +1)$ where $k$ is the number of close neighbours of $u_2$. Let $B$ be the result of turning $u_2$. Then $\tilde{u_1}$ and $\tilde{u_3}$ are isolated in $B$ (witnessed by a subword of $v_2'$) and strongly isolated in $B$ if $\FracM_r(u_2) \leq n - (5\tau + k\cdot(\tau +1))$
\end{lemma}
\begin{proof}
By Lemma~\ref{turn_res_complement} (ii) we have $\Lambda_r(v_2')>n - \Lambda(u_2) - 2\tau$ and by Lemma~\ref{influence_on_neighbours_complement} we know that  $\Lambda_r(\tilde{u_1}),\Lambda_r(\tilde{u_2})\geq \tau+1$ and so the corresponding occurrences are well-defined. Since $\tilde{u}_i$, $i = 1, 3$, can have an overlap with $\hat{v_2}$ only if $u_i$ is a close neightbour of $u_2$, the claim follows from Corollary~\ref{overlaps_size}.
\end{proof}

\subsection{Commuting turns of rank~$r$}
\label{consecutive_turns_section}

Our next aim is to show that the result of turning maximal occurrences of appropriate measures in $A\in\Can_{r-1}$ is independent  of the order in which we perform these turns.

\begin{corollary}\label{cor:isolated turns are independent}
Let $A =Lu_1 M u_2 R\in \Can_{r - 1}$, where $u_1$, $u_2$ are maximal occurrences of rank~$r$ in $A$ isolated from each other. If $\tau \leq \FracM_r(u_1),  \FracM_r(u_2)  \leq n - 2\tau$, then the result $C$  of turning $u_1$ and $u_2$ is independent of the order  in which we perform these turns and we have
\[C\equiv L\cdot a_1^{-n}\cdot u_1M\cdot a_2^{-n}\cdot u_2R\quad \mod \llangle \Rel_0, \ldots, \Rel_{r - 1}\rrangle.\]
\end{corollary}
\begin{proof}
By  Lemma~\ref{isolated_stays_after_turn} we know that the result  $B_1$ of turning $u_1$ in $A$ satisfies:
\begin{align*}
B_1=\can_{r-1}(Lv_1Mu_2R)&= L'D_3v_1'E_3M'u_2R \\
&\equiv L\cdot a_1^{-n}\cdot u_1M u_2R \quad \mod \llangle \Rel_0, \ldots, \Rel_{r - 1}\rrangle.
\end{align*}
It follows from this that for $M=M_0M'$ we have
\begin{equation}\label{eq:isolated}
 L\cdot a_1^{-n}\cdot u_1M_0\equiv  L'D_3v_1'E_3 \quad \mod \llangle \Rel_0, \ldots, \Rel_{r - 1}\rrangle.
\end{equation}
Similarly, for the result $C$ of turning $u_2$ in $B_1$ we obtain 
\[
\can_{r-1}(L'D_3v_1'E_3M'v_2R) \equiv L'D_3v_1'E_3M' \cdot a_2^{-n} \cdot u_2R \/ \mod \llangle \Rel_0, \ldots, \Rel_{r - 1}\rrangle.
\]
Combining this with \eqref{eq:isolated} we thus  obtain
\[
C \equiv  L\cdot a_1^{-n}\cdot u_1M_0M'\cdot a_2^{-n}\cdot u_2R \quad \mod \llangle \Rel_0, \ldots, \Rel_{r - 1}\rrangle.
\]

By considering inverses and using Remark~\ref{rem:turn induction beginning} we see that first turning $u_2$ and then $u_1$ yields the same result. So the claim follows from IH~\ref{IH idempotent} and IH~\ref{IH can of equivalent words}.
\end{proof}


\begin{definition}\label{def:solid}
Let  $A \in \Can_{r - 1}$ and let 
  $u_1, \ldots, u_t$ be maximal occurrences of rank~$r$ in $A$ enumerated from left to right with $\tau\leq\Lambda_r(u_i) \leq n - (3\tau + 1)$. Put $\Z=\{u_1,\ldots, u_t\}$.
    We call an occurrence~$u\notin \Z$, solid in $A$ with respect to $\Z$ if  after turning any subset of $\Z$ in any order the remainder of  $u$ (see the definition after Lemma~\ref{influence_on_neighbours_complement})  is an occurrence of measure $\geq\tau+1$. 
  We call the sequence   $u_1, \ldots, u_t$ solid if each $u_i, 1\leq i\leq t,$ is solid in $A$ with respect to $\Z\setminus\{u_i\}$.
 
  We say that the sequence $(u_0,\ldots, u_t)$ has a \emph{gap} at $i$ if  $u_i$ and $u_{i+1}$ are strongly isolated.
\end{definition}

The conditions imply in particular that all $u_i, v_i, i=1,\ldots, t,$   have $\Lambda_r$-measure $\geq \tau+1$ (and hence their maximal prolongations are unique).

 Note that for a solid set of occurrences each turn of one of the occurrences is of Type 2 and has an inverse turn. Clearly any subset of a solid set is again solid.

\begin{lemma}
\label{lem:close neighbour turns are independent}
Let $A=L\ul u_1, u_2\ur R \in \Can_{r - 1}$, where $u_1, u_2$ is a solid sequence of maximal occurrences of rank~$r$ in $A$. Then the result of turning $u_1$ and $u_2$ is independent of the order of the turns.
\end{lemma}
\begin{proof}
As in the proof of Corollary~\ref{cor:isolated turns are independent} the statement follows directly from the fact that both results of turns are equivalent $\mod \llangle \Rel_0, \ldots, \Rel_{r - 1}\rrangle$ and IH~\ref{IH can of equivalent words}.
\end{proof}

We now write $\epsilon = 2\tau + 1$.
\begin{proposition}
\label{prop:many turns}
Let  $A \in \Can_{r - 1}$ and let 
  $u_1, \ldots, u_t$ be maximal occurrences of rank~$r$ in $A$ enumerated from left to right and suppose that $u_i$ is an initial segment of $a_i^n\in\Rel_r$,  $i=1,\ldots, t$. Assume
\begin{enumerate}[label=(c\arabic*), ref=(c\arabic*)]
\item
\label{c2}
  all occurrences $u_i, 1\leq i\leq t$, are solid in $A$ with respect to $u_1,\ldots, u_t$;

\item\label{c1}
$\tau+1\leq \FracM_r(u_i) \leq n - (4\tau + 1 + k\cdot\epsilon)$ if $u_i$ has $k$ close neighbours\footnote{Note that  $0\leq k\leq 2$ by Condition~\ref{c2}.} among $u_j$, $j\neq i$.
\end{enumerate}
Then the result of turning (the occurrences corresponding to) the $u_i, i=1,\ldots, t$, (in the sense of Remark~\ref{rem:turn remainder}) is well-defined and independent of the order of the turns.
\end{proposition}

\begin{proof}
We know from Corollary~\ref{cor:isolated turns are independent} and Lemma~\ref{lem:close neighbour turns are independent} that under the given assumptions the turns of any two occurrences $u_i, u_j, i\neq j$, commute. Therefore the result follows once we establish that after turning an occurrence $u_i$ the maximal occurrences corresponding to the remaining occurrences $u_j, j\neq i$,  still satisfy the assumptions of this proposition.

By assumption on $\Lambda_r(u_i)$ and Lemma~\ref{turn_res_complement} Type~\ref{non_empty_v_part} we have $\Lambda_r(v_i')>\tau+1+k\cdot\epsilon$ where $k$ is the number of close neighbours of $u_i$ among the $u_j$. Furthermore, if $v_i'$ has overlap with $\tilde{u_m}$, $m\in\{i-1,i+1\}$, then by Lemma~\ref{isolated_stays_after_turn}, $u_i$ was a close neighbour of $u_m$. Since the overlap of $\tilde{u_m}$ and $v_i'$ has $\Lambda_r$-measure bounded by $\tau+1$, we see that after turning $u_i$, the occurrences $\tilde{u_{i-1}}, \tilde{u_{i+1}}$ are isolated from each other. Furthermore, $\Lambda_r(\tilde{u_j})=\Lambda_r(u_j)$ if $u_j$ and $u_i$ were isolated from each other, and
$\Lambda_r(\tilde{u_j})<\Lambda_r(u_j)+\epsilon$ if $u_m$ and $u_i$ were close neigbours for $m\in\{i-1,i+1\}$, in which case the number of close neighbours of $\tilde{u_m}$ among the $\tilde{u_j}, j\neq i,$ is exactly one less than the number of close neighbours of $u_m$ among the $u_j$. Since by Condition~\ref{c2} we have $\Lambda_r(\widetilde{u_j}) \geq\tau + 1$ we see that Condition \ref{c1} holds for $\{\tilde{u_j},1\leq j\neq i\leq t\}$.

Since after turning $u_i$ the occurrences $\tilde{u_{i-1}}, \tilde{u_{i+1}}$ are isolated from each other, clearly Condition~\ref{c2} continues to hold for $\{\tilde{u_j},1\leq j\neq i\leq t\}$ by Lemma~\ref{isolated_stays_after_turn}.
\end{proof}

\begin{definition}\label{def:stable}
We call a sequence of maximal occurrences in $A\in\Can_{r-1}$ stable if it satisfies Conditions~\ref{c2} and \ref{c1} from Proposition~\ref{prop:many turns}.
\end{definition}
\begin{remark}\label{rem:c1 and c2}
Note that if $A\in\Can_{r-1}$ and $\X$ is a set of maximal occurrences of rank~$r$ in $A$ where for each $u\in \X$ we have $5\tau+3\leqslant \Lambda_r(u)\leq n - 8\tau - 3$, then $\X$ is stable. Furthermore, $u_i, u_j\in\X$ are isolated for $|i-j|\geq 2$.
\end{remark}

To simplify notation we may now use the following convention:
\begin{convention}\label{simpler notation for occurrences 2}
If $A\in\Can_{r-1}$ and  $u_1,\ldots, u_t$ is a stable sequence of maximal occurrences  of rank~$r$ in $A$ with complements $v_1,\ldots, v_t$ and $B$ is the result of turning a subset of $\Z=\{u_1,\ldots, u_t\}$, then by Proposition~\ref{prop:many turns} we may simply denote  the maximal occurrences $\tilde{u_i}$ or $\hat{v_i}$ in $B$ by $u_i, v_i$, respectively.

Turning an occurrence $u_i$ can be considered as choosing the side $v_i$ in the relator $u_iv_i^{-1}$. Hence for choices $f_i\in\{u_i, v_i\}, 1\leq i\leq t,$ we will write $B=L'\ul f_1\ldots f_t\ur R'$ for the result of turning members of $\{u_i\in\Z \colon f_i=v_i\}$ in $A$, extending Convention~\ref{simpler notation for occurrences}.
\end{convention}

Note that the proof of Proposition~\ref{prop:many turns} shows the following important property:
\begin{corollary}
\label{cor:stability keeps}
Let  $A \in \Can_{r - 1}$, let $u_1, \ldots, u_t$ be a stable sequence and let $B$ be the result of turning $u_i$. Then the sequence (of occurrences corresponding to) $u_j$, $j \neq i$, is a stable sequence in $B$.
\end{corollary}

Informally speaking, solid occurrences prevent a turn of an occurrence on one side to influence occurrences on the other side:
\begin{lemma}
\label{lem:no accumulation}
Let $A \in \Can_{r - 1}$ and $u_1, \ldots, u_t$ be a stable sequence of maximal occurrences of rank~$r$ in $A$ and let $w$ be a solid maximal occurrence of rank~$r$ in $A$ with respect to $u_1, \ldots, u_t$. Then there exists  a unique maximal occurrence $\tilde{w}$ that corresponds to $w$ in the result of the turns of $u_1, \ldots, u_t$. Furthermore if $w$ is between $u_i, u_{i + 1}$ and not isolated from $k$ of them, then $|\Lambda_r(w) - \Lambda_r(\tilde{w})| < k\epsilon$ if $k \neq 0$ and $w = \tilde{w}$ if $k = 0$. 
\end{lemma}
\begin{proof}
By Proposition~\ref{prop:many turns}, we can turn $u_1, \ldots, u_t$ in any order. 
We do induction on $k=0, 1, 2$. The case $k=0$ follows from Lemma~\ref{isolated_stays_after_turn}.
Now assume $k>0$ and let $u_i$ be a close neighbour of $w$. Then after turning $u_i$ we have $|\Lambda_r(w) - \Lambda_r(\tilde{w})| < \epsilon$ by Lemma~\ref{influence_on_neighbours_complement}. Note that $\tilde{w}$ is uniquely defined because $w$ was assumed to be solid with respect to $u_1,\ldots, u_t$ and $\tilde{w}$ isolated from $u_{i-1}$ by Condition~\ref{c1} of a stable sequence. Furthermore,  $\tilde{w}$ lies between $ u_{i-1}, u_{i+1}$, is solid with respect to the stable sequence
$u_1,\ldots, u_{i-1}, u_{i+1},\ldots, u_t$ and is not isolated from $k-1$ of them. Thus the claim follows by induction.
\end{proof}

\begin{remark}
\label{rem:no accumulation}
 Lemma~\ref{lem:no accumulation} shows that $\tilde{w}$ only depends on the turns of $u_i, u_{i + 1}$.
\end{remark}

\begin{lemma}
\label{lem:stable with complements}
Let $A\in \Can_{r - 1}$ and $S=(q_1, \ldots, q_t)$ be a stable sequence of maximal occurrences of rank~$r$ in $A$. Let $S_0=(u_1, \ldots, u_s)$ be a subsequence with complements $v_1, \ldots, v_s$ such that $\Lambda_r(\hat{v}_i) \geq 5\tau + 3$. Let $B$ be the result of turning $u_1, \ldots, u_s$ and assume that all maximal occurrences in  $B$ have $\Lambda_r$-measure $\leq n - 8\tau - 3$. Then the maximal occurrences $\{Q\in S\setminus S_0\}\cup\{\hat{v_i}\mid i=1,\ldots, s\}$   form a stable sequence of rank~$r$ in $B$.
\end{lemma}
\begin{proof}
It suffices to verify that these occurrences are solid  in $B$. By Proposition~\ref{prop:many turns}, Corollary~\ref{cor:stability keeps}, Lemma~\ref{lem:no accumulation} and Remark~\ref{rem:no accumulation} it is enough to check that $q_i$ is solid in $A$ with respect to $q_{i - 1}, q_{i + 1}$ and $\widehat{v}_j$ is solid in the result of turning $u_j$ in $A$ with respect to the occurrences corresponding to $q_{j - 1}, q_{j+ 1}$. This follows from the initial assumptions.
\end{proof}

\begin{lemma}\label{lem:turns in power words}  
Let $A=LC^NR\in\Can_{r-1}$. Suppose $C=\ul u_1\ldots u_k\ur$ where $u_1,\ldots, u_k$ is a stable sequence of maximal occurrences $\geq 5\tau+3,  k\geq 2$. Then for $i=1,\ldots k,$ the result $B_i$ of turning all periodic shifts of $u_i$ in $C^N$ is of the form 
\begin{align*}
B_1 & =L'D_3v_1'E_3\ul u_2\ldots u_k\ur(\ul v_1 u_2\ldots u_k\ur)^{N-1} R;\\
B_i &=L(\ul u_1\ldots u_{i-1}v_iu_{i+1}\ldots u_k\ur)^N R \quad \mbox{ for } i\neq 1, k;\\
B_k &=L(\ul u_1\ldots u_{k-1} v_k\ur)^{N-1} \ul u_1\ldots u_{k-1}\ur D_3 v_k'E_3 R'.
\end{align*}
Furthermore, if $C=\ul u\ur$  contains a single maximal occurrence $u$ with $ 5\tau+3\leq\Lambda_r(u)\leq n-(3\tau+2)$, i.e. $A=L\ul u\ur^N R$, the result $B$ of turning all periodic shifts of $u$ inside $C^N$ is of the form \[B=L'D_3v'E_3(\ul v\ur)^{N-2}F_3v''G_3R'\] where $v'$ and $v''$ are respective remainders of the complements of maximal prolongations of $u$. 
\end{lemma}
\begin{proof}
This follows from Corollaries \ref{cor:prefix turn} and ~\ref{cor:turn protection} and their corresponding right version.
\end{proof}

\subsection{$\lambda$-semicanonical forms of rank~$r$}
\label{semican_section}
Recall that $\epsilon = 2\tau + 1$.
\begin{definition}[$\kappa\geq \frac{n}{2}$] 
\label{semicanonical_word_def}
A word in $\Can_{-1}$ is $\kappa$-bounded of rank~$r$ if all occurrences of rank~$r$ have $\Lambda_{j}$-measure $\leq \kappa$.
A $\kappa$-bounded word  from $\Can_{r - 1}$  is called \emph{$\kappa$-semicanonical of rank~$r$} and $\Scan_{\kappa, r}$ denotes the set of all $\kappa$-semicanonical words of rank~$r$.

\label{semicanonical_form_def}
If $A, A' \in \Can_{r - 1} $, $A' \in \Scan_{\kappa, r}$ and $A'$ and $A$ represent the same element of the group $\Fr/ \llangle\Rel_1,\ldots, \Rel_r\rrangle$, then $A'$ is called \emph{a $\kappa$-semicanonical form of rank~$r$ of~$A$}. 
\end{definition}
We emphasize that  $\kappa$-semicanonical forms of rank~$r$  are not  unique and that, by definition, $\Scan_{\kappa, r}\subseteq \Can_{r-1}$. Eventually we will have $\Can_r\subset\Scan_{\fracn+3\tau+1,r}$.

\begin{definition}\label{def:seam occurrence}
Let $A, C\in\Can_{r-1}$ and suppose that  either  $Z=L'QR'\in \Can_{r-1}$ is the result of a turn of a maximal occurrence $u$ of rank~$r$ in $A=LuR$ where $Q$ is $(3\tau+1)$-free or $Z=A\cdot_{r-1}C=A'QC'$ where $Q$ is $\tau$-free of rank~$r$. Suppose that $L', R'$ and $A', C'$ are $\kappa$-bounded.
If there is a unique maximal occurrence $w$ of $\Lambda_r$-measure $\geq\kappa+\epsilon$
in $Z$, then we call $w$ a \emph{seam occurrence} (with respect to $\kappa$).
A \emph{seam turn} is a turn of a seam occurrence.
\end{definition}
We collect a number of useful observations:

\begin{lemma}[$\kappa\geq\frac{n}{2}+\tau$]\label{lem:prolongation after turn}
Let $A=LuR \in \Can_{r - 1}$, where $L, R$ are $\kappa$-bounded in rank~$r$ and $u$ is a maximal occurrence of rank~$r$ in $A$  with $\kappa\leq\Lambda_r(u)<n$. Let $B=L'QR'$ be the result of turning $u$ in $A$.

\medskip\noindent
(i) If $\Lambda_r(v')\geq \tau+1$, then $B\in\Scan_{\kappa+\epsilon,r}$.

\medskip\noindent
(ii) If $B\notin\Scan_{\kappa+\epsilon,r}$, then $\Lambda_r(u)> n-(3\tau+1)$.

\medskip\noindent
(iii) If $B$ contains a maximal occurrence $w$ of $\Lambda_r$-measure $\geq \kappa+\tau+1$ containing $Q$, then $w$ is the  unique maximal occurrence of $\Lambda_r$-measure $\geq \kappa+\tau+1$.

\medskip\noindent
(iv) $B$ contains at most two maximal occurrences $w$ of $\Lambda_r$-measure $\geq \kappa+\tau+1$, one from the left of $Q$ and one from the right.

\medskip\noindent
(v)  If $B\notin\Scan_{\kappa+\epsilon}$, then  $B$  contains a unique occurrence  of $\Lambda_r$-measure $>\kappa+\epsilon$.

\end{lemma}
\begin{proof} We first note that $Q$ is $\kappa$-free. This is clear if the turn $A\mapsto B$ is of Type 3 or of Type 2 with $\Lambda_r(v') < \tau+ 1$ and follows from $ n-\kappa\leq \frac{n }{2}-\tau\leq \kappa-2\tau$ and Lemma~\ref{v_prolongation}  in case it is of Type 2  with $\Lambda_r(v')\geq \tau+ 1$.

\medskip\noindent
(i)  We have  $\tau+1\leq \Lambda_r(v')\leq \kappa-2\tau$ by Lemma~\ref{turn_res_complement}. Thus, $\Lambda_r(\hat{v})<\kappa$ by Lemma~\ref{v_prolongation} and the $\Lambda_r$-measure of  any other maximal occurrence in $L$ and $R$ can increase from the turn of $u$ at most by $\Lambda_r$-measure $<\epsilon$  by Lemma~\ref{overlap_type_occurrence}. Hence $B\in\Scan_{\kappa+\epsilon}$.

\medskip\noindent
(ii) If $B\notin\Scan_{\kappa+\epsilon,r}$, then $\Lambda_r(v')<\tau+1$ by part (i). Hence $\Lambda_r(u)>n-(3\tau+1)$.

\medskip\noindent
(iii) Let $w$ be a maximal occurrence in $B$ of $\Lambda_r$-measure $>\kappa+\tau+1$ containing  $Q$. If $\Lambda_r(Q)\geq\tau+1$, then $w$ is the unique occurrence containing $Q$ by Lemma~\ref{common_part_in_rank_r}.  If $\Lambda_r(Q)<\tau+1$, then $w$ contains a suffix of $L'Q$ and a prefix $QR'$ each of $\Lambda_r$-measure $\geq\tau+1$ and hence $w$ is the unique occurrence containing $Q$ by Lemma~\ref{common_part_in_rank_r}.
If $w'$ is another maximal occurrence in $B$ of $\Lambda_r$-measure $\geq\kappa+\tau+1$, then $w'$ is properly contained in $L'Q$ or $QR'$
and the overlap of $w'$ with $Q$ must have $\Lambda_r$-measure $\geq\tau+1$. Hence the overlap of  $w$ and $w'$ has $\Lambda_r$-measure $\geq\tau + 1$ again contradicting Lemma~\ref{common_part_in_rank_r}.

\medskip\noindent
(iv) If there are at least two maximal occurrences  in $B$ of $\Lambda_r$-measure $\geq\kappa+\tau+1$, then by part (iii) they must be contained in $L'Q$ or $QR'$. Since $L', R'$ and $Q$ are $\kappa$-bounded, any maximal occurrence  of $\Lambda_r$-measure $\geq\kappa+\tau+1$ in $L'Q$ contains both a suffix of $L'$ and a prefix of $Q$ of  $\Lambda_r$-measure $\geq\tau+1$ (and similarly for $QR'$). Hence such an occurrence is unique by Lemma~\ref{common_part_in_rank_r}.

\medskip\noindent
(v)  If  $B\notin\Scan_{\kappa+\epsilon}$, then $\Lambda_r(v')<\tau+1$-free by part~(ii) and so the turn is of Type 1 or 2 with $\Lambda_r(v')<\tau+1$ or of Type 3.
 If the turn $A\mapsto B$ is of Type 3, then $Q$ is $2\tau$-free. Hence any maximal occurrence  of $\Lambda_r$-measure $\geq\kappa+\epsilon$ contains $Q$ and so is unique  by part~(ii). 

Now assume that  the turn $A\mapsto B$ is of Type 1 or 2 and  $Q=Dv'E$ where $D, E$ are $\tau$-free and $\Lambda_r(v')<\tau+1$ so that $Q$ is $3\tau+1$-free.
 In this case any maximal occurrence $w$ of $\Lambda_r$-measure $\geq\kappa+\epsilon$ must contain  a prefix of $Dv'ER'$ and a suffix of $L'Dv'E$ of  $\Lambda_r$-measure $>\epsilon$ and hence there can only be one such $w$ by Lemma~\ref{common_part_in_rank_r}.\end{proof}

\begin{lemma}[$\kappa\geq\fracn+\tau$]\label{lem:product of lambda canonical}
Suppose $A, B\in \Scan_{\kappa,r}$. Then 
\[Z=\can_{r-1}(A\cdot B)=A'D_3B'\in\Scan_{\kappa+\epsilon,r}\]
unless $Z$ contains a  seam occurrence. 
\end{lemma}
\begin{proof} Since  $D_3$ is $\tau$-free and $A, B$ are $\kappa$-semicanonical, any  occurrence of $\Lambda_r$-measure $>\kappa+\epsilon$  in $A'D_3B'$ must contain both a suffix of $A'$ and a prefix of $B'$ of $\Lambda_r$-measure $\geq\tau+1$, and hence is unique.
\end{proof}

For a sufficiently big constant $\mu$ the natural greedy algorithm of turning occurrences of $\FracM_r$-measure $> \mu$ converges and leads to a $\mu$-semicanonical form of rank~$r$ of the word:
\begin{lemma}[$\mu\geq \frac{n}{2} + 9\tau $, $\alpha = 5\tau + 3$]
\label{lem:scan_greedy_alg}
If   $A \in \Can_{r - 1}$ and  $A=LuR \mapsto B$ is the turn  of the maximal occurrence $u$ of rank~$r$ in $A$ with $\FracM_r(u) > \mu$, then $d(B) < d(A)$ where  $d(X)$ denotes the sum of the $\FracM_r$-measures of all maximal occurrences of rank~$r$ in $X$ of $\FracM_r$-measure $\geqslant \alpha$ for $X\in\Can_{r-1}$.
\end{lemma}
\begin{proof}
Note that $d(A)-d(B)\geq\Lambda_r(u)- S$ where $S$ is the sum of $\Lambda_r$-measures of maximal occurrences
in $B=L'QR'$ that did not contribute to $d(A)$ but count for $d(B)$. These arise from maximal occurrences in $B$ of $\Lambda_r$-measure $\geq\alpha$ having nontrivial overlap with $Q$. Note that if $w=\ell q$ is a maximal occurrence in $A$ that contained in  $L'Q$, where $\ell$ is a suffix of $L'$ with $\Lambda_r(\ell)\geq\alpha$, then only $\Lambda_r(q)$ may contribute to $S$, and similarly for maximal occurrences contained in $QR'$. So in order to compute an upper bound for $S$, we may assume in such cases that $\Lambda_r(\ell)<\alpha$ and hence $\Lambda_r(\ell q)<\alpha+\Lambda_r(q)$.  Note that by Lemma~\ref{turn_res_complement} (iii) and Lemma~\ref{overlap_type_occurrence} we have $\Lambda_r(q) < 3\tau + 1$.

If $Q$ is $3\tau+1$-free, any maximal occurrence in $B$ that contributes to $S$ must contain a suffix of $L'$ or prefix of $R'$ (or both) and at least one of the suffix or prefix must have $\Lambda_r$-measure $\geq\tau+1$. By Lemma~\ref{common_part_in_rank_r} there can be at most one such occurrence from either side of $Q$ and only one if both overlaps with $L'$ and $R'$ are of $\Lambda_r$-measure $\geq\tau + 1$. Hence we can  estimate the contributions in $S$ by $S<2(\alpha+4\tau+2)=18\tau+10<\mu<\Lambda_r(u)$ (because $n > 18\tau + 20$).

By Lemma~\ref{turn_res_complement} (iii) it remains to consider the case that the turn is of Type 2, so $Q=D_3v'E_3$ (where $D_3, E_3$ may be empty) with $\Lambda_r(v')\geq\tau+1$.

Here contributions to $S$ can arise from $\hat{v}$ and, as before, from maximal occurrences containing a suffix of $L'$ or a prefix of $R'$. By Lemma~\ref{v_prolongation} we have $\Lambda_r(\hat{v})< \Lambda_r(v)+2\tau$, and so $\Lambda_r(u)-\Lambda_r(\hat{v})> \Lambda_r(u)-(n-\Lambda_r(u)+2\tau)  >2(9\tau)-2\tau=16\tau$. Furthermore the occurrences containing a suffix of $L'$ or a prefix of $R'$ may contribute at most $2\alpha+2(2\tau+1)=14\tau+8$. Hence again $\Lambda_r(u) - S > 16 \tau - 14\tau - 8 > 0$, and this finishes the proof.
\end{proof}

The previous lemma immediately implies:
\begin{corollary}[$\mu =\frac{n}{2} + 9\tau \geq n-7\tau-3$]
\label{scan_greedy_alg}
Any $A \in \Can_{r - 1}$ can be transformed into  a $\mu$-semicanonical form of rank~$r$ of~$A$ by a sequence of  turns of occurrences of rank~$r$ of $\FracM_r$-measure $> \mu$, starting from $A$.
\end{corollary}

While the previous algorithm is the most intuitive way to obtain a semi\-canonical form, the bound $\mu=\frac{n}{2}+9\tau$ will not be good enough for our purpose. Therefore we will further improve this bound below.

For future reference we record the following observation:
\begin{lemma}[$\kappa\geq\frac{n}{2}+\tau$]\label{lem:canonical implies turn restrictions}
If $A\mapsto B$ with $A, B\in\Scan_{\kappa,r}$ is the turn of a maximal occurrence $u$ in $A$ of Type 2, then $ n-\kappa-2\tau <\Lambda_r(u)\leq\kappa$.
\end{lemma}
\begin{proof} 
By Lemma~\ref{turn_res_complement} (ii) we have $n-\Lambda_r(u)-2\tau<\Lambda_r(v')\leq \kappa$.

\end{proof}

The following lemma will be used in Section~\ref{scan_equiv_section} to define an auxilliary group structure:
\begin{lemma}[$ \frac{n}{2}+\tau\leq \kappa\leq n-7\tau-3 $]\label{lem:mu' semicanonical}
Let $A,C\in\Scan_{\kappa,r}$ and  $Z=\can_{r-1}(A\cdot C)=A_1D_3C_1$. Then there is a sequence of seam turns  
\[Z=A_1D_3C_1\mapsto Z_2=A_2Q_2C_2\mapsto \ldots\mapsto Z_k=A_kQ_kC_k=Z'\in\Scan_{\kappa+(3\tau+1),r}\]
such that  $A_{i+1}, C_{i+1}$, $i\leq k-2$, are proper prefix and suffix of $A_i, C_i$, respectively,  $Q_i$ is $(3\tau+1)$-free of rank~$r$ for $i<k$, the last turn has $\Lambda_r$-measure $>\kappa+3\tau+1$ and all other turns have $\Lambda_r$-measure $>n-(3\tau+1)$. We write $Z'=\prod_{\kappa+(3\tau+1),r}(A\cdot C)$.  
\end{lemma} 
\begin{proof}
If $Z=A_1D_3C_1\notin\Scan_{\kappa+\epsilon,r}$, then, by Lemma~\ref{lem:product of lambda canonical}, $Z$ contains a seam occurrence $w$ of $\Lambda_r$-measure $>\kappa+\epsilon$ that properly contains $D_3$. Hence the result of turning $w$ in $Z$ is of the form $Z_2=A_2Q_2C_2$ where $A_2, C_2$ are proper prefix and suffix of $A_1, C_1$, respectively. If $Z_2\notin\Scan_{\kappa+(3\tau+1),r}$, then $\Lambda_r(w)>n-(3\tau+1)$  by Lemma~\ref{lem:prolongation after turn} (ii) and $Z_2$ contains a unique maximal occurrence $w_2$ of $\Lambda_r$-measure $\geq\kappa+(3\tau+1)$ by Lemma~\ref{lem:prolongation after turn} (v). Since $n - 3\tau - 1 > \kappa + 3\tau + 1)$, $w_2$ has non-trivial overlap both with $A_2$ and $C_2$.
Let $Z_3$ be the result of turning $w_2$ in $Z_2$. If $\Lambda_r(w_2)\leq n-(3\tau+1)$, then $Z_3\in\Scan_{\kappa+\epsilon,r}$  by Lemma~\ref{lem:prolongation after turn} (ii), and we are done. Otherwise  $w_2$ is the seam occurrence in $Z_2$ (of $\Lambda_r$-measure $>n-3\tau-1$).
We continue until $Z_k=Z'\in \Scan_{\kappa+3\tau +1,r}$. Since at each step we obtain a proper prefix  of $A_i$ and $B_i$  by the description in Lemma~\ref{turn_res_complement}, this process stops with $Z'=\prod_{\kappa+3\tau+1,r,}(Z)$ after finitely many turns.
\end{proof}

\begin{lemma}[$\frac{n}{2}+3\tau+1  \leq  \mu_2 \leq \mu_1\leq n-7\tau-3$]\label{scan_exists_simple_case}
For $A \in \Scan_{\mu_1,r}$  there exists a sequence of turns of rank~$r$ and $\FracM_r$-measure $> \mu_2 - \varepsilon$
\[A = C_1 \mapsto C_2\mapsto \ldots \mapsto C_l \in \Scan_{\mu_2,r}\] such that all $C_i$ are $(\mu_1 + \epsilon)$-semicanonical of rank~$r$.
\end{lemma}
\begin{proof}
Let $(u_0,\ldots, u_m)$ be an enumeration of all maximal occurrences in $A$ of $\Lambda_r$-measure $> \mu_2-2\epsilon$ enumerated from left to right. Note that this forms a stable sequence (see Remark~\ref{rem:c1 and c2}) and hence $u_i, u_j$ are isolated for $|i-j|\geq 2$.
Let $u=u_i$ be the left-most maximal occurrence of rank~$r$ in $A$ of $\FracM_r$-measure $>\mu_2 - \epsilon$ and write $A = LuR$. Then  $L$ is $(\mu_2 - \epsilon)$-bounded of rank~$r$ and $R$ is $\mu_1$-bounded. Since $ n - (3\tau + 1)\geq \mu_1\geq \Lambda_r(u)>\mu_2-\epsilon$, the result $B=L'\hat{v}R'$ of turning $u$ belongs to $\Scan_{\mu_1+\epsilon}$ by Lemma~\ref{lem:prolongation after turn}(i) and $\Lambda_r(\hat{v})<\mu_2-\epsilon$ by Lemma~\ref{turn_res_complement}. Furthermore,  $L'\hat{v}$ is $\mu_2$-bounded by Lemma~\ref{influence_on_neighbours_complement}.

 We consider in $B$ the left-most maximal occurrence $w$ in $\hat{v}R'$ of $\Lambda_r$-measure $> \mu_2-\varepsilon$. Then $w$ corresponds to $u_j$, $j > i$, in the original sequence. Since $\Lambda_r(u) \leqslant \mu_1 \leqslant n - (3\tau + 1)$, by Lemma~\ref{lem:prolongation after turn}~(i) $B \in \Scan_{\mu_1 + \varepsilon}$. Moreover, there exists at most one maximal occurrence of $\Lambda_r$-measure $> \mu_1$ in $B$ and if it exists, it is contained in $\hat{v}R'$ and must agree with $w$. Write $B = L_1wR_1$. Then $R_1$ is $\mu_1$-bounded. Now we turn $w$ and argue as above. Let $L_1'\widehat{z}R'$ be the result of the turn. Although $L_1$ is not $\mu_2 - \epsilon$-bounded now, since $\Lambda_r(\widehat{v}) \geqslant 3\tau + 2$, all maximal occurrences in $B$ from the left of $\widehat{v}$ stay unchanged in the result of the turn. Therefore $L_1'\hat{z}$ is $\mu_2$-bounded. We continue to argue in the same way until we reach the end of the sequence $(u_0,\ldots, u_m)$.
\end{proof}

{\bf From now on we fix $\lambda=\fracn+3\tau+1$.}

\begin{corollary}
\label{scan_exists}
Every  $A\in\Can_{r - 1}$ has a $\lambda$-semicanonical form 
of rank~$r$.
\end{corollary}
\begin{proof}
By Lemma~\ref{scan_greedy_alg} every $A\in\Can_{r-1}$ has a $(\frac{n}{2}+9\tau)$-semicanonical form. Now apply Lemma~\ref{scan_exists_simple_case} with $\lambda=\mu_2\leq\mu_1=\frac{n}{2}+9\tau\leq n-7\tau-3$. 
\end{proof}

While by Proposition~\ref{prop:many turns}  the result of a number of turns is independent of the order of the turns, the properties of the intermediate results  may depend on the order.

\begin{lemma}[$\kappa \geqslant 5\tau + 3$]
\label{scan_constant_for_consecutive_turns}
Let $A\in\Scan_{\kappa,r}$ and let $u_1, \ldots, u_t$ be a stable sequence of maximal occurrences of rank~$r$ in $A$ enumerated from left to right. Let $B$ be the result of turning all $u_i, i=1,\ldots, t$, and assume $B\in\Scan_{\kappa,r}$.

If $A = X_0 \mapsto  X_1 \mapsto \ldots \mapsto X_t=B$ is the  sequence of  turns from left to right, then $X_i\in\Scan_{\kappa+\epsilon,r}$ for $i=1,\ldots, t$.
\end{lemma}
\begin{proof}

Assume towards a contradiction that $X_i \notin \Scan_{\kappa+\epsilon,r}$ for some $i$. Let $w$ be a maximal occurrence of rank~$r$ in $X_i$ that does not correspond to $u_{i + 1}$. If there are occurrences to the left of $u$ that are not yet turned in $X_i$, then $w$ is to  the right of $\widetilde{u_{i + 1}}$. Since $u_{i + 1}$ is solid with respect to $u_1, \ldots, u_i$, the occurrence corresponding to $w$ in $A$  is equal to $w$ as a word. So, $\Lambda_r(w) \leqslant \kappa$.

If $\Lambda_r(w) > \kappa + \varepsilon$, then by Lemma~\ref{lem:no accumulation} $w$ does not correspond to any $u_j$, $j \geqslant i + 1$. So all occurrences that are not yet turned in $X_i$ are to the right of $w$ and $w$ is solid with respect to them. By Lemma~\ref{lem:no accumulation} the occurrence corresponding to $w$ in $B$ has $\Lambda_r$-measure $> \kappa$, a contradiction.
\end{proof}

\begin{corollary}[$\mu= n - 8\tau-3$]
\label{sequence_of_turns_to_consecutive_turns}
Let $X_0 \mapsto X_1 \mapsto \ldots \mapsto X_l$ be a sequence of turns of rank~$r$ and  $\FracM_r$-measure $\geq 9\tau+5$,  where  $X_0 \in \Scan_{\mu, r}$ and $X_i \in \Can_{r - 1}$. Then there exists  a stable sequence of maximal occurrences  $u_1, \ldots, u_t$ of rank~$r$ and $\Lambda_r$-measure $\geqslant 5\tau + 3$ in $X_0$  such that the result of the  corresponding turns is equal to~$X_l$.
\end{corollary}
\begin{proof}
The proof is by induction on $l$. If $l = 1$, there is nothing to prove. So consider $l>1$ and assume inductively that
there exists a stable sequence $(q_1, \ldots, q_s)$ of maximal occurrences of  $\FracM_r$-measure $\geq 9\tau+5$ in $X_0$ whose turns results in $X_{l - 1}$. If we turn them from left to right, then by Lemma~\ref{lem:no accumulation} every turn is of $\Lambda_r$-measure $< \mu + \varepsilon = n - 6\tau - 2$. Hence the maximal occurrence that contain the remainder of the complement has $\Lambda_r$-measure $> 4\tau + 2$. So again by Lemma~\ref{lem:no accumulation} the corresponding maximal occurrence in $X_{l - 1}$ is well defined has $\Lambda_r$-measure $> 2\tau + 1$. Let $z_1, \ldots, z_s$ be maximal occurrences in $X_{l - 1}$ that correspond to the remainders of the complements of $q_1, \ldots, q_s$.

Assume that $X_{l - 1} \mapsto X_l$ is the turn of an occurrence $\tilde{u}$. Then either $\tilde{u}$ coincides with some $z_i$, or $\tilde{u}$ lies between $z_i, z_{i + 1}$ for some $i$.
If $\tilde{u}$ coincides with some $z_i$, we put $u = q_i$. Now consider the second possibility. By the initial assumptions, $\FracM_r(\tilde{u}) \geq 9\tau + 5$. Hence $q_i$ and $q_{i + 1}$ are isolated in $X_0$ and $X_0 = Lq_iMq_{i + 1}R$. When we first turn the $q_j$, $j \neq i, i + 1$, then by Lemma~\ref{lem:no accumulation} the result is of the form $L_1\widetilde{q_i}M\widetilde{q_{i + 1}}R_1$, where $\Lambda_r(\tilde{q_i}), \Lambda_r(\tilde{q_{i + 1}}) < \mu + \epsilon = n - 6\tau - 2$. Since $\tilde{u} \geq 9\tau + 5$, we obtain that $X_{l - 1} = L_2{w'}_iE_3M'F_3{w'}_{i + 1}R_2$, where $w_i, w_{i + 1}$ are complements of $\tilde{q_i}, \widetilde{q_{i + 1}}$ and ${w'}_i, {w'}_{i + 1}$ are their remainders. Hence the common part of $M'$ and $\tilde{u}$ has $\Lambda_r$-measure $>5\tau + 3$. So its maximal prolongation $u$ in $X_0$ is unique. We denote it by $u$. 

If $u = q_i$ for some $i \in\{ 1, \ldots, s\}$, then we claim that $\{ q_1, \ldots, q_s\} \setminus \{ q_{i_0}\}$ is the required set of occurrences. Clearly, they form a stable sequence. By Proposition~\ref{prop:many turns} we may assume that turning $q_i$ is  the last turn and hence the turn $X_{l - 1}\mapsto X_l$ is its inverse. Therefore, $X_l$ is the result of turning the occurrences in $\lbrace q_1, \ldots, q_s\rbrace \setminus \lbrace q_i\rbrace$.

If $u\notin\{ q_1, \ldots, q_s\}$, then $\{ q_1, \ldots, q_s\} \cup \{ u\}$ is the required set of occurrences. Indeed, since $5\tau + 3 \leqslant \Lambda_r(q_i), \Lambda_r(u) \leqslant n - 8\tau - 3$, they form a stable sequence and clearly $X_l$ is the result of their turns.
\end{proof}

We need the following lemma  in Section~\ref{scan_equiv_section}.
\begin{lemma}[$\mu= n-8\tau-3,\lambda=\fracn+3\tau+1$]
\label{two_sequences_to_consecutive_turns}
Let $A \in \Scan_{\mu, r}$ and $A_1, A_2 \in \Scan_{\lambda + \tau, r}$. Assume that $A_1$ and $A_2$ are obtained from $A$ by sequences of turns of $\Lambda_r$-measure $\geq 9\tau + 5$. Then $A_1$ can be obtained from $A_2$ by a sequences of turns where all intermediate words are in $\Scan_{\lambda + \tau + \varepsilon, r}$.
\end{lemma}
\begin{proof}
By Lemma~\ref{sequence_of_turns_to_consecutive_turns}, there exist stable sets $\X_1, \X_2$ of maximal occurrences of rank~$r$ in $A$ of $\Lambda_r$-measure $\geqslant 5\tau + 3$ such that $A_i$ is the result of turning the occurrences in $\X_i, i=1,2$. Since all occurrences in $\X_i, i = 1, 2$ satisfy the restrictions $5\tau + 3 \leqslant \Lambda_r(u) \leqslant \mu$, $\X= \X_1 \cup \X_2$ is a stable set. Therefore Lemma~\ref{lem:stable with complements} implies that the maximal occurrences in $A_1$ that correspond to $\X$ (and to remainders of complements of turned occurrences) form a stable set $\Y$ in $A_1$. Clearly turns of some subset of $\Y$ in $A_1$ give $A_2$. Using Lemma~\ref{scan_constant_for_consecutive_turns}, we turn them from left to right and obtain the required sequence of turns.
\end{proof}

We next aim to show that turns commute (under suitable conditions) with the multiplication of canonical words. This will be  used in Section~\ref{scan_equiv_section}.

\begin{lemma}
\label{mergings_and_cancellations1}
Let $A = LuX, C \in \Can_{r - 1}$ for a maximal occurrence $u = a^ka_1, a^n\in\Rel_r,$ and let $B$ be the result of turning $u$. Assume that $\can_{r - 1}(A\cdot C) = A'DC'$ where $A', C'$ are  prefix and suffix of $A$ and $ C$, respectively, and $D$ is $\tau$-free of rank~$r$. Let $u'$ be the (possibly empty) common part of $u$ and $A'$. Then the following holds:

\medskip\noindent
(1) If $\Lambda_r(u')\geq\tau+1$, and $\tilde{u}$ is the maximal occurrence containing $u'$ in $A'D_3C'$, then the following diagram commutes:

\begin{center}
\begin{tikzpicture}
\node at (-.5, 0) {$A=LuX$};
\draw[thick, arrow = 0.99] (0.5, 0) to node [midway, above] {{\tiny turn of $u$}} (1.5, 0);
\node at (2, 0) {$B$};
\node at (-1, -2) {$\can_{r - 1}(A\cdot C)$};
\node at (3, -2) {$\can_{r - 1}(B\cdot C)$};
\draw[thick, arrow = 0.99] (.5, -2) to node [midway, above] {{\tiny turn of $\tilde{u}$}} (1.5, -2);
\draw[thick, arrow = 0.99] (0, -0.5) to (0, -1.5);
\draw[thick, arrow = 0.99] (2, -0.5) to (2, -1.5);
\end{tikzpicture}
\end{center}

\medskip\noindent
(2)  If $\Lambda_r(u) \geqslant \frac{n}{2} - 3\tau - 1$ and $\Lambda_r(u') < 2\tau+1$, then $C = X^{-1}c\inv R$, where $X\inv c\inv$ is the maximal cancellation in $uX\cdot C$. Then $c\inv$ is a fractional power of $a\inv$  with $\Lambda_r(c\inv)>\frac{n}{2} - (6\tau + 2)$. Let $\hat{w}$ be the 
maximal occurrence in $\can_{r - 1}(B\cdot C)$ corresponding to $c\inv$ (note that $\hat{w}$ then also corresponds to  $\hat{v}$, if this is defined). Then $\Lambda_r(\widehat{w}) > n - (4\tau + 1)$ and
the following diagram commutes:

\begin{center}
\begin{tikzpicture}
\node at (-.5, 0) {$A=LuX$};
\draw[thick, arrow = 0.99] (0.5, 0) to node [midway, above] {{\tiny turn of $u$}} (1.5, 0);
\node at (2, 0) {$B$};
\node at (-1, -2) {$\can_{r - 1}(A\cdot C)$};
\node at (3, -2) {$\can_{r - 1}(B\cdot C)$};
\draw[thick, arrow = 0.99] (1.5, -2) to node [midway, above] {{\tiny turn of $\hat{w}$}} (.5, -2);
\draw[thick, arrow = 0.99] (0, -0.5) to (0, -1.5);
\draw[thick, arrow = 0.99] (2, -0.5) to (2, -1.5);
\end{tikzpicture}
\end{center}

\end{lemma}
\begin{proof}
In the first case we can write $\can_{r - 1}(A\cdot C) = L\widetilde{u}R_1$ and let $Z$ be the result of the turn of $\widetilde{u}$. Then $Z \equiv  L a^{-n}\cdot \widetilde{u}R_1 \equiv L a^{-n}\cdot uX\cdot C \equiv \can_{r - 1}(L a^{-n}\cdot uX) \cdot C = B\cdot C \mod \llangle \Rel_0, \ldots, \Rel_{r - 1}\rrangle$. So the first part follows from IH~\ref{IH idempotent} and IH~\ref{IH can of equivalent words}.

For the second case if $\Lambda_r(u) \geqslant \frac{n}{2} - 3\tau - 1$ and $\Lambda_r(u')< 2\tau+1$, then $\Lambda_r(c) > \frac{n}{2} - (6\tau + 2)$. So $\can_{r - 1}(B\cdot C) = \can_{r - 1}(La^{-n}\cdot u \cdot c^{-1}R) = \can_{r - 1}(LwR)$, where $w = a^{-n}\cdot u \cdot c^{-1}$ is a fractional power of $a^{-1}$ with $\Lambda_r(w) > n - (3\tau + 1)$.

Now  write $LwR = La^N \cdot a^{-N}wR$. Since $\Lambda_r(u), \Lambda_r(c) \geqslant \tau$, IH~\ref{can_periodic_extension_constriction_hyp} and then Lemma~\ref{lem:product preparation} are applicable and imply that $\can_{r - 1}(LwR) = L'F_3w'R$, where $w'$ is a suffix of $w$ with $\Lambda_r(w') > \Lambda_r(w) - \tau > n - (4\tau + 1)$. So $w'$ has a unique maximal prolongation $\hat{w}$ and by Remark~\ref{rem:cyclic shift} and IH~\ref{IH can of equivalent words} the result of the turn of $\hat{w}$ is equal to $\can_{r - 1}(La^n\cdot wR) = \can_{r - 1}(La^n\cdot (a^{-n}\cdot u \cdot c^{-1})R) = \can_{r - 1}(A\cdot C)$. 
\end{proof}

\medskip\noindent
Note that by symmetry, if both $\tilde{u}$  and $\hat{w}$ are defined in $\can_{r-1}(A\cdot C)$ and $\can_{r-1}(B\cdot C)$, respectively, then both diagrams in the above situation commute.
Furthermore, if $\hat{v}$ is defined in Lemma~\ref{mergings_and_cancellations1}, then $\hat{w}$ comes from merging of $\hat{v}$ and $c^{-1}$ (this effect is described in Remark~\ref{influence on neighbour small Q}).

Similarly to Lemma~\ref{mergings_and_cancellations1}~(2) the following extension of Lemma~\ref{lem:close neighbour turns are independent} holds.
\begin{lemma}
\label{lem:mergings and cancellations 2}
Let $A = L\ul u_1 u_2\ur R\in \Can_{r - 1}$ be such that $\Lambda_r(u_1) \geqslant \fracn - \tau$, $u_1$ is a fractional power of $a^n \in \Rel_r$, $u_2$ is solid with respect to $u_1$ and $u_1$ is not solid with respect to $u_2$. Then $A = L\ul u_1 u_2\ur Xc^{-1}R_1$, where $c^{-1}$ is a fractional power of $a^{-1}$ with $\Lambda_r(c^{-1}) > \fracn - (9\tau + 3)$. Let $B_i$ be the result of turning $u_i$ in $A$ and $C$ be the result of turning the remainder of $u_2$ in $B_1$. Then there exists $\widehat{w}$ a maximal occurrence of rank~$r$ in $C$ that corresponds to $c^{-1}$ (and to $\widehat{v_1}$ if it is defined) such that $\Lambda_r(\widehat{w}) > n - (4\tau + 1)$, and the following diagram commutes:
\begin{center}
\begin{tikzpicture}
\node at (0, 0) {$A$};
\draw[thick, arrow = 0.99] (0.5, 0) to node [midway, above] {\tiny turn of $u_1$} (1.5, 0);
\node at (2, 0) {$B_1$};
\node at (0, -2) {$B_2$};
\node at (2, -2) {$C$};
\draw[thick, arrow = 0.99] (1.5, -2) to node [midway, above] {\tiny turn of $\hat{w}$} (.5, -2);
\draw[thick, arrow = 0.99] (0, -0.5) to node[midway, left] {\tiny turn of $u_2$} (0, -1.5);
\draw[thick, arrow = 0.99] (2, -0.5) to node[midway, right, text width=1.9cm] {\tiny turn of~$u_2'$} (2, -1.5);
\end{tikzpicture}
\end{center}
\end{lemma}
\begin{proof}
We can write $A = Lu'_1Mu_2R$, where $u'_1$ is a prefix of $u_1$ with $\Lambda_r(u'_1) > \Lambda_r(u_1) - \tau - 1$, $u'_1 = u_1$ if $M\neq 1$. Let $Z_1 = \can_{r - 1}(Lu'_1M)$, $Z_2 = \can_{r - 1}(u_2R)$, $W_1$ be the result of turning the occurrence that corresponds to $u'_1$ in $Z_1$ (which has $\Lambda_r$-measure $> \Lambda_r(u_1) - 2\tau - 1 \geqslant \fracn - 3\tau - 1$), and $W_2$ be the result of tuning the occurrence that corresponds to $u_2$ in $Z_2$. Then by IH~\ref{IH can of equivalent words} and Remark~\ref{rem:cyclic shift} $B_1 = \can_{r - 1}(W_1\cdot Z_2)$, $B_2 = \can_{r - 1}(Z_2\cdot W_2)$ and $C = \can_{r - 1}(W_1\cdot W_2)$. So the result follows from Lemma~\ref{mergings_and_cancellations1}~ (2) applied to $Z_1, W_1, W_2$.
\end{proof}

\section{Defining the canonical form of rank~$r$}

\subsection{Determining winner sides}
\label{can_of_semican_section}
Recall that $\tau =  15$, $\epsilon = 2\tau + 1$  and let 
$\mu =  n - 8\tau - 3\geq \frac{n}{2}+9\tau.$

In this section we define \emph{the canonical form of rank~$r$} of $A \in \Scan_{\mu, r}\subseteq \Can_{r - 1}$.
For this we consider all maximal occurrences of rank~$r$ in $A\in \Scan_{\mu, r}$ of $\FracM_r$-measure $\geq 5\tau + 3$.  Since  $\mu = n - (4\tau + 1 + 2\varepsilon)$, these occurrences form a stable sequence in $A$ and their complements are defined. Hence the result of turning any subset of these occurrences in $A$ is well-defined by Proposition~\ref{prop:many turns} and the canonical form $\can_r(A)$ is the result of turning a specific subset of these occurrences. Roughly speaking, for every maximal occurrence of rank~$r$ in $A$ of $\FracM_r$-measure $\geqslant \tau + 1 + 2\varepsilon$ we decide whether or not to turn it using  a threshold  of $\Lambda_r$-measure roughly $\fracn$.
For every maximal occurrence $u$ in $A$ at least one of $u$ or its complement will be below this threshold (see Corollary~\ref{cor:one side always certified}). 
\subsection{Rank $\mathbf{r = 1}$.} This case  is much simpler than the general case because relators in $\Rel_1$ are of the form $x^n$, where $x$ is a single letter, so maximal occurrences of rank~$1$ have no overlaps.
Since canonical triangles of rank~$0$ are trivial (i.e. all sides are equal to $1$), a turn of a rank~$1$ occurrence consists simply of replacing an occurrence by its complement.
Furthermore, for a maximal occurrence $u$ of rank~$1$ we have $\Lambda_1(u)=|u|$. Since the exponent $n$ is odd,  either $u$ or its complement has $\Lambda_1$-measure $<\frac{n}{2}$.

Now for $A \in \Scan_{\mu, 1}$, the canonical form of $A$ of rank~$1$, denoted by $\can_1(A)$, is defined as the word obtained from $A$ by replacing all maximal occurrences of rank~$1$ of $\Lambda_1$-measure $>\frac{n}{2}$ by their respective complements.

\begin{lemma}
\label{can_one_turn_rank_1}
Let $A\in \Scan_{\mu, 1}$, and let $A\mapsto B$ be a turn of rank~$1$. Assume that $B \in \Scan_{\mu, 1}$. Then $\can_1(A) = \can_1(B)$.
\end{lemma}
\begin{proof}
Since a turn of rank~$1$  just consist of replacing an occurrence by its complement, it does not change any other maximal occurrences and so this follows directly from the definition of the canonical form of rank~$1$.
\end{proof}

\medskip
\subsection{Rank $\mathbf{r \geq 2}$.}  From now on  until the end of Section~\ref{can_of_semican_section} we consider the general case, namely, rank~$r \geqslant 2$ and we fix the following set-up:  

Let $A$ be in $\mu$-semicanonical form, and let $u$ be a maximal occurrence of rank~$r$ in $A=LuR$ of $\FracM_r$-measure $\geqslant \tau + 1 + 2\varepsilon$ where $u = a^ta_1$, for some $a^n \in \Rel_r$, $a = a_1a_2$ ($a_1$ can be empty). 

We now state conditions whether or not to turn $u$ when we construct $\can_r(A)$. 
Let $\lambda_1$ and $\lambda_2$ be two constants with the following properties:
\begin{enumerate}[label=($\lambda$\arabic*), ref=($\lambda$\arabic*)]
\item 
\label{lambda_condition1}
$n-(11\tau+5)\geq\lambda_1>\lambda_2\geq \frac{n}{2} + 5\tau+2$.
\item
$\lambda_1 - \lambda_2 \geq \epsilon$
\end{enumerate}
For $n>36\tau + 16$ the interval $[ \frac{n}{2} + 5\tau+2, n-(11\tau+5)]$  has length $\geq 2\tau+1$ and hence such $\lambda_1, \lambda_2$ exist. 

\bigskip

We will use the fact that there exist sequences $m:\mathbb{N}\longrightarrow\{1, 2\}$  without subsequences of the form $BBb$ \cite{ChKa} where $b$ is a nontrivial initial segment of $B$.  By Proposition~\ref{scan_exists} we know that we can obtain $\lambda_2$-semicanonical forms by making a number of appropriate turns. In the certification process  we test whether a given occurrence can be made short enough by appropriate turns without significantly increasing other occurrences.
 The cubic-free sequence given by $m$ will ensure that we are not creating new power words (of higher rank)  in the process.

We first note the following:
\begin{lemma}\label{lem:coincide}
Let $u_1, \ldots, u_k$ be maximal occurrences of rank~$r$  in a word $W=D\ul u_1\ldots u_k\ur E$ such that $u_i, u_{i + 1}$ are not essentially isolated and $\FracM_r(u_i) \geqslant \tau + 1$ for all~$i$. Suppose that  $D, E$ are $\tau$-free of rank~$r$.  Let $u$ be a maximal occurrence of rank~$r$ in $W$ with $\Lambda_r(u)\geq 5\tau+2$. Then $u$  coincides with one of the $u_i$.
\end{lemma}
\begin{proof}
Clearly, if $u$ has nontrivial overlap with $D$ or $E$, then $u=u_1$ or $u_k$ respectively by Corollary~\ref{overlaps_size}.
If $u$ has a common part with the gap between $u_i, u_{i + 1}$ for some $1\leq i\leq k-1$, this common part has $\Lambda_r$-measure $<3\tau$. Since $\Lambda_r(u)\geq 5\tau+2$, the overlap of $u$ with  $u_i$ or $u_{i+1}$ has $\Lambda_r$-measure $\geq \tau+1$ and hence $u$ coincides with that occurrence.
\end{proof}

Recall that in a stable sequence any turn of a member of the sequence has an inverse turn by Lemma~\ref{inverse_turn_formal_existence}.

\begin{definition}[certification sequence]
\label{certification_def}
Let $A\in\Scan_{\mu,r}$. Then a stable sequence $(u=u_0, u_1, u_2,\ldots, u_t), t\geq 1,$ of maximal occurrences of $\Lambda_r$-measure $\geq 5\tau+2$ in $A=L\ul u_0\ldots u_t\ur R$ (enumerated from left to right) with complements $v=v_0, v_1, v_2,\ldots, v_t$ is called a \emph{certification sequence} in $A$ to the right of $u$ (with respect to $m:\mathbb{N}\to\{1,  2\}$) if the following holds

\begin{enumerate}

\item $u_1$ is essentially non-isolated from $u_0$;
\item there is a choice $f_i\in\{u_i, v_i\}, 0\leq i\leq t,$ such that in $W=L'\ul f_0f_1\ldots f_t\ur R'$ the maximal occurrences (corresponding to) $f_i$ for $i = 1, \ldots, t$
satisfy  $\Lambda_r(f_i)\leq\lambda_{m(i)}$.

\item After turning $f_0$ in $W$ and denoting the occurrence corresponding to $f_1$ in the result by $\tilde{f_1}$ we have $\Lambda_r(\tilde{f_1}) \geqslant \Lambda_r(f_1)$. Moreover if $\Lambda_r(\tilde{f_1}) = \Lambda_r(f_1)$, then $f_0 = u_0$.

\item For $2\leq i\leq t$, after turning $f_i$ in $W$ the occurrence  corresponding to $f_{i-1}$ has $\Lambda_r$-measure $>\lambda_{m(i-1)}$.

\item If there is a maximal occurrence $w$ in $A$ of  $\Lambda_r$-measure $\geq 5\tau+3$ to the right of $u_t$, then  after turning (the occurrence corresponding to) $w$ in $W$, in the resulting word we still have $\Lambda_r(f_t)\leq\lambda_{m(t)}$.
\end{enumerate}

\noindent
We say that the sequence \emph{certifies} $f_1$ to the right of $u$, i.e. either $f_1=u_1$ or $f_1=v_1$ is certified by the sequence.
$W$ is called the witness for the certification (of $u_1$ or $v_1$, respectively), exhibiting the choices $f_i\in \{u_i, v_i\}$.

We let $\Y_R(u)=\Y_R(u,A)$ denote the set of sides $f_1$ which are certified by a certification sequence to the right of $u$. (Note that if $f_1=v_1$ this is not an occurrence in $A$.)  Similarly we define $\Y_L(u,A)$ as the set  of inverses of $\Y_R(u\inv,A\inv)$ and put $\Y(u)=\Y_L(u)\cup\Y_R(u)$. Note that $\Y_L(u), \Y_R(u)$ contain at most two elements and are empty if there are no maximal occurrences of $\Lambda_r$-measure $\geqslant 5\tau + 3$ essentially non-isolated from $u$ from the left or right, respectively.

\noindent
We say that a stable sequence $(u=u_0,u_1, u_2,\ldots, u_t), t\geq 1,$  is an \emph{un-certification sequence} if it satisfies 1., 3. and 4. above and in place of 2. and 5. it satisfies the following:
\begin{enumerate}
\item[{2'.}] there is a choice $f_i\in\{u_i, v_i\}, 0\leq i\leq t,$ such that in $W=L'\ul f_0f_1\ldots f_t\ur R'$ the maximal occurrences (corresponding to) $f_i$ for $i = 1, \ldots, t-1$
satisfy  $\Lambda_r(f_i)\leq \lambda_{m(i)}$ and  $\Lambda_r(f_t)>\lambda_{m(t)}$.

\item[{5'.}] If there is a maximal occurrence $w$ in $A$ of  $\Lambda_r$-measure $\geq 5\tau+3$ to the right of $u_t$, then  after turning (the occurrence corresponding to) $w$ in $W$, in the resulting word we still have $\Lambda_r(f_t)>\lambda_{m(t)}$.
\end{enumerate}
\medskip 

Similarly we define (un-)certification sequences to the left in the obvious way by considering inverses.
We then say that a maximal occurrence $w$  or its complement of $\Lambda_r$-measure $\geq 5\tau+3$  contained in $uR$ is certified (or uncertified, respectively) in $A$ to the left of $u$ by a stable sequence $(u_t, \ldots, u = u_0)$ enumerated from right to left if $(u_0\inv, \ldots, u_t\inv)$ is an (un-)certification sequence  for $w\inv$ to the rigth of $u\inv$  in $A\inv$.  

If there is no maximal occurrence of $\Lambda_r$-measure $\geq 5\tau + 3$ to the right of $u$ and essentially non-isolated from $u$, then we say that $(u = u_0)$ is both the certification and uncertification sequence to the right of $u$.
\end{definition}

\begin{remark}\label{rem:certification non-isolated}
Let $A=LuR\in\Scan_{\mu,r}$ and let  $(u=u_0,\ldots, u_s)$ be an enumeration  of all maximal occurrences of $\Lambda_r$-measure $\geq 5\tau+3$ in $uR$ enumerated from left to right. Then $u_i, u_{i + 2}$ are strictly isolated, and hence, essentially isolated from each other. Combining Conditions 1, 2, 4, we see that  (un-)certification sequences have no gaps. Therefore, any (un-)certification sequence to the right of $u$ is an initial segment of $(u=u_0,\ldots, u_s)$.

For an (un-)certification sequence to the right of $u$ it suffices to check Conditions~5 and 5'  for the left most occurrence $w$ to the right of $u_t$ with $\Lambda_r(w) \geqslant 5\tau + 3$ because all maximal occurrences of rank~$r$ in $A$ to the right of $w$ are strictly isolated from~$u_t$.
\end{remark}

We first record the following remarks, which follow directly from Definition~\ref{certification_def}:

\begin{remark}\label{rem:certification}
Let $A=LuR\in\Scan_{\mu,r}$ and let  $(u=u_0,\ldots, u_t)$ be an (un-) certification sequence to the right of $u$ in $A$.
\begin{enumerate}
\item\label{prefix}
By  Condition~4, a proper prefix of  $(u=u_0,\ldots, u_t)$ can be neither a certification \ nor an un-certification sequence. 
\item For $i=0,\ldots t-1$, the members $u_i$ and $u_{i+1}$ are essentially non-isolated by Condition 4 and $\lambda_2-2\epsilon< \Lambda_r(u_i)< \mu$.
In particular, any (un-)certification sequence is stable.

\item\label{long_certification_for_one_side} If in $A$ we have $\Lambda_r(u_1)<\lambda_{m(1)}-\epsilon$ for a maximal occurrence $u_1$ not essentially isolated from $u$, then by Condition 3 and Lemma~\ref{influence_on_neighbours_complement}, $u_1$ is certified with the sequence $(u, u_1)$ and witness $A$.
Since $\lambda_2-\epsilon> \frac{n}{2}+\tau$, at least one of $u_1$ and $v_1$ is certified in~$A$ by Corollary~\ref{cor:one side always certified} with certification sequence $(u, u_1)$. So for at most one of $u_1$ and $v_1$ we have a certification or un-certification sequence that contains $> 2$ occurrences.

\item Suppose $(u_0=u, u_1, u_2,\ldots, u_t)$ is an (un-)certification sequence  to the right of $u$ with witness $W=L'\ul f_0\ldots f_t\ur R'$. If in $W$ we have  $\Lambda_r(f_i, W) \leq \lambda_{m(i)}-\epsilon$ (or  $\Lambda_r(f_i, W) \geq \lambda_{m(i)}+\epsilon$, respectively) for some  $1\leq i\leq t$, then $i=t$ by Conditions~2, 2' and~4.

\item\label{rem:necessary for certified} If $y$ is certified in $A$ to the right of $u$, then by Lemma~\ref{lem:no accumulation}
 $\Lambda_r(y)<\lambda_{m(1)}+k\epsilon$ where $k$ is the number of close neighbours of $y$ among $u_0, \ldots, u_t$ .
\end{enumerate}
\end{remark}

\begin{lemma}\label{lem:cert sequence exist}
Let $A = LuR\in\Scan_{\mu,r}$ where $\Lambda_r(u)\geq 5\tau + 3$. Let $u_1$ be a maximal occurrence of rank~$r$ in $uR$ essentially non-isolated from $u$ with $\Lambda_r(u_1) \geqslant 5\tau + 3$. Then for any choice $f_1 \in \{ u_1, v_1\}$ either there exists a unique certification sequence or  a unique un-certification sequence for $f_1$. In either case, the witness $W$ is unique.
\end{lemma}
\begin{proof}
Let $(u_0=u, u_1,\ldots, u_s)$ be an enumeration of all maximal occurrences of $\Lambda_r$-measure $\geq 5\tau + 3$ in $uR$ enumerated from left to right. Any (un-)certification sequence for $u$ is an initial segment of this sequence by Remark~\ref{rem:certification non-isolated}. 

Clearly there exists a unique choice for $f_0 \in \lbrace u_0, v_0\rbrace$ such that $(u_0, u_1)$ (for the choice for $f_1$) satisfies either Conditions~1--3, or Conditions~1, 2', 3. If it also satisfies one of Conditions~5 and~5', then $(u_0, u_1)$ is a certification or an un-certification sequence, respectively. If it satisfies neither Condition~5, nor Condition~5', then there exists $u_2$ such that $(u_0, u_1, u_2)$ satisfies either Conditions~1--4, or Conditions~1, 2', 3, 4. Therefore adding $u_i$ one by one, we eventually obtain either a certification, or an un-certification sequence. Moreover, by Conditions~2 and~4 the choice of $f_i \in \lbrace u_i, v_i\rbrace$ for every added occurrence, $i > 1$, is unique.
\end{proof}

\begin{remark}
Lemma~\ref{lem:cert sequence exist} implies that $f_1\in \lbrace u_1, v_1\rbrace$ cannot be both certified and ``uncertified''. So if there exists an un-certification sequence for $f_1$, then $f_1$ is not certified to the right of $u$. 
\end{remark}

The proof of Lemma~\ref{lem:cert sequence exist} shows that  certification sequences are equivariant under turns in the following sense:

\begin{corollary}
\label{cor:turn outside}
Let $A = L\ul u_0\ldots u_t\ur R \in \Scan_{\mu, r}$ where $(u_0 = u, \ldots, u_t)$, $t \geqslant 1$, is the (un-) certification sequence in $A$ for $f_1 \in \lbrace u_1, v_1\rbrace$. Let $w$ be a maximal occurrence in $u_tR$ of $\Lambda_r$-measure $\geq 5\tau+3$ with complement $y$ and assume $C \in \Scan_{\mu, r}$ is the result of turning $w$ in $C$.
Then the following holds:

\medskip\noindent
(i) If $\Lambda_r(u_t, C)<5\tau+3$, then $(u_0 = u, \ldots, u_{t-1})$,  is the certification sequence for $f_1$ in $C$. 

\medskip\noindent
(ii) If $\Lambda_r(u_t, C)\geq 5\tau+3$ and $C$ does not contain an occurrence  between $u_t$ and $y$ of $\Lambda_r$-measure $\geq 5\tau+3$, then $(u_0 = u, \ldots, u_t)$ is the (un-)certification sequence for $f_1$ in $C$.

\medskip\noindent
(iii) If $\Lambda_r(u_t, C)\geq 5\tau+3$ and $C$ contains an occurrence $z$ between $u_t$ and $y$ of $\Lambda_r$-measure $\geq 5\tau+3$, then $(u_0, \ldots, u_t)$ or  $(u_0 = u, \ldots, u_t, z)$ are  the (un-)certification sequence for $f_1$ in $C$. 

\medskip\noindent
Furthermore, $f_1$ is certified in $C$ to the right of $u$ if and only if this holds in $A$.

\medskip\noindent
(iv) If $t = 0$ and $C$ contains  a maximal occurrence $z$  essentially non-isolated from $u$ with $\Lambda_r(z, C)\geq 5\tau + 3$, then the complement of $z$ is not certified to the right of $u$ in $C$ by  Remark~\ref{rem:certification}(\ref{rem:necessary for certified}).
\end{corollary}
\begin{proof}
(i) and (ii) follow directly from the definition and
 the proof of Lemma~\ref{lem:cert sequence exist}. 
 
For part (iii) assume that there exists  a maximal occurrence $z$ of rank~$r$ in $C$ with $\Lambda_r(z, C) \geqslant 5\tau + 3$  with $\Lambda_r(z,A) < 5\tau + 3$. Then by Lemma~\ref{influence_on_neighbours_complement} $\Lambda_r(\widetilde{z}) < 7\tau + 4$. If $(u_0, \ldots, u_t)$ in $C$ still satisfies Condition~5 or~5', then $(u_0, \ldots, u_t)$ is a (un-)certification sequence in $C$. So assume that $(u_0, \ldots, u_t)$ in $C$ violates the corresponding condition (5 for a certification sequence and 5' for an un-certification sequence). This can happen only because of $z$. Consider the sequence $(u_0, \ldots, u_t, z)$ and the choice of $f_{t + 1} \in \{z,y\}$, where $y$ is the complement of $z$, such that this sequence satisfies Condition~2 or 2'. Since $(u_0, \ldots, u_t)$ does not satisfy Condition~5 or~5', we see that $(u_0, \ldots, u_t, z)$ satisfies Condition~4. If $f_{t + 1} = z$, then both $(u_0, \ldots, u_t)$ in $A$ and $(u_0, \ldots, u_t, \widetilde{z})$ in $C$ satisfy Conditions~2 and~5, so they both are certification sequences.

Assume that $f_{t + 1} = y$. Since $(u_0, \ldots, u_t)$ in $A$ satisfies Conditions~2' and ~5', the occurrence that corresponds to $f_{t + 1} = y$ after turning $u_i$ such that $f_i = v_i$, $0 \leqslant i \leqslant t$, has $\Lambda_r$-measure $> n - \Lambda_r(z) - 2\tau - \epsilon > n - (11\tau + 5) \geq \lambda_1$. Hence $(u_0, \ldots, u_t, z)$ in $C$ satisfies Conditions~2'. To see that it satisfies also Condition~5', let $A'$ be the result of turning $z$. Then $\Lambda_r(y, A') > n - (5\tau + 3) - 2\tau \geq \lambda_1 + 2\epsilon$. Let $g$ be the complement of $w$. Then $g$ is essentially non-isolated in $C$  from $z$ and $\Lambda_r(g) \geq 5\tau + 3$. So, we check Condition~5' for $(u_0, \ldots, u_t, z)$ in $C$ using $g$. Thus the occurrence that corresponds to $f_{t + 1} = y$ after turning $u_i$ such that $f_i = v_i$, $0 \leqslant i \leqslant t$, and the turn of the occurrence (corresponding to) $g$ has $\Lambda_r$-measure $> \Lambda_r(y, A') - 2\epsilon > \lambda_1$. Therefore, $(u_0, \ldots, u_t, z)$ in $A$ satisfies Condition~5' as required.

 For part (iii) we see  that $(u_0, \ldots, u_t)$ satisfies Conditions (1) -- (4) of Definition~\ref{certification_def}. If it also satisfies (5) or (5'), then $(u_0, \ldots, u_t)$ is the (un-)certification sequence for $f_1$ in $W$. Now suppose it does not satisfy either (5) or (5') and $(u_0, \ldots, u_t)$ is a certification sequence for $f_1$ in $A$. Then after turning $z$ in $C$ the occurrence corresponding to $f_t$ has $\Lambda_r$-measure $>\lambda_{m(t)}$. Since $\Lambda_r(f_t, C)\leq\lambda_{m(t)}$ and $\Lambda_r(z,C)<7\tau+4$, it follows that $(u_0, \ldots, u_t, z)$ is the certification sequence for $f_1$ in $C$.
 
On the other hand, if $(u_0, \ldots, u_t)$ is an un-certification sequence for $f_1$ in $A$, let $x$ be the complement of $z$ and $D$ be the result of turning $z$ in $C$. Then  $\Lambda_r(f_t, D)<\lambda_{m(t)}$.  Since $\Lambda_r(z,C)<7\tau+4$, it follows that $\Lambda_r(x, D)> n-9\tau-4 $. Hence by  $(u_0, \ldots, u_t, z)$ is the certification sequence for $f_1$ in $C$.
\end{proof}

\begin{corollary}\label{cor:invariance of cert sequence}
Let $A\in\Scan_{\mu,r}$, and let $(u_0,\ldots, u_t)$ be an (un-)certification sequence to the right of \label{cor:turns inside} $u_0$ in $A$ with witness $W=L\ul f_0\ldots f_t\ur R$. Let $g_i\in\{u_i, v_i\}, i=0,\ldots, t$, let $C$ be the result of turning all occurrences $u_i$ in $A$ with $g_i=v_i$ and suppose $C\in\Scan_{\mu,r}$ and $\Lambda_r(g_0, C) \geqslant 5\tau + 3$. Then the following hold:

\medskip\noindent
(i) $\Lambda_r(g_i, C) \geqslant 5\tau + 3$ for $0 \leqslant i \leqslant t - 1$.

\medskip\noindent
(ii) If $\Lambda_r(g_t, C) < 5\tau + 3$ and $t \geqslant 2$, then the sequence (corresponding to) $(g_0,\ldots, g_{t - 1})$ in $C$ is an (un-)certification sequence to the right of $g_0$ for the side that corresponds to $f_1$.

\medskip\noindent
(iii) Assume that $t = 1$ and $\Lambda_r(g_1, C) < 5\tau + 3$. If $g_1 = u_1$, then $v_1$ is not certified in $A$ to the right of $u$. If $g_1 = v_1$, then $u_1$ is not certified in $A$ to the right of $u$.

\medskip\noindent
(iv) Assume that $\Lambda_r(g_t, C) \geqslant 5\tau + 3$. Then either the sequence (corresponding to) $(g_0,\ldots, g_t)$ in $C$, or $(g_0,\ldots, g_t, \widetilde{z})$ is a (un-)certification sequence to the right of $g_0$ for the side that corresponds to $f_1$, where $\widetilde{z}$ corresponds to some maximal occurrence $z$ in $A$ with $\Lambda_r(z) < 5\tau + 3$.

\medskip\noindent
Furthermore, the choices for $g_i$ in (un-)certification sequences in (ii) and (iv) agree with the choices $f_i$ in the initial sequence $(u_0,\ldots, u_t)$.
\end{corollary}
\begin{proof}
If $t = 1$, then (i) immediately holds. So let $t \geqslant 2$ and assume towards a contradiction that $\Lambda_r(g_i, C) < 5\tau + 3$ for some $1 \leqslant i \leqslant t - 1$. If $g_i = f_i$, then $\Lambda_r(f_i, W) < \Lambda_r(g_i, C) + 2\varepsilon < 9\tau + 5 < \lambda_2 - \varepsilon$,  contradicting Remark~\ref{rem:certification}~(iii). If $g_i \neq f_i$ (i.e. $g_i \in \lbrace u_i, v_i\rbrace \setminus \lbrace f_i\rbrace$), then $\Lambda_r(f_i, W) > n - \Lambda_r(g_i, C) - 2\tau - 2\varepsilon > \lambda_1$, which contradicts to Condition~2 or~2'.

(ii)--(iv) are proved as in Corollary~\ref{cor:turn outside}.
\end{proof}

\begin{lemma}\label{lem:context introduction} 
Let $A\in\Scan_{\mu,r}$, let $u$ be a maximal occurrence of $\Lambda_r$-measure $\geq 5\tau + 3$ and let $(u=u_0,\ldots, u_t), t\geq 1,$ be an (un-)certification sequence for $f_1\in\{u_1,v_1\}$ in $A$ to the right of $u$. Let $\kappa=\Lambda_r(f_t, W)$.  Write (in the notation of Convention~\ref{simpler notation for occurrences})  $A=L\ul u_0\ldots u_t\ur MR$ and let  $B= L'v^{\prime}_tEM'R$ be the result of turning $u_t$ in $A$, where $M^{\prime}$ is a suffix of $M$ and $E$ is $\tau$-free of rank~$r$, such that 
\begin{itemize}
\item if $\kappa\leq\lambda_{m(t)}-\epsilon$ or  $\kappa\geq\lambda_{m(t)}+\epsilon$, then $M$ or $EM'$ contain an occurrence $a^\tau M_0b^\tau$, $a^n, b^n \in \Rel_r$;
\item if $\lambda_{m(t)} - \epsilon < \kappa < \lambda_{m(t)} + \epsilon$, then $M$ or $EM'$ contain a  strong separation word (see Definition~\ref{isolated_occurrences}).
\end{itemize} 

Then for any $A'=L\ul u_0\ldots u_t\ur MR'\in\Scan_{\mu,r}$ the corresponding sequence  $(u_0,\ldots, u_t), t\geq 1,$ in $A'$ is still an (un-)certification sequence for the corresponding $f_1$ in $A'$ (with the same choices for all $f_i \in \lbrace u_i, v_i\rbrace$).
\end{lemma}
\begin{proof}
This follows directly from the definition, Corollaries~\ref{cor:turn protection aux} and~\ref{cor:turn protection} and Definition~\ref{isolated_occurrences}.
\end{proof}

\begin{definition}
\label{can_context_def}
Let $A = LuMR \in \Scan_{\mu, r}$, where $u$ is a maximal occurrence of rank~$r$  with $\frac{n}{2} - 5\tau - 2 < \FracM_r(u) < \frac{n}{2} + 5\tau +~2$. We say that  $uM$  is a \emph{right context} for $u$ in $A$ if any (un-)certification sequence on the right of $u$ in $A$  is properly contained in $uM$ and for any word $A'=LuMR' \in \Scan_{\mu, r}$ the sequence of corresponding occurrences is an (un-)certification sequence to the right of $u$ in $A'$ for the same $f_1$.
\end{definition}
Note that such $M$ might not exist and in this case the right context is not defined.
\smallskip

Let $A=LuR\in\Scan_{\mu.r}$ for a maximal occurrence $u$ of rank~$r$.
When we  decide for an occurrence $u$ in a word $A\in\Scan_{\mu,r}$ whether to turn it (and thus  replace it essentially by the complement) we need to take into account the effect the turn has on the neighbouring occurrences because we want the canonical form to be invariant under certain turns. We therefore make this decision after also considering the neighbours of $u$ and their possible turns.
We first note that an occurrence with sufficiently small $\Lambda_r$-measure will always be shorter than its complement no matter which neighbours we turn,  and, conversely, if the $\Lambda_r$-measure of an occurrence is sufficiently large, then no matter which neighbours we turn, it will always be the longer than its complement:

\begin{lemma}
\label{too_small_and_too_big_side} Let $A=L\ul y_1, u, y_2\ur R\in\Scan_{\mu, r}$, $u, y_1, y_2$ be maximal occurrences of rank~$r$ with complements $v, z_1, z_2$, respectively,  and assume that $\Lambda_r(y_i)\geq 3\tau~+~2, i=1, 2$. Let $B$ be the result of turning $u$ in $A$ and for choices
$f_i, g_i\in \{y_i, z_i\}, i=1, 2,$ let $A'=L'\ul f_1, u, f_2\ur R'$,  $B'=L''\ul g_1, v, g_2\ur R''$,
\begin{enumerate}
\item  If \quad $5\tau + 3\leq \FracM_r(u) \leq\frac{n}{2} -5\tau - 2$, then  $\Lambda_r(u, A') <\Lambda_r(v, B')$.

\item If \quad $\frac{n}{2} +5\tau+2\leq \FracM_r(u) \leq\mu$, then $\FracM_r(u, A') > \FracM_r(v, B')$.
\end{enumerate} 
\end{lemma}
\begin{proof}
1. If $\FracM_r(u) \leq \frac{n}{2} - 2\epsilon-\tau$, then $\FracM_r(v) = n - \Lambda_r(u) \geq  \frac{n}{2} + 2\epsilon+\tau$.
So after possibly turning neighbours of $u$ by Lemma~\ref{influence_on_neighbours_complement} the corresponding occurrence $u$ satisfies
$\Lambda_r(u, A')<\Lambda_r(u)+2\epsilon\leq \frac{n}{2}-\tau$
whereas $\Lambda_r(v, B')> \Lambda_r(v, B) - 2\varepsilon > \Lambda_r(v) - 2\tau - 2\varepsilon \geqslant \frac{n}{2} - \tau$.  

\medskip\noindent
2. If $\FracM_r(u) \geq\frac{n}{2} + 2\epsilon+\tau$, then $\FracM_r(v, B) \leq  \frac{n}{2} - 2\epsilon - \tau$. By Lemma~\ref{influence_on_neighbours_complement} we have that $\FracM_r(u, A') > \Lambda_r(u, A) - 2\varepsilon \geqslant \frac{n}{2} +\tau$ and $\Lambda_r(v, B') < \Lambda_r(v, B) + 2\varepsilon < \Lambda_r(v) + 2\tau + 2\varepsilon \leqslant \frac{n}{2} + \tau$.
\end{proof}
So in these cases, no matter which neighbours we turn, the occurrence corresponding to $u$ remains  shorter (or longer, respectively) than the one corresponding to $v$.
Thus, according to our definition we never turn an occurrence $u$ of $\Lambda_r$-measure $\leq \frac{n}{2} - 5\tau - 2$ and we always turn  an occurrence $u$ of $\Lambda_r$-measure $\geq \frac{n}{2} + 5\tau + 2$.
Therefore we can now restrict our attention to occurrences $u$ with  $ \frac{n}{2} - 5\tau - 2< \FracM_r(u) < \frac{n}{2} + 5\tau +2$. Note that in this situation all occurrences to the left of $u$ are strictly isolated from all occurrences to the right of $u$. Therefore
we can consider the left and right side separately. 

We now define $\tilde{u}$ to be the shortest occurrence among the occurrences corresponding to $u$ 
when we turn neighbours of $u$ of $\Lambda_r$-measure $\geq 5\tau + 3$ according to the certified sides in $\Y(u)$.
If $\Y(u)=\emptyset$, then $\tilde{u}=u$.
Also we define $\tilde{v}$ to be the shortest occurrence among the occurrences corresponding to $\hat{v}$ using the same set $\Y(u)$.

If $|\tilde{u}|\neq |\tilde{v}|$, then we choose the shorter occurrence as the winner side.
Since the canonical form is equivariant with respect to inversion, we need to make sure that the choice for $A=LuR$ is consistent with the choice for $A\inv$. Therefore we use the following more intricate procedure to determine the winner side in case $|\tilde{u}|= |\tilde{v}|$:

Consider $a^n\in\Rel_r$ as a cyclic word. Let $I_u$ be the starting point of $u$, $F_u$ be the end point of $u$.
Since $u$ contains at least one period of $a$, $I_u$ and $F_u$ are fixed up to a cyclic shift by some  number of periods of $a$. Following the construction of $\tilde{u}$ and $\tilde{v}$ from $u$, we mark the initial and the final points of $\tilde{u}$ and $\tilde{v}$ in $a^n$ with respect to the points $I_u$ and $F_u$. Denote them by $I_{\tilde{u}}$,  $F_{\tilde{u}}$, $I_{\tilde{v}}$ and $F_{\tilde{v}}$ respectively. Notice that $\tilde{u}$ and $\tilde{v}$ may or may not have overlaps in the cyclic word $a^n$ and so the overlaps or gaps between  $\tilde{u}$ and $\tilde{v}$ have measure $<3\tau+1$.

Consider the subword of $a^n$ with endpoints $[I_{\tilde{u}}, I_{\tilde{v}}]$ of $\Lambda_r$-measure $<3\tau+1$ and let $c$ denote the middle letter if the length of this is odd, otherwise let $c$ mark the mid point between the two middle letters. Similarly,  consider the subword of $a^n$ with endpoints $[F_{\tilde{u}}, F_{\tilde{v}}]$ and define $d$ in the same way. We denote the segment corresponding to $c$ and $d$, respectively, by the (possibly empty) intervals $[P_1, P_1']$ and $[P_2, P_2']$ (see diagram below).
  Let $u_0$ be the subword of $a^n$ starting at  $P_1'$ and ending  at $P_2$, and let $v_0$ be the subword of $a^{-n}$ starting at $P_1$ and ending at $P_2'$. So, we have a partition of $a^n$ into four segments: $u_0$, $d$, $v_0\inv $, $c$, where $|c|, |d|\leq 1$. 

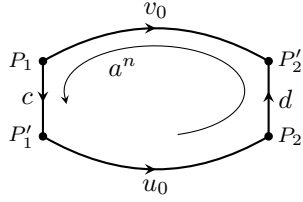
\begin{center}
\begin{tikzpicture}[scale = 2]
\draw[black, thick, arrow=0.5] (1.5+4, 0) to [bend right] node[midway, below] {$u_0$} (3+4, 0);
\node[left] at (1.5+4, 0) {\footnotesize $P_1'$};
\node[left] at (1.5+4, 0.5) {\footnotesize $P_1$};
\filldraw (1.5+4, 0) circle (0.7pt);
\filldraw (1.5+4, 0.5) circle (0.7pt);
\draw[black, thick, arrow=0.5] (1.5+4, 0.5) to [bend left] node[midway, above] {$v_0$} (3+4, 0.5);
\draw[black, thick, reversearrow=0.6] (1.5+4, 0) to node[midway, left] {$c$} (1.5+4, 0.5);
\draw[black, thick, arrow=0.6] (3+4, 0) to node[midway, right] {$d$} (3+4, 0.5);
\node[right] at (3+4, 0) {\footnotesize $P_2$};
\node[right] at (3+4, 0.5) {\footnotesize $P_2'$};
\filldraw (3+4, 0) circle (0.7pt);
\filldraw (3+4, 0.5) circle (0.7pt);
\draw[reversearrow=0.01] (5.65, 0.25) arc (190:-75: 0.6 and 0.3) node[pos=0.3, below] {$a^n$};
\end{tikzpicture}
\end{center}

 Now we are ready to specify the conditions for turning  $u$ in $A$ in order to construct $\can_r(A)$.

\begin{remark}
\label{no_square_relators}
Note that for $a^n \in \Rel_r$ we have $a^n \neq Z^2$ for all $Z \in \Canc_0$:
if $a^n=Z^2$, then $Z\in Cen(a^n)$. By definition of $\Rel_r$ we have $Cen(a^n)=\langle  a\rangle $, so $Z=a^k$ for some $k\in\mathbb{Z}$. Then $Z^2=a^{2k}\neq a^n$ since $n$ is odd.
\end{remark}

The following is well-known:
\begin{lemma}
\label{common_part_with_inverse}
Let $b\neq 1$ be cyclically reduced. If $b=xy$ where $|x|\geq \frac{1}{2}|b|$, then no cyclic shift of $b$ contains $x\inv$ as a subword.
\end{lemma}
\begin{proof}
Suppose otherwise. Then either there exists an occurrence of $x\inv$ in $b$, or $b = x_1\inv zx_2\inv$, where $x= x_1x_2$. Since $|x| \geq \frac{b}{2}$ and $b$ is a reduced word, in either case we obtain that $x$ and $x\inv$ have an overlap which is impossible.
\end{proof}

Recall that $a^n \in \Rel_r$ as a cyclic word is separated into four parts $u_0$, $d$, $v_0\inv $, $c$, where $c$ and $d$ are either empty, or a single letter (independently from each other).

\begin{lemma}
\label{different_sides}
For $a^n\in \Rel_r, r\geq 2,$ the sets of words
\begin{align*}
&\lbrace u_0,\ u_0\inv ,\ cu_0,\ u_0d,\ u_0\inv c\inv ,\ d\inv u_0\inv \rbrace,\\
&\lbrace v_0,\ v_0\inv ,\ c\inv v_0,\ v_0d\inv ,\ v_0\inv c,\ dv_0\inv \rbrace.
\end{align*}
are not equal to each other.
\end{lemma}
\begin{proof}
Since $\Rel_r$ is invariant under cyclic shifts by IH~\ref{IH Rel}, we may assume $a^n=u_0dv_0\inv c$. 
Now assume to the contrary that
\begin{equation*}
\lbrace u_0,\ u_0\inv ,\ cu_0,\ u_0d,\ u_0\inv c\inv ,\ d\inv u_0\inv \rbrace = \lbrace v_0,\ v_0\inv ,\ c\inv v_0,\ v_0d\inv ,\ v_0\inv c,\ dv_0\inv \rbrace.
\end{equation*}

The words $u_0$ and $u_0\inv $ are the shortest  in the left-hand set, and similarly $v_0$ and $v_0\inv $ are shortest on the right-hand side. Hence, either $u_0 = v_0$, or $u_0 = v_0\inv $. Since $v_0$ is a subword of $a^{-n}$ and $\FracM_r(u_0), \FracM_r(v_0) \geqslant 1$, by construction, we have $u_0\neq v_0$ by Lemma~\ref{common_part_with_inverse} and so $u_0 = v_0\inv $.
If both $c$ and $d$ are empty, then $a^n=u_0v_0\inv=u_0^2$, contradicting Remark~\ref{no_square_relators}. So, we may assume that at least one of $c$ and $d$ is not empty. By symmetry assume that $d\neq 1$.

For the sets to be equal, we must have $u_0d\in\{c\inv v_0,  v_0d\inv,  v_0\inv c, dv_0\inv\}$. Since $u_0=v_0\inv$,  we have $u_0d\notin\{c\inv v_0, v_0d\inv\}$   by Lemma~\ref{common_part_with_inverse}. 
If $u_0d=v_0\inv c=u_0 c$, then $c=d$ and hence $a^n=(u_0d)(v_0\inv c)=(u_0d)^2$, contradicting Lemma~\ref{no_square_relators}.
And finally if $u_0d=dv_0\inv =du_0$, then $d\in Cen(u_0)$ and hence $u_0\in\llangle  d\rrangle $. Since $\Lambda_r(u_0)\geq 1$, we also have $a\in  \llangle  d\rrangle $. However, since $r\geq 2$ we have   $|a|>1$  by definition. This contradiction proves the lemma.
\end{proof}

\begin{definition}[$\deglex$ order]
\label{deglex_order}
Fix an ordering on the set of letters. For reduced words $C_1, C_2$  we say that $C_1 <_{\deglex} C_2$ if either $|C_1|< |C_2|$, or $|C_1| = |C_2|$ and $C_1$ is lexicographically smaller than $C_2$ with respect to the order that we fixed on the letters. For finite sets of words $\U\neq \V$, we put $\U <_{\deglex} \V$ if the minimal element of $\U\cup\V\setminus(\U\cap \V)$ belongs to $\U$.
\end{definition}
Now let $\U=\{u_0,\ u_0\inv ,\ cu_0,\ u_0d,\ u_0\inv c\inv ,\ d\inv u_0\inv \}$ and $\V=\{v_0,\ v_0\inv ,\ c\inv v_0,\ v_0d\inv ,\ v_0\inv c,\ dv_0\inv\}$.  By Lemma~\ref{different_sides} we have $\mathcal{U} \neq \mathcal{V}$. If $\U<_{\deglex} \V$, we do not turn $u$ and call $u$ the winncer side, otherwise we turn $u$ and $v$ is called the winner side.

Thus $v$ is the winner side if and only if either  $|\tilde{u}|>|\tilde{v}|$ (as defined above) or, in case $|\tilde{u}|=|\tilde{v}|$, if $\mathcal{V}$ is smaller than $\mathcal{U}$ with respect to the $<_{\deglex}$ order.

\begin{lemma}\label{rem:winners are certified}
Let $A=LuR\in\Scan_{\mu,r}$ for a maximal occurrence $u$ of rank~$r$ and suppose that $q$ is a maximal occurrence of $\Lambda_r$-measure $\geq 5\tau+3$ essentially non-isolated from $u$ to the right of $u$. Then the winner side for $q$ is certified to the right of  $u$.
 \end{lemma}
 \begin{proof}
Indeed, if the winner side for $q$ is not certified, then by turning occurrences not essentially isolated from $q$ of $\Lambda_r$-measure $\geqslant 5\tau + 3$, the  maximal occurrence corresponding to the winner side can be made $> \lambda_2 = \frac{n}{2} + 5\tau + 2$ contradicting Lemma~\ref{too_small_and_too_big_side}.
 \end{proof}

\begin{definition}[canonical form for $\lambda$-semicanonical words, $\lambda=\fracn+3\tau+1$]\label{def:can_r lamda-semi}
For $A\in\Scan_{\lambda,r}$,   consider the set of all maximal occurrences  of $\Lambda_r$-measure $\geq \frac{n}{2}-5\tau-2$ and turn each one of them according to the decision process described above. The result is denoted by $\can_r(A)$.
\end{definition} 
By Proposition~\ref{prop:many turns}, the result $\can_r(A)$ does not depend on the order in which we perform the necessary turns.  

\begin{lemma}\label{can subset scan}
For $A\in\Scan_{\mu,r}$ we have $\can_r(A) \in \Scan_{\lambda,r}$.
\end{lemma}
\begin{proof}
Let $A\in\Scan_{\mu,r}$  and $u$  some maximal occurrence of rank~$r$ in $A$, let $f \in \lbrace u, v\rbrace$ be the winner side for $u$ and let $q_1$, $q_2$ be maximal occurrences in $A$ of $\Lambda_r$-measure $\geqslant 5\tau + 2$ not essentially isolated from $u$ on the left and right, respectively. (If no such $q$ exists, the statement follows from Corollary~\ref{cor:one side always certified} and for only one such $q$, the proof is essentially the same as here.)

Assume towards a contradiction that $\Lambda_r(f,\can_r(A)) > \lambda = \frac{n}{2} + 3\tau + 1$. Then the shortest occurrence that corresponds to the side $f$ has $\Lambda_r$-measure $> \frac{n}{2} + 3\tau + 1 - 2\epsilon = \frac{n}{2} - \tau - 1$. Since the winner sides for $q_1, q_2$ are certified by Lemma~\ref{rem:winners are certified}, the shortest occurrence that corresponds to the complement of $f$ has $\Lambda_r$-measure $< n - \Lambda_r(\widetilde{f}) + 2\tau = \frac{n}{2} - \tau - 1$ contradicting our assumption that $f$ is the winner side.
\end{proof}

\begin{proposition}
\label{can_stability_after_one_turn}
Let $A, C \in \Scan_{\mu, r}$ and let $A\mapsto C$ be the turn of a maximal occurrence $f$ in $A$ of rank~$r$ and $\FracM_r$-measure $\geqslant \tau + 1$. Then $\can_r(A) = \can_r(C)$.
\end{proposition}

\begin{proof}
Since  $A, C \in \Scan_{\mu, r}$, we have $5\tau + 3 \leq \FracM_r(f) \leq \mu$ by Lemma~\ref{lem:canonical implies turn restrictions}. Hence  the turn $A\mapsto C$ is of Type 2 with inverse turn $C\mapsto A$ of the maximal occurrence $\hat{g}$ where $g$ is the complement of $f$. By symmetry we also have  $5\tau + 3 \leq \FracM_r(\hat{g}) \leq \mu$. Let $u$ be a maximal occurrence in $A$ with $\Lambda_r(u) \geqslant 5\tau + 3$. Assume that at least one of $\Lambda_r(u, A), \Lambda_r(u, C) \leqslant \frac{n}{2} - 5\tau - 2$ or at least one of $\Lambda_r(u, A), \Lambda_r(u, C) \geqslant \frac{n}{2} + 5\tau + 2$. Then as in Lemma~\ref{too_small_and_too_big_side} we obtain that both in $A$ and $C$ the winner side is $u$ or $v$, respectively. So we now assume that this does not happen.

Assume that there exists a maximal occurrence $q$ to the right of $u$ essentially non-isolated from $u$ with $\Lambda_r(q) \geqslant 5\tau + 3$ and the complement $z$. We need to show that in $A$ and in $C$  the certified sides of $q$ are the same. This is clear if $f$ is to the left of $u$, so we assume that either $f = u$, or $f$ is to the right of $u$. Then the result follows from Corollaries~\ref{cor:turn outside} and~\ref{cor:turns inside}. If such occurrence $q$ does not exist in $A$ but exists in $C$, then we consider the inverse turn $C\mapsto A$ and the result follows.
\end{proof}

\section{An auxilliary group structure}
\label{scan_equiv_section}

In order to prove IH~\ref{IH can of equivalent words} we will need to show that for $A, B\in\Can_{-1}$ with $A\equiv B\modRr$ we have $\can_r(A)=\can_r(B)$. We begin with showing this for $A, B\in\Scan_{\lambda + \tau + \varepsilon,r}$ by introducing a group structure on equivalence classes on $\Scan_{\lambda + \tau + \varepsilon,r}$ where (as always) $\lambda =  \frac{n}{2}+3\tau+1$. We will show that the equivalence relation coincides with  equality in $\Fr/ \llangle \Rel_0, \ldots, \Rel_r\rrangle$. Using this equivalence relation, we will then define the canonical form of rank~$r$ of arbitrary words in Section~\ref{can_of_arbitrary_words_section}.

\begin{definition}
\label{semican_equiv_def}
For $A_1, A_2 \in \Scan_{\kappa + \varepsilon, r}$ we define $A_1 \sim_{\kappa, r} A_2$ if and only if there exists a (possibly empty) sequence of turns of rank~$r$ of $\FracM_r$-measure $\geq \tau$ 
\begin{equation*}
A_1 = C_1 \longmapsto C_2 \longmapsto \ldots \longmapsto C_t = A_2
\end{equation*}
with $C_i\in\Scan_{\kappa +\epsilon,r}$ for $2 \leq i \leq t - 1$. The same $\sim_{\kappa, r}$ is defined also for $A_1, A_2 \in \Scan_{\kappa, r}$.
\end{definition}

\noindent
In this section we will use the relation $\sim_{\kappa, r}$ with $\kappa = \lambda + \tau = \fracn + 4\tau + 1$.

\begin{remark}
\label{min_measure_in_equivalence_sequence}
Since all $C_i$ in this sequence belong to $\Scan_{\lambda + \tau + \epsilon,r}$, by Lemma~\ref{lem:canonical implies turn restrictions} any occurrence $u$ which is turned in the sequence satisfies $\frac{n}{2} - 8\tau - 2 < \Lambda_r(u) \leq \lambda + \tau + \epsilon$.  Since $n-(\lambda + \tau + \epsilon) = \frac{n}{2} - 6\tau - 2 > 12\tau$, such a turn is of Type~2 by Lemma~\ref{turn_res_complement} (ii) with remainder $v'$ of $\Lambda_r$-measure $> 10\tau$ and hence has an inverse turn of $\Lambda_r$-measure $> 10\tau$ by  Lemma~\ref{inverse_turn_formal_existence}. Thus, if
\[A_1 = C_1 \mapsto C_2 \mapsto \ldots \mapsto C_t = A_2\]
is a sequence witnessing  $A_1 \sim_{\lambda + \tau, r} A_2$, then  $A_2 \sim_{\lambda + \tau, r} A_1$ is witnessed by its inverse sequence
\[A_2 = C_t \mapsto C_{t - 1} \mapsto \ldots \mapsto C_1 = A_1.\]
\end{remark}

By this observation we obtain:

\begin{corollary}
\label{semican_equiv}
The relation $\sim_{\lambda + \tau, r}$ is an equivalence relation on $\Scan_{\lambda + \tau + \varepsilon, r}$ and on $\Scan_{\lambda + \tau, r}$ with finite equivalence classes. Moreover every equivalence class in $\Scan_{\lambda + \tau r}$ has a representative in $\Scan_{\lambda, r}$.
\end{corollary}
\begin{proof}
By Lemma~\ref{lem:canonical implies turn restrictions} every turned occurrence in the sequence witnessing $A_1 \sim_{\lambda + \tau, r} A_2$ has $\Lambda_r$-measure $> \frac{n}{2} - 8\tau - 2 \geqslant 9\tau + 5$. Hence Corollary~\ref{sequence_of_turns_to_consecutive_turns} implies that there exists a stable sequence of occurrences in $A_1$ such that their turns give $A_2$. Therefore the members of an equivalence class of $A_1$ correspond to choices of sides in maximal occurrences $u$ in $A$ such that $\Lambda_r(u) \geqslant \tau + 1$, and there are only finitely many of these.

Lemma~\ref{scan_exists_simple_case} implies that every equivalence class in $\Scan_{\lambda + \tau, r}$ has a representative in $\Scan_{\lambda, r}$.
\end{proof}

The equivalence class of a word $A \in \Scan_{\lambda + \tau + \varepsilon, r}$ is denoted by $[A]$.
 Recall that in Section~\ref{can_of_semican_section} we defined $\can_r$ for $\mu$-semicanonical words, where $\mu = n -  (8\tau + 3)$. Since  $\lambda +  \tau + \varepsilon = \frac{n}{2} + 6\tau + 2 \leqslant \mu$, Proposition~\ref{can_stability_after_one_turn}  now implies:

\begin{corollary}
\label{scan_equiv_same_can}
If $A_1, A_2 \in \Scan_{\lambda + \tau + \varepsilon, r}$ are $\sim_{\lambda + \tau, r}$-equivalent, then $\can_r(A_1) = \can_r(A_2)$.
\end{corollary}

We will now define an (auxilliary) group structure on $\Scan_{\lambda + \tau, r}/\sim_{\lambda + \tau,r}$ and establish that 
for $A_1, A_2 \in \Scan_{\lambda + \tau, r}$ we have $A_1 \sim_{\lambda + \tau, r} A_2$ if and only if $A_1$ and $A_2$ represent the same element in $\Fr / \llangle \Rel_0, \ldots, \Rel_r\rrangle$. Since we were not able to show directly that different $\lambda + \tau$-semicanonical forms of a given word are $\sim_{\lambda + \tau, r}$-equivalent, we need this group structure to show that we obtain a well-defined canonical form of rank~$r$ of an arbitrary word using an arbitrary $\lambda + \tau$-semicanonical form for it.

We first define the multiplication $\times_{\lambda + \tau,r}$ on $\Scan_{\lambda + \tau,r}$. For technical reasons we define it on larger set $\Scan_{\lambda + \tau + \epsilon,r}$.
\begin{definition}
\label{semican_mult_def}
For $A, C \in \Scan_{\lambda + \tau + \epsilon, r}$, let $Z'=\prod_{\lambda+6\tau+6,r}(A\cdot C)\in \Scan_{\lambda +6\tau+2, r}$ and let $Z''$  be a $\lambda + \tau$-semicanonical form of rank~$r$ of $Z'$ obtained by turns of $\Lambda_r$-measure $\geq \fracn + 2\tau$ (such sequence of turns exists by Lemma~\ref{scan_exists_simple_case}). Define $A\times_{\lambda + \tau, r} C = [Z'']\in \Scan_{\lambda + \tau, r}/ \sim_{\lambda + \tau, r}$.
\end{definition}

\begin{remark}
\label{rem:scan mult}
In Definition~\ref{semican_mult_def} we define the multiplication~$\times_{\lambda + \tau, r}$ in two steps: for $A, C \in \Scan_{\lambda + \tau + \varepsilon, r}$ we find $Z=\can_{r - 1}(A\cdot C)$,  compute a \emph{specific} $(\lambda+6\tau + 2)$-semicanonical form $Z'$ of $Z$ and then find a $\lambda + \tau$-semicanonical form $Z''$  of $Z'$.

Note that every word in $\Scan_{\lambda + 6\tau + 2, r}$ has a $\lambda + \tau$-semicanonical form by Corollary~\ref{scan_exists}. Since $\lambda + 6\tau + 2 = \fracn + 9\tau + 3 \leqslant \mu = n - 8\tau - 3$, Lemma~\ref{two_sequences_to_consecutive_turns} implies that the resulting equivalence class does not depend on the particular $\lambda + \tau$-semicanonical form of $Z'$ as long as the descent is obtained from turns of $\Lambda_r$-measure $\geq 9\tau + 5$.
\end{remark}

We will just  write $\times$ for $\times_{\lambda + \tau, r}$ if the parameters are clear from the context. We emphasis that $A\times_{\lambda + \tau, r} B$ is not a single word, but an equivalence class in $\Scan_{\lambda + \tau, r} / \sim_{\lambda + \tau, r}$.

\begin{remark}\label{rem:group Scan}
It follows directly from the definition that for $A,C \in \Scan_{\lambda + \tau, r}$ (rather than $\Scan_{\lambda + \tau + \epsilon, r}$)  we have
\begin{enumerate}
\item $A\times_{\lambda + \tau, r} 1=[A]$;
\item $ A\times_{\lambda +\tau, r} A\inv=[1]$; and
\item $(A\times_{\lambda + \tau, r} C)\inv=[C\inv\times_{\lambda + \tau, r}A\inv]$.
\end{enumerate}
\end{remark}

We next show that the multiplication factors through $\sim_{\lambda + \tau, r}$:
\begin{proposition}\label{prop:times on classes}
For $[A], [C]\in\Scan_{\lambda + \tau, r}/\sim_{\lambda + \tau, r}$ the multiplication $[A]\times  [C]=[A\times C]$ is well-defined.
\end{proposition}

The crucial step for the proof is contained in the following lemma:  
\begin{lemma}\label{lem:times with turn}
Let $A, B\in\Scan_{\lambda + \tau + \epsilon, r}$. If $A\mapsto B$ is a turn of a maximal occurrence $u$ in $A$ with $\Lambda_r(u)\geq \tau$, then for any $C\in\Scan_{\lambda + \tau, r}$ we have $A\times C = B\times C$.
\end{lemma}
\begin{proof}
By Remark~\ref{min_measure_in_equivalence_sequence} the inverse turn $B\mapsto A$ is defined. Therefore by Corollary~\ref{cor:one side always certified} we can assume that $\Lambda(u) > \frac{n}{2} - \tau$. Hence by Lemma~\ref{mergings_and_cancellations1} we may assume that the turn of $u$ commutes with the product with $C$. That is, we can assume that $\can_{r-1}(A\cdot C) \mapsto \can_{r-1}(B\cdot C)$ is a turn of measure $\geqslant \tau + 1$ of maximal occurrence $\widetilde{u}$.

Now consider the sequences of seam turns
\begin{align*}
X_0 &=\can_{r-1}(A\cdot C)\mapsto X_1\mapsto\ldots\mapsto X_m & =\prod_{\lambda +6\tau+2,r}(A\cdot C) & \quad \mbox{ and }\\
Y_0&=\can_{r-1}(B\cdot C)\mapsto Y_1\mapsto\ldots\mapsto Y_k & =\prod_{\lambda +6\tau+2,r}(B\cdot C).
\end{align*}
First assume that $\widetilde{u}$ is isolated from the seam occurrences in $X_i$ for all $0 \leqslant i \leqslant m$. Then $\Lambda_r(\widetilde{u},X_m) = \Lambda_r(\widetilde{u}, X_0) > \fracn - \tau$. Recall from Lemma~\ref{lem:mu' semicanonical} that all but possibly the last turns in each sequence have $\Lambda_r$-measure $>n-(3\tau+1)$ and the last turn has $\Lambda_r$-measure $>\fracn + 8\tau + 3$. At each step the turn of $\widetilde{u}$ commutes with the seam turn by Lemma~\ref{lem:close neighbour turns are independent} and hence we see that $m=k$ and we obtain $Y_k$ from $X_k$ by turning the occurrence corresponding to $\widetilde{u}$ in $X_k$. Since $X_, Y_k\in\Scan_{\lambda+6\tau+2,r}$ and $\lambda +6\tau+2 = \fracn + 9\tau + 3 \leqslant n - 8\tau - 3$, by Lemma~\ref{two_sequences_to_consecutive_turns} $A\times C = B\times C$.

Now we consider the general case. If $X_0 \mapsto Y_0$ is the seam turn (which is unique by definition), then $X_1 = \can_{r-1}(B\cdot C)$. Since all seam turns are uniquely defined, we obtain that $X_m = Y_k$. So as above the result follows from Lemma~\ref{two_sequences_to_consecutive_turns}.

Assume that $X_0 \mapsto Y_0$ is not a seam turn.  If $u$ or $v$ correspond to the seam occurrences in $X_1$ or $Y_1$, respectively, then by Lemma~\ref{lem:close neighbour turns are independent} we see that $X_2 = Y_1$ or $X_1 = Y_2$. Hence as above $X_m = Y_k$ and the result follows from Lemma~\ref{two_sequences_to_consecutive_turns}.

If $u$ does not correspond to the seam occurrence in $X_1$ and its complement $v$ does not correspond to the seam occurrence in $Y_1$, then $X_1 \mapsto Y_1$ is the turn of the occurrence corresponding to $u$ by Lemma~\ref{lem:close neighbour turns are independent}. By Lemma~\ref{cor:one side always certified} we can assume that it has $\Lambda_r$-measure $> \fracn - \tau$ (otherwise we switch to the inverse turn) and denote the turning occurrence by $u_1$. Repeating the above argument for the turn of $u_1$ until we reach the end of one of the sequences, we obtain the required result.

Assume that the second sequence of seam turns is empty and the first sequence of seam turns is not empty. This means that there are no seam occurrences in $Y_0$, so the occurrence in $Y_0$ that corresponds the the seam occurrence in $X_0$ is not a seam occurrence in $Y_0$. In particular it has $\Lambda_r$-measure $< \lambda + 6\tau + 2$ and the corresponding seam occurrence in $X_1$ has $\Lambda_r$-measure $< \lambda + 6\tau + 2 + \varepsilon < n - (3\tau + 1)$. So the first sequence is of length one and $X_1 \in \Scan_{\lambda + 6\tau + 2, r}$. However, by Remark~\ref{rem:scan mult} we still can turn this occurrence in $Y_0$ and after that take its $\lambda + \tau$-semicanonical form. As above then $X_1 \mapsto Y_1$ is the turn of the occurrence corresponding to $u$. Therefore the result follows from Lemma~\ref{two_sequences_to_consecutive_turns}.
\end{proof}

\begin{proof} (of Proposition~\ref{prop:times on classes})
We denote the operation $\times_{\lambda + \tau, r}$ for short by $\times$.

Let $A_1, A_2, C \in \Scan_{\lambda + \tau, r}$ with $A_1\sim_{\lambda + \tau,r} A_2$. By Remark~\ref{rem:group Scan} (iii) it suffices to prove that $A_1 \times C = A_2\times C$. 
Since $A_1 \sim_{\lambda + \tau, r} A_2$, there exists a sequence of turns of rank~$r$
\begin{equation*}
A_1 = X_1 \longmapsto X_2 \longmapsto \ldots \longmapsto X_t = A_2
\end{equation*}
such that all $X_i \in \Scan_{\lambda + \tau + \epsilon, r}$  for $2 \leq i \leq t - 1$
and all turns have $\Lambda_r$-measure $\geq \tau$.
By Lemma~\ref{lem:times with turn} we have $X_i\times C = X_{i+1}\times C$ and hence $A_1\times C=A_2\times C$.
\end{proof}

\begin{remark}
\label{scan_of_prefix}
Let $C \in \Scan_{\kappa, r}$.

\medskip\noindent
(i) If $C_1$ is a prefix of $C $, then by IH~\ref{can_of_prefix_and_suffix_hyp} we have  $\can_{r - 1}(C_1)=C_1'D$ for a prefix $C_1'$ of $C_1$ and a $\tau$-free side of a canonical triangle $D$, so $\can_{r - 1}(C_1)\in \Scan_{\kappa + \tau, r}$.

\medskip\noindent
(ii) If $z=\can_{r-1}(z)$ is a single letter, then  $\can_{r - 1}(C \cdot z) =C'Dz'$  for a prefix $C'$ of $C$, a $\tau$-free side $D$ of a canonical triangle and $z'\in\{1,z\}$, so $\can_{r - 1}(C \cdot z) \in \Scan_{\kappa+ \tau}$. 
\end{remark}

\begin{lemma}
\label{collect_by_letters}
Let $C= z_1\cdots z_t$ be a reduced word such that for every initial segment $C_s=z_1\cdots z_s$
we have $\can_{r-1}(C_s) \in \Scan_{\lambda+\tau, r}$. Then
\[
(\ldots (z_1\times z_2)\times z_3 )\times\ldots )\times z_t=[\can_{r-1}(C)]
\]
where $\times$ is $\times_{\lambda + \tau, r}$ and $z_i$, $1 \leqslant i \leqslant s$ are single letters.
\end{lemma}

\noindent
Note that by Remark~\ref{scan_of_prefix} this applies in particular if $C=z_1\cdots z_t\in\Scan_{\lambda, r}$.
\begin{proof}
We prove $(\ldots(z_1\times z_2)\times z_3 )\times\ldots )\times z_s = [\can_{r - 1}(z_1\cdots z_s)]$ for all $s \leq t$ by induction on $s$.

For $s=1$ there is nothing to prove, so assume inductively
\[
(\ldots(z_1\times z_2)\times z_3 )\times\ldots\times z_{s-1} )\times z_s = \can_{r - 1}(z_1\cdots z_{s - 1})\times z_s.
\]
By Corollary~\ref{can_of_can_product_hyp} and Remark~\ref{scan_of_prefix} we have \[\can_{r - 1}(\can_{r - 1}(z_1\cdots z_{s - 1}) \cdot z_s)=\can_{r-1}(z_1\ldots z_s)\in\Scan_{\lambda + \tau,r}.\]
So $\can_{r - 1}(z_1\cdots z_{s - 1})\times z_s =  [\can_{r - 1}(z_1\cdots z_s)]
$
and for $s = t$ we obtain the required result.
\end{proof}

In order to establish that $\Scan_{\lambda+\tau, r}/ \sim_{\lambda + \tau, r}$  is a group with respect to  $\times_{\lambda + \tau, r}$, we now show that the multiplication is associative.

\begin{lemma}
\label{associativity_in_scan}
Let $A, B, C \in \Scan_{\lambda+\tau, r}$. Then $(A\times B)\times C = A\times (B\times C)$, where $\times$ is $\times_{\lambda + \tau, r}$.
\end{lemma}
\begin{proof}
First we prove the statement for $C = z$ a single letter. By Corollary~\ref{semican_equiv}, we can assume that $B\in \Scan_{\lambda, r}$. Then by Remark~\ref{scan_of_prefix}, $B\cdot_{r - 1} z \in \Scan_{\lambda + \tau, r}$, therefore $B\times C = [B\cdot_{r - 1} z]$. By definition $A\times B$ is calculated using a sequence of turns
\[A\cdot_{r - 1} B = X_0 \mapsto\ldots\mapsto X_m \in \Scan_{\lambda+\tau,r},\]
and $A\times (B\times z) = A\times (B\cdot_{r - 1} z)$ is calculated using a sequence of turns
\[ A \cdot_{r - 1} (B\cdot_{r - 1} z) = Y_0 \mapsto\ldots\mapsto Y_k \in \Scan_{\lambda+\tau,r}.\]
Recall that all turns in these sequences are of $\Lambda_r$-measure $> \frac{n}{2} + 2\tau$. When we multiply the first sequence by $z$,  Lemma~\ref{mergings_and_cancellations1} implies that $X_i \cdot_{r - 1} z \mapsto X_{i + 1} \cdot_{r - 1} z$ are turns of rank~$r$ and by Remark~\ref{scan_of_prefix} they have $\Lambda_r$-measure $> \frac{n}{2} + \tau$.

If in at least one of the sequences there are no seam turns, then all maximal occurrences in $X_0, Y_0$ are of $\Lambda_r$-measure $< \lambda + \tau + 3\tau + 1 + \tau$ by Remark~\ref{scan_of_prefix}, so $X_0, Y_0 \in \Scan_{\frac{n}{2} + 8\tau + 2, r}$. Hence the result follows from Lemma~\ref{two_sequences_to_consecutive_turns}. Assume that the first sequence has seam turns and let $X_i \mapsto X_{i + 1}$ be a seam turn and $X_i \cdot_{r - 1} z = Y_i$. Then either the corresponding turn $X_i \cdot_{r - 1} z \mapsto X_{i + 1} \cdot_{r - 1} z$ is the seam turn in $Y_i$ so $X_{i + 1} \cdot_{r - 1} z = Y_{i + 1}$, or the corresponding occurrence in $Y_i$ has $\Lambda_r$-measure $< \lambda + \tau + 3\tau + 1$. In the latter case $X_i \in \Scan_{\frac{n}{2} + 8\tau + 2, r}$ by Remark~\ref{scan_of_prefix} so the result again follows from Lemma~\ref{two_sequences_to_consecutive_turns}. In the first case we repeat the argument until we exhaust all seam turns in one the sequences.

Now let $C = z_1\cdots z_s \in \Scan_{\lambda, r}$, $C_t$ be a prefix of $C$ of length $t$ and $Z_t = \can_{r - 1}(C_t)$. Then Lemma~\ref{collect_by_letters} implies that $Z_{t - 1}\times z_t = [Z_t]$. Therefore using induction on $t$ we have
\begin{align*}
(A\times B)\times Z_t &= (A\times B) \times (Z_{t - 1}\times z_t) = ((A\times B)\times Z_{t - 1})\times \ z_t =\\
&= (A\times (B \times Z_{t - 1})) \times z_t = A \times ((B\times Z_{t - 1})\times z_t) =\\
&= A \times (B \times (Z_{t - 1} \times z_t)) = A\times (B \times Z_t),
\end{align*}
and for $t = s$ we obtain the final result.
\end{proof}
By Proposition~\ref{prop:times on classes} and Lemma~\ref{associativity_in_scan}, we now have
\begin{corollary} $(\Scan_{\lambda + \tau, r}/ \sim_{\lambda + \tau, r}, \times_{\lambda + \tau, r})$ is a group.
\end{corollary}

The main statement of this section is
\begin{proposition}\label{prop:main isomorphism}
For $A_1,A_2\in\Scan_{\lambda+\tau,r}$ the following are equivalent:
\begin{enumerate}
\item $A_1\sim_{\lambda + \tau,r} A_2$ ;
\item $A_1$ and $A_2$ represent the same element in $\Fr/\llangle \Rel_0, \ldots, \Rel_r\rrangle$;
\item $\can_r(A_1)=\can_r(A_2)$.
\end{enumerate}
\end{proposition}

\medskip\noindent
For the proof we  define an epimorphism 
\[\phi : \Fr=\langle  x_1, \ldots, x_m\rangle  \longrightarrow \Scan_{\lambda+\tau, r}/ \sim_{\lambda + \tau, r}\colon\quad x_i\mapsto  [x_i].\] 
By the universal property of free groups  $\phi$  is well-defined.
 For a reduced word $C=z_1\cdots z_t \in \Can_0$  we then have
\[
\phi(C) = \phi(z_1\cdots z_t) = \phi(z_1)\times \ldots \times \phi(z_t) = z_1\times\ldots \times z_t.
\]
\begin{remark}\label{rem:phi surj}
By Lemma~\ref{collect_by_letters} we have $\phi(C)=[C]$ for $C\in\Scan_{\lambda,r}$. Therefore  $\phi$ is surjective by Corollary~\ref{semican_equiv}.
\end{remark}

\begin{remark}\label{rem:turn equivalent}
We will repeatedly make use of the observation that,  by the very definition of a turn, if $A\mapsto B$ is a turn of rank~$i \leq r$, then $A\equiv B\mod  \llangle \Rel_0, \ldots, \Rel_i\rrangle$.
\end{remark}

\begin{lemma}\label{lem:can for relator halfs}
Let $R=a^n=z_1\cdots z_t\in \Rel_i$, $i\leq r$. Write $R_s=z_1\cdots z_s$ and let $V_s$ denote its complement. Then for $s\leq t$ the following holds:
\begin{enumerate}
\item If $\Lambda_i(R_s)\geq 3\tau + 1$, then $\can_{i-1}(R_s)=DR_s'E\mapsto \can_{i-1}(V_s)$  is a turn of rank~$i$ of the maximal prolongation of $R_s'$. 
\item If $i<r$, then \begin{equation*}
\can_{r - 1}(R_s) = \can_{r - 1}(V_s) =
\begin{cases}
\can_{i - 1}(R_s) &\textit{ if } \can_{i - 1}(R_s) \in \Can_i,\\
\can_{i - 1}(V_s) &\textit{ if } \can_{i - 1}(V_s) \in \Can_i.
\end{cases}
\end{equation*}
Furthermore $\can_{r - 1}(R_s)$ is $6$-semicanonical of rank~$r$.

\end{enumerate} 
\end{lemma}
\begin{proof}
Part 1: By IH \ref{can_of_prefix_and_suffix_hyp} we have $\can_{i-1}(R_s)=DR'_sE$ for an appropriate subword 
$R'_s$ of $R_s$ and $\tau$-free of rank~$i$ words $D, E$. If $\Lambda_i(R_s)\geq 3\tau+1$, then $\Lambda_i(R_s')\geq \tau+1$ and so the maximal prolongation of $R_s'$ is a maximal occurrence. By Remark~\ref{can_preserves_coset_hyp}, IH~\ref{IH can of equivalent words} and Remark~\ref{rem:cyclic shift}, we have 
\begin{equation*}
\can_{i - 1}(D\cdot \widehat{a}^{-n}\cdot R_s'E) = \can_{i - 1}(a^{-n}\cdot R_s) = \can_{i - 1}(V_s)
\end{equation*}
 where $\widehat{a}$ is the corresponding cyclic shift of $a$. 

Part 2:
If $\can_{i-1}(R_s)\notin\Can_i$, then only the maximal prolongation of $R'_s$ can have $\Lambda_i$-measure $>3\tau+1$. Hence by Part 1, turning $R_s'$ in $DR'_sE$ yields $\can_i(R_s)=\can_{i - 1}(V_s)$. The word $\can_{i - 1}(V_s)$ is $6$-free of rank $>i$ by IH~\ref{rel_common_parts_in_diff_ranks_hyp} and hence $\can_j(R_s)=\can_i(R_s)$ for $j>i$ by IH~\ref{IH immediate}.
\end{proof}

\begin{corollary}\label{x for relators}
Let $R=a^n=z_1\cdots z_t\in \Rel_i, i\leq r, R_s=z_1\cdots z_s$. Then for $s\leq \half$ we have 
\[z_1\times\ldots\times z_s=[\can_{r-1}(R_s)].\]
\end{corollary}
\begin{proof}
Since $\can_{r-1}(R_s)$ is $6$-semicanonical of rank~$r$ for $i < r$ by Lemma~\ref{lem:can for relator halfs} and $\can_{r-1}(R_s)$ is $\fracn + 2\tau + 1$-semicanonical of rank~$r$ for $i = r$, the result follows immediately from Lemma~\ref{collect_by_letters}.
\end{proof}

For the proof of Proposition~\ref{prop:main isomorphism} we need
\begin{lemma}\label{lem:kernel phi}
$\llangle \Rel_0, \ldots, \Rel_r\rrangle\leq \ker\phi$.\footnote{In fact, one can show that $ \ker\phi=\llangle \Rel_0, \ldots, \Rel_r\rrangle$.}
\end{lemma}

\begin{proof}
Let $R=a^n=z_1\cdots z_t\in \Rel_i, i\leq r, R_s=z_1\cdots z_s$ and  $T_s=z_{s +1}\cdot\ldots\cdot z_t$.
By Corollary~\ref{x for relators} we have
\[Z_1=z_1\times\ldots\times z_\half= [\can_{r-1}(R_\half)]\] and similarly, by considering the appropriate cyclic shift,
\[Z_2=z_{\half +1}\times\ldots\times z_t=[\can_{r-1}(T_\half)].\]
Then, by Definition~\ref{semican_mult_def},
$Z_1\times Z_2=z_1\times\ldots\times z_t$ is computed from 
\[\can_{r-1}(Z_1\cdot Z_2)=\can_{r-1}(R_\half\cdot_{r-1}T_\half)=\can_{r-1}(R)\] by first taking seam turns. However, if $R\in\Rel_i, i<r$, then $\can_{r-1}(R)=1$ by Remark~\ref{can_preserves_coset_hyp} and hence in this case no seam turn is necessary and $\varphi(R) = \can_{r-1}(R) = 1$.
In  case $R\in\Rel_r$, we have $\can_{r-1}(R)=DR'E$ for some $\tau$-free of rank~$r$ $D, E$, so by Lemma~\ref{lem:product of lambda canonical} the maximal prolongation of $R'$ is a seam occurrence and the result of the seam turn
is equal to $1$ by Lemma~\ref{lem:can for relator halfs} part~1.
\end{proof}

\begin{proof}(of Proposition~\ref{prop:main isomorphism})
Let $A_1,A_2\in\Scan_{\lambda + \tau,r}$.

\medskip\noindent
1.$\Rightarrow$ 2. and 1.$\Rightarrow$ 3.:  If $A_1\sim_{\lambda + \tau,r} A_2$, 
there is a sequence of turns of rank~$r$ such that $A_1=X_0\mapsto\ldots\mapsto X_k=A_2$. Thus we have  $A_1\equiv A_2 \mod \llangle\Rel_0, \ldots, \Rel_r\rrangle$ by Remark~\ref{rem:turn equivalent} and  $\can_r(A_1) = \can_r(A_2)$ by Corollary~\ref{scan_equiv_same_can}.
 
\medskip\noindent
2.$\Rightarrow$ 1.: Suppose  $A_1\equiv A_2\mod  \llangle \Rel_0, \ldots, \Rel_r\rrangle$. By Corollary~\ref{semican_equiv} we may assume $A_1, A_2\in\Scan_{\lambda,r}$. Hence by Remark~\ref{rem:phi surj} and Lemma~\ref{lem:kernel phi} we have $[A_1]=\phi(A_1)=\phi(A_2)=[A_2]$ and  hence  $A_1\sim_{\lambda + \tau,r} A_2$.

\medskip\noindent
3.$\Rightarrow$ 2.:  Suppose $\can_r(A_1) = \can_r(A_2)$. Since $\can_r(A_i)$ is obtained from $A_i$, $i=1,2$, by turns of rank~$r$,  we have 
\[A_1\equiv \can_r(A_1)=\can_r(A_2)\equiv A_2\quad \mod \llangle\Rel_0, \ldots, \Rel_r\rrangle\]
\end{proof}

\begin{corollary}
\label{different_scan_equiv}
 If $A, B$ are $\lambda$-semicanonical forms of rank~$r$ of some $C\in\Can_{r-1}$, then 
$\can_r(A) = \can_r(B)$.
\end{corollary}
\begin{proof}
If $A, B$ arise from $C$ by turns of rank~$r$, then
$A\equiv C\equiv B \mod \llangle\Rel_0, \ldots, \Rel_r\rrangle$, so the claim follows Proposition~\ref{prop:main isomorphism}.
\end{proof}

\subsection{Canonical form of rank~$r$ of arbitrary words}
\label{can_of_arbitrary_words_section}

Recall that since
\[
\mu= n - (8\tau + 3)> \lambda = \frac{n}{2} + 3\tau + 1
\]
for our choice of the exponent $n$, every $\lambda$-semi\-canonical form of rank~$r$ is also $\mu$-semicanonical.

\begin{definition}[canonical form of rank~$r$]
\label{can_r_def}
For $A\in\Can_{r - 1}$  we define \emph{the canonical form of rank~$r$ of $A$}, $\can_r(A)$, in two steps as follows:
\begin{enumerate}
\item
choose  a $\lambda$-semicanonical form $A'$ of rank~$r$ for $A$. 
\item
put $\can_r(A)=\can_r(A')$ as defined in
 Section~\ref{can_of_semican_section}. 
\end{enumerate}
For $A\in \Can_{-1}$ we define $\can_r(A)=\can_r(\can_{r - 1}(\ldots\can_0(A)\ldots))$.
\end{definition}
Note that $A\equiv A' \mod \llangle \Rel_0, \ldots, \Rel_r\rrangle$ by construction and $\can_r(A)$ does not depend on $A'$  by Corollary~\ref{different_scan_equiv}.

\begin{lemma}[$\mu = n - (8\tau + 3)$]
\label{can_r_def_agree}
If $A \in \Scan_{\mu - \varepsilon_2, r}$, then Definition~\ref{can_r_def} and choosing winner sides from Section~\ref{can_of_semican_section} give the same $\can_r(A)$.
\end{lemma}
\begin{proof}
By Lemma~\ref{scan_exists_simple_case} there exists a $\lambda$-semicanonical form $A'$ of $A$ that is obtained from $A$ by a sequence of turns with all intermediate words $\mu$-semicanonical. Thus the result follows from Proposition~\ref{can_stability_after_one_turn}.
\end{proof}

\begin{lemma}
\label{can_preserves_coset}
Let $A \in \Can_{-1}$. Then $A \equiv \can_r(A) \mod \llangle \Rel_0, \ldots, \Rel_r\rrangle$.
\end{lemma}
\begin{proof}
Let $\can_{r - 1}(A) = B$. By definition, $\can_r(A) = \can_r(B)$. By Remark~\ref{can_preserves_coset_hyp}, we have $B \equiv A \mod \llangle  \Rel_0, \ldots, \Rel_{r - 1}\rrangle $. By definition, $\can_r(B)$ is obtained from $B$ by a sequence of turns of rank~$r$, therefore, $B \equiv \can_r(B) \mod \llangle  \Rel_0, \ldots, \Rel_r\rrangle $. Thus,
\[
A \equiv  \can_r(B) = \can_r(A) \mod \llangle  \Rel_0, \ldots, \Rel_r\rrangle .
\]
\end{proof}

We now verify some further induction hypotheses:
\begin{proposition}
\label{can_of_equivalent_words}
IH~\ref{IH can of equivalent words} and IH~\ref{IH idempotent} hold for rank~$r$:
for $A, B \in \Can_{-1}$ we have $\can_r(A) = \can_r(B)$ if and only if $A\equiv B \mod \llangle  \Rel_0, \ldots, \Rel_r\rrangle $.  In particular, $A=\can_r(A)$ for $A\in\Can_r$.

\end{proposition}
\begin{proof}
If $\can_r(A) = \can_r(B) = C$, then Lemma~\ref{can_preserves_coset} implies that $A \equiv C \equiv B \mod \llangle  \Rel_0, \ldots, \Rel_r\rrangle $.

Conversely, assume that $A\equiv B \mod \llangle  \Rel_0, \ldots, \Rel_r\rrangle $. Let $A', B'$ be  $\lambda$-semicanonical forms of rank~$r$ of $A, B$, respectively. Then
\begin{equation*}
A' \equiv A \equiv B \equiv B' \mod \llangle  \Rel_0, \ldots, \Rel_r\rrangle .
\end{equation*}
Hence, Corollary~\ref{scan_equiv_same_can} implies that $\can_r(A) = \can_r(B)$.

If $A=\can_r(B)\in\Can_r$, then $A\equiv B \mod \llangle  \Rel_0, \ldots, \Rel_r\rrangle $ and hence $\can_r(A)=\can_r(B)=A$ by the previous.
\end{proof}

\begin{proposition}
\label{can_elem_prop}
Inductive Hypotheses~\ref{IH can unique}-- \ref{IH can of equivalent words}  hold for $\can_r$.
\end{proposition}
\begin{proof}
IH~\ref{IH can unique}--\ref{IH Rel} follow trivially from the definition of $\can_r$.

IH~\ref{IH immediate} follows from Lemma~\ref{too_small_and_too_big_side}.

IH~\ref{IH small cancellation} and IH~\ref{rel_common_parts_in_diff_ranks_hyp} are proved in Corollary~\ref{ind_hyp_holds_for_rel}.

IH~\ref{IH idempotent} and \ref{IH can of equivalent words} are proved in Proposition~\ref{can_of_equivalent_words}.

IH~\ref{IH can inverse} follows from IH~\ref{IH can inverse} for rank~$r - 1$, and the decision process in Section~\ref{can_of_semican_section}.

\end{proof}
\begin{remark}
\label{rem:winner sides idempotent}
If $A \in \Can_r \subseteq \Scan_{\lambda, r}$, then for every maximal occurrence $u$ of rank~$r$ in $A$  it follows from IH~\ref{IH idempotent} that $u$ is the winner side in the process from Section~\ref{can_of_semican_section}.
\end{remark}

\begin{lemma}
\label{scan_no_collapses}
If $A \in \Scan_{\lambda+\tau, r}\setminus\{1\}$ for $\lambda = \fracn=3\tau+1$, then  $A\notin \llangle \Rel_0, \ldots, \Rel_r\rrangle$.
\end{lemma}
\begin{proof} 
If $A \equiv 1 \mod \llangle  \Rel_0, \ldots, \Rel_r\rrangle $, then $A \sim_{\lambda, r} 1$ by Proposition~\ref{prop:main isomorphism}. However, this is not possible, because the equivalence class of $1$ consists only of $1$ itself.
\end{proof}

\begin{corollary}
IH~\ref{can_non_trivial_subwords_hyp} holds for rank~$r$.
\end{corollary}
\begin{proof}
Let  $U$ be a subword of $A \in \Can_r$. Then $\can_{r - 1}(U) = DU'E$, where $D$ and $E$ are sides of canonical triangles of rank~$r - 1$. Since $D$ and $E$ are $\tau$-free of rank~$r$, it follows from Lemma~\ref{can subset scan} that $\can_{r - 1}(U)$ is $\left(\frac{n}{2} + 5\tau + 1\right)$-semicanonical. 
Using Lemma~\ref{scan_exists_simple_case} we find a $\lambda$-semicanonical form $U_1$  of $\can_{r - 1}(U)$. In particular the lemma implies that $U_1 \neq 1$. Then byLemma~\ref{scan_no_collapses} $U_1 \notin \llangle \Rel_0, \ldots, \Rel_r \rrangle$, and neither is $U$. 
\end{proof}

\section{Power subwords in canonical words of rank~$r$}
\label{turns_in_periodic_words_section}

We start with the following natural definition:

\begin{definition}
\label{closed_under_shifts}
Let $A = LX^KX_1R \in \Can_{r - 1}$ where $X$ is primitive, $X^n \notin \Rel_r$, and $X_1$ is a prefix of $X$.
If  $u$ is a maximal occurrence of rank~$r$ in $X^KX_1$, a  \emph{periodic shift of $u$ in $X^KX_1$} is a shift of $u$ by $\pm k\abs{X}$ in $X^KX_1$ contained in $X^KX_1$. 

\end{definition}
If $u$ is a maximal occurrence properly contained in $X^KX_1$, then clearly $u$ is also a maximal occurrence of rank~$r$ in $A$ and so are all periodic shifts of $u$ that are properly contained in $X^KX_1$. However, if a periodic shift of $u$ is a prefix or a suffix of $X^KX_1$, it may have a prolongation in $A$.

\begin{remark}
\label{occurrences_in_periodic_word}
Let $x, a$ be primitive, not cyclic shifts of each other and let $K\geq 2$. If $x^K$ contains a maximal occurrence $u$ which is a fractional power of $a$ with $\Lambda_a(u)\geq  2$, then $u$ is a prefix of a cyclic shift of $x^K$. Hence, 
if $|x|\geq|u|$, then clearly a cyclic shift of $x$ contains $u$. Otherwise by Lemma~\ref{common_part_of_powers} we have $|x|<|u|<|x|+|a|$ and hence a cyclic shift of $x$ contains $a^m$ where $m=\lfloor\Lambda_a(u) \rfloor$. In particular, $|u|<2|x|$ and so $u$ is a proper subword of $x^3$. Hence
 there exist $\geqslant K - 2$ different periodic shifts of $u$ in $x^K$.
Note that  if (and only if) $u$ is a prefix or suffix of $x^K$, the periodic shifts of $u$ may have proper prolongations with respect to $a$ in $x^K$.
If there exist precisely $K - 2$ different periodic shifts of $u$ in $x^K$, then $u$ is not a subword of $x^2$ and so $u=u_1xu_2$ where $u_1, u_2$ are nonempty suffix and prefix, respectively, of $x$ with $0<\Lambda_a(u_1)+\Lambda_a(u_2)<1$ and $\Lambda_a(x)>\Lambda_a(u)-1$. In particular, all periodic shifts are proper subwords of~$x^K$ and are maximal occurrences in~$x^K$.
\end{remark}

We state these observations for further applications in the following form:

\begin{corollary}
\label{periodic_shifts_cor1}
Let $A = LX^KR \in \Can_{r - 1}$, where $K \geqslant 3$ and $X$ is primitive, $X^n \notin \Rel_r$, and let $u$ be a maximal occurrence of rank~$r$ in $A$ that is contained in $X^K$, $\FracM_r(u) \geqslant 2$. There exist $\geqslant K - 2$ different periodic shifts of $u$ in $X^K$ that are maximal occurrences of rank~$r$ in $A$. Moreover, there exist precisely $K - 2$ such periodic shifts of $u$ in $X^K$ if and only if $u=u_1Xu_2$ and $0 \leqslant \Lambda_r(u_1)+\Lambda_r(u_2)<1$.
\end{corollary}

Further we use the notations from Definition~\ref{certification_def}.

\begin{lemma}
\label{same_measure_in_domino}
Let $(u_0,\ldots, u_t)$ be an (un-)certification sequence in $A=L\ul u_0\ldots u_t\ur R$ to the right of $u_0$ and let $i, j\in\{1,\ldots, t\}$ with $\Lambda_r(u_i)=\Lambda_r(u_j)$. If $t\in\{i, j\}$ assume that $\lambda_{m(t)}-\epsilon<\Lambda_r(f_t, W)<\lambda_{m(t)}+\epsilon$. Then $f_i=u_i$ if and only if $f_j=u_j$.
\end{lemma}
\begin{proof}
Suppose $f_i=u_i$ and $f_j=v_j$.
Then in $A$ by Condition~2 of Definition~\ref{certification_def}  we have 
\[ \lambda_{m(i)}+2\epsilon>\Lambda_r(u_i)=\Lambda_r(u_j)> \lambda_{m(i)}-2\epsilon\geq \frac{n}{2}+\tau.\]
Write $W_j=L'\ul f_0\ldots f_j u_{j+1}\ldots u_t\ur R'$ for the result of turning the necessary occurrences $u_i, i \leqslant j$.
Then in $W_j$ we have
$\Lambda_r(v_j, W_j)< \frac{n}{2}+3\tau+1\leq\lambda_2-\epsilon$. Thus $j=t$ by Condition 4, contradicting our assumption on $f_t$.
\end{proof}

\begin{remark}\label{no_cubic_sequences_in_sdomino}
Suppose $W=L\ul f_0\ldots f_t\ur R$ is the witness of an (un-)cer\-tification sequence to the right of $u = u_0$ in $A$. Let $i, j\in\{1,\ldots, t - 1\}$. Then $\Lambda_r(f_i, W)=\Lambda_r(f_j, W)$ implies $m(i)=m(j)$ since, by definition, $\lambda_{m(i)} \geqslant \Lambda_r(f_i, W)=\Lambda_r(f_j, W)>\lambda_{m(i)}-\epsilon$ and the intervals $[\lambda_1, \lambda_1-\epsilon)$ and $[\lambda_2, \lambda_2-\epsilon)$ are disjoint. In particular, by the choice of the function $m$,
 there are no subsequences of the form $BBb$ in $\FracM_r(f_1, W), \ldots, \FracM_r(f_{t-1},W)$.
\end{remark}

We will need the following refinement of Lemma~\ref{common_part_of_powers}:
\begin{lemma}\label{lem:tiny guys preparation}
Let $u$ be a fractional power of $B$ with $B^n\in\Rel_r$. Let $C=wM$ be primitive, and assume that $w=a^m a_0=a_0(a_1a_0)^m, m\geq\tau,$ is a maximal occurrence of rank~$r$ in $MwM=MC$, where  $a^n\in\Rel_r$, $a=a_0a_1$ and $C^n\notin \Rel_0\cup\ldots\cup\Rel_r$, $M\neq 1$. 
Then the following holds:
\begin{enumerate}
\item If $w$ is a subword of $u$ and $a$ is not a cyclic shift of $B$, then $m=\tau$ and  a cyclic shift of $B$ is of the form $a^{\tau-1}a_2$ for a prefix $a_2$ of $a$;
\item $wMa^2$ is not a subword of $u$;
\item $(a_1a_0)^2Mw$ is not a subword of $u$.
\end{enumerate} 
\end{lemma}
\begin{proof}
1. If $a^m a_0, m\geq \tau,$ is contained in $u$ we see from Lemma~\ref{common_part_of_powers} that $|a^m a_0|<|B|+|a|$. Hence by Corollary~\ref{rem:no nesting in same rank}, we may assume (after taking a cyclic shift) that $B=a^{\tau-1}a_2$ for a proper prefix $a_2$ of $a$. Hence  $m=\tau$ by Lemma~\ref{common_part_of_powers} and $|B|<|C|$.

2. If $wMa^2$ is a subword of $u$ and $a$ is a cyclic shift of $B$, then $w$ is not a maximal occurrence in $wM$, contradicting to our assumption. So by Part 1, $B=a^{\tau-1}a_2$, hence for $Ca^2 = wMa^2$ to be a subword of $u$, $C$ must be a prefix of $B^K$ for some $K$. Hence we may write $C=B^kM'$ where $k$ is maximal possible and $M'$ is not empty. Then $M'$ is a proper prefix of $B$ and a non-empty suffix of $C$.  Hence $M'$ and $M$ have a common suffix.  Notice that occurrences of $a^2$ in $B^K$ arise only inside the maximal prolongations of $a^\tau$.
If $Ca^2=B^kM'a^2$ is a prefix of $B^K$, then $M'a^2$ is a prefix of $B^{K - k}$. However, this implies that $M'$ has a common suffix with $a$, since $\abs{M'} < \abs{B}$. Then $M$ has a common suffix with $a$ contradicting our assumption that $w$ is maximal in $MwM$.

3. follows from 2. by considering the inverses $w\inv, u\inv$ and $M\inv w\inv M\inv$.
\end{proof} 
 
\begin{lemma}\label{lem:fine periods}
Let $a, B, C$ be primitive and $B=a^sa_1\neq C=a^ta_2$ for $4\leq s\leq t\leq s+1$ and $a_1, a_2$  nontrivial prefixes of $a$. If $D$ is a common prefix of $B^m, C^m, |B|\leq|C|$,  then $|D|<|C|+|a|$.
\end{lemma} 
\begin{proof}
Suppose $|D|\geq |C|+|a|$, then $Ca$ is a prefix of $B^2$. So $Ba$ is a prefix of $C$. Therefore $Ba$ is a common prefix of $B^K$ and $a^K$, and we get a contradiction to Lemma~\ref{common_part_of_powers} applied to $B$ and $a$.
\end{proof}

\begin{remark}
\label{rem:good nesting power}
Let $C\in\Canc_{r - 1}$ is primitive and $C^n\notin\Rel_1\cup\ldots\cup\Rel_r$. Then by Corollary~\ref{lem:nesting_occurrences_rel_r} $C$ cyclically contains $\widehat{a}^{\tau}$ for some $\widehat{a}^n \in \Rel$. One can see that there exist $a$ and $C_0$ cyclic shifts of $\widehat{a}$ and $C_0$, respectively, such that $C_0 = a^m a_0C_1$, $m \geqslant \tau$, $a = a_0a_1$, and either $C_1 = 1$, or $a^m a_0$ is a maximal occurrence of rank~$r$ in $C_1a^m a_0C_1$.
\end{remark}

\begin{lemma}
\label{lem:powers in between}
Let $A = LuMzR \in \Can_{r - 1}$, where $u, w$ are maximal occurrences of rank~$r$ with $\Lambda_r(u), \Lambda_r(z) \geqslant \tau + 1$. Let $C^N$ be a subword of $uMz$ where $C\in\Canc_{r - 1}$ is primitive and $C^n\notin\Rel_1\cup\ldots\cup\Rel_r$. Let $C_0 = a^m a_0 C_1$, $m \geqslant\tau$, $a^n \in \Rel_r$, be a cyclic shift of $C$ as in Remark~\ref{rem:good nesting power}. Then $M$ contains $\geqslant N - 4$ occurrences of $a^m a_0$.
\end{lemma}
\begin{proof}
Let $a^m a_0 = w$ and assume the contrary. Clearly $uMz$ contains $\geqslant N - 1$ occurrences of $w$. By Corollary~\ref{common_part_in_different_ranks} each of $u$ and $z$ contain at most one occurrence of $w$, so $uMz$ contains precisely $N - 1$ occurrence of $w$ and $u$ and $z$ properly contain occurrences of $w$. Then $wC_1a^2$ is a subword of $z$ or $(a_1a_0)^2C_1w$ is a subword of $u$ (since $C^N$ is a subword of $uMz$ and $m > 4$). If $C_1 \neq 1$, this contradicts to Lemma~\ref{lem:tiny guys preparation}, otherwise this contradicts to Lemma~\ref{lem:fine periods}.
\end{proof}

\begin{lemma}
\label{no_5_shifts_in_domino}
Let $A=LDu_0\ul u_1\ldots u_t\ur R\in\Scan_{\mu,r}$ where $(u_0,\ldots, u_t)$ is an (un-)certification sequence to the right of $u_0$ in $A$ and $D$ is $\tau$-free of rank~$r$. Let $C^N$ be a subword of $Du_0\ul u_1\ldots u_t\ur$ where $C\in\Canc_{r - 1}$ is primitive and $C^n\notin\Rel_1\cup\ldots\cup\Rel_r$. If $C^N$ cyclically contains $a^{2\tau}$ with $a^n \in \Rel_r$, then $N \leqslant 5$, otherwise $N\leq 6$.
\end{lemma}
\begin{proof}
First assume that $C^N$  does not contain any of $u_i$, $0 \leqslant i \leqslant t$. Since there is no gap in the sequence $(u_0,\ldots, u_t)$, it follows from Lemma~\ref{lem:powers in between} that $N \leqslant 6$. Assume moreover that $C^N$ cyclically contains $a^{2\tau}$ with $a^n \in \Rel_r$. Then also by Lemma~\ref{lem:powers in between} $N \leqslant 5$.

Hence we may assume that $C^N$ contains some $u_i$. Let $W=L'E\ul f_0 f_1\ldots f_t\ur R'$  be the witness of the certification sequence (where $L'$ is a prefix of $LD$, $E$ is $\tau$-free of rank~$r$ if $f_0=v_0$ and $L'E = LD$ if $f_0=u_0$). Let $i\leq t$ be minimal such that $u_i$ is contained in $C^N$ and all its periodic shifts in $C^N$ are maximal occurrences in $C^N$. Let $C_0$ be a cyclic shift of $C$ such that $u_i$ is a prefix of $C_0^2$. So $C_0^{N-1}C_0'$ is a subword of $Du_0\ul u_1\ldots u_t\ur$, where $C_0'$ is a prefix of $C_0$. By Convention~\ref{simpler notation for occurrences} we write $C_0 = \ul u_i\ldots u_{i+k-1}\ur$ for some $k\geq 1$. Then $C_0' = \ul u_i\ldots u_{i+k-2}\ur$ for $k \geqslant 2$, here $u_{i+k-2}$ in $C_0'$ is a periodic shift of a prefix of $u_{i+k-2}$ with $\Lambda_r$-measure $> 4\tau + 2$. Since $\Lambda_r(u_j) \geqslant 5\tau + 3$ for all $0 \leqslant j \leqslant t$, all periodic shifts of $u_s$ that are maximal occurrences in $C^N$ are equal to some $u_j$ by Lemma~\ref{lem:coincide} (except possibly shifts that are a prefix and a suffix of $C^N$). Thus $u_{i + j + sk}$ are equal to each other for $0 \leqslant s \leqslant N - 1$ with a fixed $0 \leq j \leq k - 3$, and are equal to each other for $0 \leqslant s \leqslant N - 2$ with a fixed $j \in \lbrace k - 2, k - 1\rbrace$. Moreover $i + (N - 1)k + k - 3 \neq t$ and $i + (N - 2)k + k - 1 \neq t$ for $k \geqslant 2$.

From now on we assume additionally that $i \neq 0$. If $C^N$ contains only $u_0$, then clearly $N \leqslant 2$, so the result follows.

If $k\geq 2$, then by Lemma~\ref{same_measure_in_domino} the choices in $W$ for $f_{i + j + sk}$ are the same for $0 \leqslant s \leqslant N - 1$ with a fixed $0 \leq j \leq k - 3$, and and are the same for $0 \leqslant s\leqslant N - 2$ with a fixed $j \in \lbrace k - 2, k - 1\rbrace$. Hence by Lemma~\ref{lem:turns in power words} $\Lambda_r(f_{i + sk}, W)$ are the same for $1\leqslant s\leqslant N - 2$ and $\Lambda_r(f_{i + j + sk}, W)$ are the same for $0 \leqslant s\leqslant N - 3$ with a fixed $1 \leqslant j \leqslant k - 1$ (we put these indices in order to consider cases $k = 2$ and $k \geqslant 3$ simultaneously). So by Remark~\ref{no_cubic_sequences_in_sdomino} $m(i + sk)$ are equal to each other for $1\leqslant s\leqslant N - 2$, and $m(i + j + sk)$ with a fixed $1 \leqslant j \leqslant k - 1$ are equal to each other for $0 \leqslant s \leqslant N - 3$. Since $m$ is $BBb$-free, we must have $N - 2 \leqslant 2$, so $N \leqslant 4$.



If $k=1$, we write $C_0=\ul u_i\ur$, then $u_i = u_{i+1} = \ldots = u_{i + N - 3}$, $i + N - 3\neq t$, and $u_{i + N - 2}$ has a prefix equal to $u_i$. So the choices $f_j$ in the witness $W$ are the same for $u_j$ with $i \leqslant j \leqslant i + N - 3$. Hence by Lemma~\ref{lem:turns in power words} and Remark~\ref{no_cubic_sequences_in_sdomino} the $m(j)$ are the same for $i + 1 \leqslant j \leqslant i + N - 4$. It follows from Corollary~\ref{cor:prefix turn} that the turn of $u_{j + 1}$ have the same influence on $u_j$ for all $i \leqslant j \leqslant i + N - 2$. If $N - 2 \geqslant 4$, then  for one of $u_j$ with $i + 1 \leqslant j \leqslant i + N - 3$ any choice for $u_{j + 1}$ fits for a witness, since $\lambda_1 - \lambda_2 \geqslant \varepsilon_2$. This contradicts Condition~4 of Definition~\ref{certification_def}, so $N - 2 \leqslant 3$ and $N \leqslant 5$.
\end{proof}

In fact the proof of Lemma~\ref{no_5_shifts_in_domino} shows the following
\begin{corollary}
\label{cor:periodic subwords in certification}
Let $A=LDu_0\ul u_1\ldots u_t\ur R\in\Scan_{\mu,r}$ where $(u_0,\ldots, u_t)$ is an (un-)certification sequence to the right of $u_0$ in $A$ and $D$ is $\tau$-free of rank~$r$. Let $C^6$ be a subword of $Du_0\ul u_1\ldots u_t\ur R$ that properly contains some $u_i$  where $C\in\Canc_{r - 1}$ is primitive and $C^n\notin\Rel_1\cup\ldots\cup\Rel_r$. Then $Du_0\ul u_1\ldots u_t\ur$ contains $\leqslant 3$ periodic shifts of $u_i$ that are maximal occurrences of rank~$r$ in $A$ different from $u_0$ and $u_t$.
\end{corollary}

\begin{corollary}
\label{cor:powers in context}
Let $A = LuR \in \Scan_{\mu, r}$, where $u$ is a maximal occurrence of rank~$r$ with $\Lambda_r(u) \geqslant 5\tau + 3$. Let $C \in \Canc_{r - 1}$ be primitive and $C^n\notin\Rel_1\cup\ldots\cup\Rel_r$. Let $C^N$ be a subword of $DuR$, where $D$ is $\tau$-free of rank~$r$ suffix of $L$, $N = 6$ if $C$ cyclically contains $a^{2\tau}$ with $a^n \in \Rel_r$, and $N = 7$ otherwise. Then the prefix of $R$ that ends at the end point of $C^N$ is a right context for $u$.
\end{corollary}
\begin{proof}
Let $(u = u_0, u_1, \ldots, u_t)$ be an (un-)certification sequence from the right of $u$ for side $f_1$ and $A = L\ul u_0, \ldots, u_t\ur R_1$. Then $C^N$ is not contained in $D\ul u_0, \ldots, u_t\ur$ by Lemma~\ref{no_5_shifts_in_domino}, so $u_t$ ends strictly from the left of the end of $C^N$. Denote by $M$ a prefix of $R$ that ends at the end point of $C^N$ and let $R = MR_1$ and consider a word $B = LuMR_2 \in \Scan_{\mu, r}$. First notice that $(u = u_0, u_1, \ldots, u_t)$ is a certification or an un-certification sequence for the side $f_1$ in $B$, because it cannot be extended or shortened by Lemma~\ref{no_5_shifts_in_domino} and Conditions~4 and~5 or~5'.

It remains to show that $(u = u_0, u_1, \ldots, u_t)$ in $B$ cannot change its status from certification sequence for $f_1$ to un-certification sequence and vice versa. Let $Q$ be a suffix of $M$ that starts at the end of $u_t$. It is sufficient to show that $Q$ contains a subword of the form $a^\tau M_1 b^\tau$, $a^n, b^n \in \Rel_r$ (because the only possible problematic case is when $f_t = v_t$ and $\Lambda_r(f_t)$ is different in the witnesses for $A$ and for $B$). If $C^N$ does not contain any $u_i$, this follows from Corollaries~\ref{lem:nesting_occurrences_rel_r} and~\ref{common_part_in_different_ranks}.

If $C^N$ contains only $u_t$, then the result is clear. So we can assume that $C^N$ properly contains $u_t$. Then by Corollary~\ref{cor:periodic subwords in certification} $D\ul u_0, \ldots, u_t\ur$ contains $\leqslant 4$ periodic shifts of $u_t$ different from $u_0$. Hence $Q$ contains a periodic shift of a suffix of $u_t$ with $\Lambda_r$-measure $>5\tau + 1$. This completes the proof.
\end{proof}

\begin{corollary}\label{cor:periodic shift of certification}
Let $A= LC^NR \in \Scan_{\mu, r}$, $N \geqslant \tau$, where $C$ is primitive and $C^n \notin \Rel_1\cup\ldots \cup \Rel_r$. Let $u$ be a maximal occurrence of rank~$r$ in $A$ contained in $C^N$ with $\Lambda_r(u) \geqslant 5\tau + 3$ such that its periodic shifts are maximal occurrences in $A$ (except possibly the first and the last one). Then there exist periodic shifts of $u$ that are contained neither in $LC^6$, nor in $C^6R$. Furthermore the left and right contexts together for these the periodic shifts of $u$ are contained in $C^{13}$ and the (un-)certification sequences are periodic shifts of each other (certification sequences are shifted to certification sequences, un-certification sequences are shifted to un-certification sequences).
\end{corollary}
\begin{proof}
The existence follows from Remark~\ref{occurrences_in_periodic_word}, since $\tau \geqslant 13$. The second part follows directly from Corollary~\ref{cor:powers in context}.
\end{proof}

\begin{corollary}
\label{periodic_can_subword}
Let $A = LC^NR \in \Scan_{\mu, r}$, where $N \geqslant \tau$ is a sufficiently big positive number, $C$ is  primitive and $C^n \notin \Rel_1\cup\ldots \cup \Rel_r$. Then $\can_r(A) = \tilde{L}Y^{N - \gamma}\tilde{R}$, where $C$ and $Y$ are conjugate in rank~$r$, and $\gamma$ does not depend on~$N$.
\end{corollary}
\begin{proof}
By Corollary~\ref{cor:periodic shift of certification} the certification sequences for any maximal occurrence of rank~$r$ in $A$ that is contained in $C^N$ and is contained neither in $LC^6$, nor $C^6R$ are contained inside $C^N$ and are periodic shifts of each other. Hence the winner choice for periodic shifts is the same and the result follows from Lemma~\ref{lem:turns in power words}.
\end{proof}

\begin{remark}\label{rem:certification in Can_r}
If $A=LC^{\tau}R\in\Can_r$, then by Corollaries~\ref{cor:powers in context} and~\ref{cor:periodic shift of certification} we see that the left and the right context together for any maximal occurrence $u$ in $A$ is contained in either $LC^{12}, C^{13}$ or $C^{12}R$. In particular, there is no maximal occurrence in $A$ whose left and right contexts have nontrivial overlap with both $L$ and~$R$.
\end{remark}

\begin{corollary}
\label{can_periodic_extension_constriction_cor}
IH~\ref{can_periodic_extension_constriction_hyp} holds for rank~$r$.
\end{corollary}

For the proof we first note the following:

\begin{lemma}[$\kappa = \mu - \varepsilon = n - 10\tau - 4$]
\label{scan_changing_power_subword}
Let $A = L_1C^{N_1}R_1, B = L_2C^{N_2}R_2 \in \Scan_{\kappa, r}$, where $C$ is primitive, $C^n \notin \Rel_1\cup\ldots \cup \Rel_r$, and $N_1, N_2 \geqslant \tau$. Then  $L_1C^SR_2\in \Scan_{\kappa, r}$ for any $S\geq\tau$.
\end{lemma}
\begin{proof}
Since $A, B \in \Scan_{\kappa, r} \subseteq \Can_{r - 1}$ and $C^n \notin \Rel_1\cup\ldots \cup \Rel_r$,  we have $L_1C^SR_2\in\Can_{r - 1}$ by IH~\ref{can_periodic_extension_constriction_hyp} in rank $r-1$.
If  $L_1C^SR_2$ contains a maximal occurrence $u$ of rank~$r$ of $\Lambda_r$-measure $>\kappa$, then $u$ is contained in $L_1C^2$, in $C^2R_2$, or in $C^3$ by Remark~\ref{occurrences_in_periodic_word} and Corollary~\ref{common_part_in_different_ranks}. This is impossible since $A, B \in \Scan_{\kappa, r}$.
\end{proof}

\begin{proof}[Proof of Corollary~\ref{can_periodic_extension_constriction_cor}.]
For $r=1$ this follows directly from the definition of $\can_1$ in Section~\ref{can_of_semican_section}. So assume $r>1$.
Let $X_1= L_1C^{\tau}R_1, X_2=L_2C^{\tau}R_2\in\Can_r$. Then $L_1C^SR_2\in\Scan_{\kappa, r}$ for any $S\geq\tau$ by Lemma~\ref{scan_changing_power_subword}. If all maximal occurrences in $L_1C^SR_2$ have $\Lambda_r$-measure $\leqslant \fracn - 5\tau - 2$, the claim follows directly from Lemma~\ref{too_small_and_too_big_side}. So let $u$ be a maximal occurrence in $L_1C^SR_2$ with $\Lambda_r(u)> \fracn - 5\tau - 2$.  By Remark~\ref{rem:certification in Can_r} we see that the left and right contexts of $u$ are contained in $LC^{12}, C^{13}$ or $C^{12}R$. So they coincide with the corresponding ones in $X_1$ or in $X_2$. Since any occurrence $u$ in $X_1, X_2$ is the winner side (because $X_1, X_2\in\Can_r$), the same is true for $u$ in $L_1C^SR_2$, proving the claim.
\end{proof}

\subsection{Multiplication and canonical triangles}
\label{can_triangle_section}

We now prove that the multiplication of canonical words of rank~$r$ can be described in terms of canonical triangles of rank~$r$. 

We say that a word $W$ contains \emph{a gap} if it contains a subword of the form $a^\tau M b^\tau M' c^\tau$ where $a^n, b^n, c^n \in \Rel_r$.

\begin{lemma}\label{lem:multiplication preparation small occurrence}
Let $A=LuW_0FW_1wR\in\Can_{r-1}$ where $u, w$ are maximal occurrences of $\Lambda_r$-measure $\geq \tau+1$, and $W_0, W_1$ do not contain strong separation words (from the right and left, respectively). Assume that at least one of the following conditions holds:
\begin{enumerate}
\item
\label{gap1}
$F$ is $\tau$-free of rank~$r$;
\item
\label{gap2}
$W_0$, $W_1$ do not contain gaps and $F = DzE$, where $D, E$ are $\tau$-free of rank~$r$, and $z$ is an occurrence of rank~$r$;
\item
\label{gap3}
at least one of $W_0$, $W_1$ is $\tau$-free of rank~$r$ and $F = DzE$ as above.
\item
\label{gap4}
$W_0$ does not contain a gap, $F = DE$, where $D, E$ are $\tau$-free of rank~$r$;
\item
\label{gap5}
$W_0$ contains a subword $b^{2\tau + 1}$, $b^n \in \Rel_r$, and $F = DzE$ as above.
\end{enumerate}
Let $C^N$ be a subword of $uW_0FW_1w$  where $C\in\Can_0$ is primitive and $C^n\notin\Rel_1 \cup\ldots \cup\Rel_r$. Then $N<\tau=15$. 
\end{lemma}
\begin{proof}
If $C\notin \Canc_{r - 1}$, then by Definition~\ref{canc_def} $A$ does not contain $C^{\tau}$. So we suppose that $C\in \Canc_{r - 1}$. Hence Lemma~\ref{lem:powers in between} implies that $W_0FW_1$ contains $\geqslant N - 4$ occurrences of $a^{\tau}$, $a^n \in \Rel_r$ not overlapping with each other. If one of Conditions~\ref{gap1}---\ref{gap4} holds, then the result follows from Example~\ref{isolation_words_examples} and Corollary~\ref{common_part_in_different_ranks} for a common part of $C^N$ and $z$.

Under Condition~\ref{gap4} either $C^N$ is contained in $uW_0$, or $C$ cyclically contains $b^{2\tau}$. In the first case the result follows from the above argument. In the second case we count directly periodic shifts of $b^{2\tau}$ and obtain $N < \tau$.
\end{proof}

We now show a preliminary version of IH~\ref{can_triangle_hyp1}:
 
\begin{proposition}[$\lambda=\fracn +3\tau+1$]
\label{can_triangle_multiplication}
Let $A, B \in \Can_r$. Then $\can_r(A\cdot B)=A_1'M_3 B_1'$, $A = A_1'M_1X$, $B = X\inv M_2B_1'$, where $X\cdot X\inv$ is the maximal cancellation in $A\cdot B$, and $M_3$ is $\tau$-free modulo $r$.
\end{proposition}
\begin{proof} 
By  IH~\ref{can_triangle_hyp1} for rank~$r - 1$, we know that $\can_{r - 1}(A\cdot B) = A''EB''$, where $E$ is $\tau$-free of rank~$r$. Since by Lemma~\ref{can subset scan} we have that $A, B\in\Scan_{\lambda, r}$, we  can apply Lemma~\ref{lem:mu' semicanonical} to $A''EB''$  and obtain  $C_1=A_1QB_1\in\Scan_{\lambda+3\tau+1,r}$ for a prefix $A_1$ of $A''$ and suffix $B_1$ of $B''$ by iterated turns of $\Lambda_r$-measure $> \lambda + (3\tau+1)$. By Lemma~\ref{turn_res_complement} $Q = DzE$, where $D, E$ are $\tau$-free of rank~$r$ and $z$ is an occurrence of rank~$r$ (any part can be empty).

Since $\lambda+3\tau+1 = \fracn + 6\tau + 2 < n - (8\tau + 3) - 2\tau - 1$, by Lemma~\ref{can_r_def_agree} $C=\can_r(A\cdot B)$ is obtained from $C_1=A_1QB_1$ by choosing the winner sides in the maximal occurrences of rank~$r$ in $C_1$. Clearly we can write $C = A_1'M_3B_1'$ for some prefix $A_1'$ of $A_1$, suffix $B_1'$ of $B_1$ and some word $M_3$. We need to show that it is possible to take $M_3$ $\tau$-free modulo $r$.

To determine the prefix $A_1'$ and the suffix $B_1'$, let $u$ and $w$ be the left- and right-most occurrences, respectively, which are turned in $C_1$. Let $u=u_0,\ldots, u_s=w$ be an enumeration from left to right of all maximal occurrences in $C_1 = A_1QB_1$ of $\Lambda_r$-measure $\geq 5\tau+3$ between $u$ and $w$. By Remark~\ref{rem:certification non-isolated} an initial segment of $u_0,\ldots, u_s$ is an initial segment of a (un-)certification sequence to the right of $u$ in $C$ (and in $A$ if $u$ is a maximal occurrence in $A$), a final segment of $u_0,\ldots, u_s$ is a final segment of a (un-)certification sequence to the left of $w$ in $C$ (and in $B$ if $w$ is a maximal occurrence in $B$).

First suppose that $u$ is contained in $A_1$ and $w$ is contained in $B_1$. Since the winner side for $u = u_0$ is different in $A$ and in $C$, $A_1$ cannot contain a right context for $u_0$. Hence if a common part of $u_k$ and $A_1$ has $\Lambda_r$-measure $\geqslant \tau + 1$, then  consecutive occurrences in $(u_0, \ldots, u_k)$ are essentially not isolated. If a common part of $u_k$ and $A_1$ has $\Lambda_r$-measure $\geqslant 5\tau + 3$, then $(u_0, \ldots, u_k)$ is an initial segment of an (un-)certification sequence in $A$ and in $C_1$ by Condition~5 or 5' of Definition~\ref{certification_def}. The corresponding properties hold for $(u_k, \ldots, u_s)$.

Let $i$ be the maximal index such that $A_1$ does not contain $u_i$  as a suffix, and $j$ be the minimal index such that $B_1$ does not contain $u_j$ as a prefix. Then $0 < j - i \leqslant 4$, since $Q = DzE$.

We now turn all necessary occurrences in $C_1$ according to the choices of the winner sides. Denote the occurrences corresponding to $u_k$ in $C$ by $f_k$. Since $A_1$, $B_1$ do not contain contexts for $u_0$ and $u_s$, respectively, there exist at most $3$ maximal occurrences of $\Lambda_r$-measure $\geqslant 5\tau + 3$ in $C$ between $f_i$ and $f_j$.

Corollaries~\ref{cor:turn outside} and~\ref{cor:turns inside} imply that either $(f_0, \ldots, f_i)$ is an initial segment of an (un-)certification sequence in $C$ from the right of $f_0$, or $\Lambda_r(f_i) < 5\tau + 3$. Hence in the second case there exist at most $3$ maximal occurrences of $\Lambda_r$-measure $\geqslant 5\tau + 3$ in $C$ between $f_{i - 1}$ and $f_j$. The symmetric property holds from the other side. So we obtain the following sequence of maximal occurrences of $\Lambda_r$-measure $\geqslant 5\tau + 3$ in $C$ between $f_0$ and $f_s$: $(f_0, \ldots, f_{i_0})$ is an initial segment of an (un-)certification sequence, where either $i_0 = i - 1$, or $i_0 = i$, $(f_{j_0}, \ldots, f_s)$ is a final segment of an (un-)certification sequence, where either $j_0 = j$, or $i_0 = j + 1$, and there exist at most $3$ maximal occurrences of $\Lambda_r$-measure $\geqslant 5\tau + 3$ in $C$ between $f_{i_0}$ and $f_{j_0}$.

Let $V$ be primitive such that $V^n \notin \Rel_1\cup \ldots \cup \Rel_r$. If $V\notin \Canc_{r - 1}$, we are done by the definition of $\Canc_{r - 1}$. So assume $V \in \Canc_r$. Assume that $V^N$ contains some maximal occurrence $x$  with $\Lambda_r(x) \geqslant 5\tau + 3$, and not as a prefix or suffix. Then by Corollary~\ref{cor:periodic subwords in certification} and by the previous considerations $V^N$ contains $\leqslant 11$ periodic shifts of $x$, so $N \leqslant 13$.

Now assume that $V^N$ does not contain any occurrence with $\Lambda_r$-measure $\geqslant 5\tau + 3$. Clearly it is sufficient to consider a space between $f_{i_0}$ and $f_{j_0}$, which is of the form $W_1QW_2$, where $W_1, W_2$ do not contain strong separation words. If $W_1$ does not have a common suffix with $A_1$ or $W_2$ does not have a common prefix with $B$, then the result follows from Lemma~\ref{lem:multiplication preparation small occurrence}~(\ref{gap3}). If $Q = E$, the result follows from Lemma~\ref{lem:multiplication preparation small occurrence}~(\ref{gap1}). Assume that $Q = DzE$ with non-empty~$z$ with $\Lambda_r(z) < 5\tau + 3$. Then the last occurrence that is turned in order to obtain $C_1$ is of $\Lambda_r$-measure $> n - (7\tau + 3)$. Hence it has common parts both with $A$ and $B$ of $\Lambda_r$-measure $> \fracn - 5\tau - 2$. Denote their maximal prolongations by $w_1$ and $w_2$, respectively. If $W_1$ contains a gap, then $w_1$ is strongly isolated from $f_{i_0}$. Hence there must exist some $x$ in $A$ with $\Lambda_r(x) \geqslant 5\tau + 3$ between $f_{i_0}$ and $w_1$ (otherwise the winner side for $u_0$ is the same in $A$ and $C_1$). Since the space between $f_{i_0}$ and $w_1$ is of the form $W_1D_1$ with $D_1$ $\tau$-free of rank~$r$, $W_1$ contains $b^{2\tau + 1}$, $b^n \in \Rel_r$. The symmetric property holds for $w_2$ and $W_2$. So the result follows from Lemma~\ref{lem:multiplication preparation small occurrence} (\ref{gap2}) and~(\ref{gap5}). If $Q = DE$, then we argue in the same way but only from the left side. Then the result follows from Lemma~\ref{lem:multiplication preparation small occurrence} (\ref{gap4})

Finally we need to consider the case that $u$ or $w$ are not contained in $A_1, B_1$, respectively. Suppose that $u$ is not contained in $A_1$. Then the sequence $(u=u_0,\ldots, u_m=w)$ contains at most three maximal occurrences not contained in $B_1$ and similarly for the other side. Thus we see from the previous arguments that $M_3$ is $\tau$-free modulo rank~$r$.
\end{proof}

Now we can finish the proof of IH~\ref{can_triangle_hyp1} for rank~$r$.
\begin{corollary}
\label{can_triangle_multiplication_cor}
Let  $A, B \in \Can_r$ and $\can_{r-1}(A\cdot B) = A''E_3B''$ by IH~\ref{can_triangle_hyp1} for rank~$r - 1$. There exists a canonical triangle $(D_1, D_2, D_3)$ of rank~$r$ such that $\can_r(A\cdot B) = A_1D_3B_1$, $A = A_1D_1X$, $B = X\inv D_2B_1$, where $X\cdot X\inv $ is the maximal cancellation in $A\cdot B$ and $A_1, B_1$ are prefix and suffix of $A'', B''$, respectively. Furthermore if $A_1 = A''$ and $B_1= B''$, then $\can_r(A\cdot B) = \can_{r-1}(A\cdot B)$.
\end{corollary}
\begin{proof}
By Proposition~\ref{can_triangle_multiplication} we have $\can_r(A\cdot B) = A'M_3B'$, $A = A'M_1X$, $B = X\inv M_2B'$, where $X\cdot X\inv $ is the maximal cancellation in $A\cdot B$. By IH~\ref{can_triangle_hyp1} for rank~$r - 1$ we have that $A = A''E_1X$ and $B = X\inv E_2B''$, where $(E_1, E_2, E_3)$ is a canonical triangle of rank~$r - 1$. By construction we have that $A'$ is a prefix of $A''$ and $B'$ is a suffix of $B''$.

First suppose that $A' = A''$ and $B' = B''$. Then by construction all turns are done in $E_3$. However, since $E_3$ is $\tau$-free of rank~$r$, it does not contain any occurrences of rank~$r$ to turn. So there are no turns in $A''E_3B''$ in order to obtain $\can_r(A\cdot B)$, hence $\can_r(A\cdot B) = \can_{r-1}(A\cdot B)$ and we can put $D_i = E_i$, $i = 1, 2, 3$.

By definition of $\cdot_r$, we can write $A\cdot_r B = \can_r(A\cdot B) = C$. Therefore $A = C \cdot_r B\inv $. Now we apply Proposition~\ref{can_triangle_multiplication} to $A = C \cdot_r B\inv $ and obtain $A = C'F_3{B'}\inv$, where $F_3$ is $\tau$-free modulo rank~$r$, $C'$ is a prefix of $C$, ${B'}\inv$ is a suffix of $B\inv$. Since $X\cdot X\inv$ is the maximal cancellation in $A\cdot B$, we have that ${B'}\inv$ is a suffix of $X\inv$.

If $C'$ is a prefix of $A'$, then $M_1$ is a subword of $F_3$, since ${B'}\inv$ is a suffix of $X\inv$. So, $M_1$ is $\tau$-free modulo rank~$r$.

Otherwise $A'$ is a proper prefix of $C'$, so $C' = A'W$. Since $C'$ is left after the maximal cancellations in $C\cdot B\inv$, we obtain that $W$ is a prefix of $M_3$. Then $A'W$ is a prefix of $A$, because $C'$ is also a prefix of $A$. If $M_1$ is contained in $W$, then $M_1$ is a subword of $M_3$, so it is $\tau$-free modulo~$r$. If $W$ is a proper prefix of $M_1$, then $W$ is a common prefix of $M_1$ and $M_3$. In this case we fold $W$ in the sides $M_1$ and $M_3$ and obtain a new triangle $\tilde{M}_1$, $\tilde{M}_2 = M_2$, $\tilde{M}_3$. Then $\tilde{M}_1$ is a prefix of $F_3$ and $\tilde{M}_3$ is a suffix of $M_3$. So, $\tilde{M}_1$ and $\tilde{M}_3$ are $\tau$-free modulo rank~$r$.

Assume that after the above procedure side $E_1$ is not a suffix of $\tilde{M}_1$ anymore. Then instead of complete folding $W$ we fold it until $E_1$ plus one extra letter. Then $\tilde{M}_1 = xE_1$ for some single letter $x$. Since $E_1$ is $\tau$-free modulo rank~$r - 1$, $\tilde{M}_1$ is $\tau$-free modulo rank~$r$ by Lemma~\ref{lem:nesting_occurrences_rel_r}.

After that we deal similarly with $\tilde{M}_2 = M_2$ in the new triangle and as a result obtain the required canonical triangle $(D_1, D_2, D_3)$.

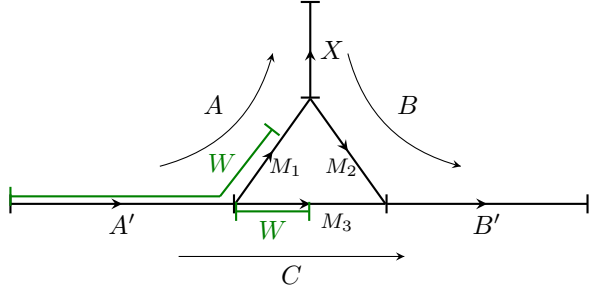
\begin{center}

\begin{tikzpicture}

\draw[|-|, thick, arrow=0.5] (0, 0) to node[midway, below] {$A'$} (3, 0);

\draw[thick, arrow=0.5] (3, 0) to node[midway, below, xshift=10] {\footnotesize $M_3$} (5, 0);

\draw[|-|, thick, arrow=0.5] (5, 0) to node[midway, below] {$B'$} (7.7, 0);

\draw[thick, arrow=0.5] (3, 0) to node[midway, right, xshift=-5, yshift=-6] {\footnotesize $M_1$} (4, 1.4);

\draw[thick, arrow=0.5] (4, 1.4) to  node[midway, left, yshift=-6, xshift=7.5] {\footnotesize $M_2$}(5, 0);

\draw[|-|, thick, arrow=0.5] (4, 1.4) to node[midway, right] {$X$} (4, 2.7);

\draw[black, arrow = 1] (-1.5 + 3.75, -0.7)--(1.5 + 3.75, -0.7) node[midway, below] {$C$};

\draw[black, arrow = 1] (-2 + 4, 0.5) to [bend right] node[midway, left, yshift=8] {$A$} (-0.5 + 4, 2);

\draw[black, arrow = 1] (0.5 + 4, 2) to [bend right] node[midway, right, yshift=8] {$B$} (2 + 4, 0.5);
\draw[thick, black!50!green, |-] (0, 0.1) to (3 - 0.2, 0.1);

\draw[thick, black!50!green, -|] (3 - 0.2, 0.1) to node[midway, left] {$W$} (3.5, 1);

\draw[thick, black!50!green, |-|] (3, -0.1) to node[midway, below] {$W$} (4, -0.1);

\end{tikzpicture}

\end{center}

In this case we can fold $W$ in the sides $D_1$ and $D_3$ and obtain a new triangle $\tilde{D}_1$, $\tilde{D}_2 = D_2$, $\tilde{D}_3$, where $\tilde{D}_1$ is a prefix of $F_3$ and $\tilde{D}_3$ is a suffix of $D_3$. So, $\tilde{D}_1$ and $\tilde{D}_3$ are $\tau$-free modulo rank~$r$. 

After that we deal similarly with $\tilde{D}_2 = D_2$ in the new triangle and as a result obtain the required canonical triangle $(D_1, D_2, D_3)$. 
\end{proof}

\subsection{Canonical form of power words}
\label{can_r_of_powers_section}

\medskip\noindent

We start with some preliminary lemmas:

\begin{lemma}
\label{accumulation_of_powers_cases}
Let $A= XWX\inv\in \Can_{r - 1}$ and  $\can_{r - 1}(A\cdot A)=X_1W_1(X_1)\inv$
where $W$ and $W_1$ are cyclically reduced. If $W$ is $\tau$-free of rank~$r$, then $W_1$ is $3\tau$-free of rank~$r$.
\end{lemma}
\begin{proof}
By IH~\ref{can_triangle_hyp1} there exists a canonical triangle $(D_1, D_2, D_3)$ of rank $r-1$ such that $XW=A'D_1$, $WX\inv=D_2A''$ and  $\can_{r-1}(A\cdot A)= A'D_3A''$. If $A'=XW'$ and $A''=W''X\inv$ (where $W'$ or $W''$ may be empty), then $X=X'$ and the claim is immediate.
Otherwise, either $A' = X'$, $A'' = W_0X_0\inv (X')\inv$, or symmetrically $A' = X'X_0W_0$, $A'' = (X')\inv$, where $X'$ is a prefix of $X$, $X_0$ is $\tau$-free subword of $X$ and $W_0$ is a subword of $W$. Then  again the claim is immediate.
\end{proof}

\begin{lemma}\label{lem:merging of powers}
Let $W', W''$ be $3\tau$-free of rank~$r$ and let $D, E$ be $\tau$-free  of rank~$r$.
If $(EW'DW'')^N$ contains an occurrence $u$ of $\Lambda_r$-measure $\geq 11\tau + 1$, then $EW'DW''=a^s$ for some $s\geq 0$ and $a^n\in\Rel_r$.
\end{lemma}
\begin{proof}
Let $u = \hat{a}^k\hat{a}_1$ for $\hat{a}^n \in \Rel_r$ and $k\geq 11\tau+1$. If $|\hat{a}|\geq |EW'DW''|$, then by Lemma~\ref{common_part_of_powers} $EW'DW''$ is a cyclic shift of $\hat{a}$.

If $|\hat{a}| < |EW'DW''|$, the assumptions on  $W', W'', D, E$ imply that $u$ contains a cyclic shift $Y$ of $EW'DW''$, and $Y$ is $11\tau$-free. So since $\Lambda_r(u) \geq 11\tau + 1$, the common part of $u$ and $(EW'DW'')^N$ has length $\geq |Y|+|\hat{a}|$. Hence by Lemma~\ref{common_part_of_powers} $EW'DW''=a^s$ for $a$ a cyclic shift of $\hat{a}$.
\end{proof}

\begin{lemma}\label{lem:power with 3tau+2}
Let $A= XWX\inv \in \Can_{r - 1}$ where $W$ is cyclically reduced and contains an occurrence $u$ of $\Lambda_r$-measure $\geq 3\tau$.
Then there is a canonical triangle $(D_1, D_2, D_3)$ of rank $r - 1$ such that \[Q=\can_{r-1}(\underbrace{A\cdot\ldots\cdot A}_{N \textit{times}})=XD_2(MD_3)^{N-1}  MD_1X\inv\]
where $W=D_2MD_1$. In particular, $MD_3$ is conjugate to $W$ in $\Fr/\llangle\Rel_0, \ldots, \Rel_{r - 1}\rrangle$.
Furthermore, if $W$ is $\kappa$-bounded of rank~$r$ for some $\kappa\geq 3\tau$, then either $(MD_3)^N$ is $2\kappa+\tau+1$-free of rank~$r$, or $MD_3 = a^s$ for $a^n \in \Rel_r$.
\end{lemma}
\begin{proof}
By IH~\ref{can_triangle_hyp1} there is a canonical triangle $(D_1, D_2, D_3)$ of rank $r - 1$ such that $\can_{r-1}(A\cdot A)=XW'D_3W''X\inv$ for some non-empty prefix $W'$ and suffix $W''$ of $W$. Since $W$ contains $u$, we can write $W=D_2M D_1$ with non-empty~$M$. Then $\tilde{W}=D_1D_2M$ is a cyclic conjugate of $W$ and we have
\[\tilde{W}=D_1D_2M\equiv D_3M\mod \llangle \Rel_0\cup\ldots\cup\Rel_{r - 1}\rrangle.\]
Since $D_1$ and $D_2$ is $\tau$-free of rank~$r$, $W'$ and $W''$ contain occurrences of $\Lambda_r$-measure $\geq 2\tau$. Hence by Corollary~\ref{can_power_context_hyp_cor1} we obtain 
\[\can_{r-1}(A\cdot A\cdot A)=XW'D_3MD_3 W''X\inv \textit{ and } W' = D_2M,\ W'' = MD_1.\] 
Inductively 
Corollary~\ref{can_power_context_hyp_cor1} yields
\[
Q = XW'\underbrace{D_3M\ldots D_3M}_{N - 2 \textit{ times}}D_3W''X\inv=XD_2(MD_3)^{N-1}MD_1X\inv.\]
The last sentence follows from Lemma~\ref{common_part_of_powers}.
\end{proof}

\begin{lemma}
\label{no_accumulation_of_powers_cor}\label{cor: no_accumulation_of_powers2}
Let $A=XWX\inv \in \Can_{r - 1}$, where $W$ is a cyclically reduced and $3\tau$-free of rank~$r$ and $W^n \notin \llangle  \Rel_1, \ldots, \Rel_r\rrangle$. By IH~\ref{can_power_hyp} for rank~$r - 1$ write
\[
\can_{r - 1}(\underbrace{A\cdot \ldots \cdot A}_{K \textit{ times}}) = T\tilde{A}^{K - \gamma}S \ \textit{ for all } K \geq \gamma,
\]
where $\tilde{A}, T, S, \gamma$ depend only on $A$ and $r$, and $A, \widetilde{A}$ are conjugate in the group $\Fr / \llangle  \Rel_1, \ldots, \Rel_{r - 1}\rrangle$. 
Then $\tilde{A}^N$ is $11\tau+1$-free of rank~$r$ for all $N \geq 1$  and hence
\[
\can_r(\underbrace{A\cdot \ldots \cdot A}_{K \textit{ times}}) = T'\tilde{A}^{K - \gamma'}S' \ \textit{ for all } K \geqslant \gamma'
\]
where $T', S'$ and $\gamma'$ only depend on $A$ and $r$. (Note that $\tilde{A}$ does not change.)
\end{lemma}
\begin{proof}
Let $A_s =  \can_{r - 1}(\underbrace{A\cdot \ldots \cdot A}_{2^s \textit{ times}}) = X_sW_sX_s^{-1}$, where $W_s$ is cyclically reduced.

First assume that $W_s$ is $3\tau$-free for all $s \geqslant 0$. For all $K = 2^s$ we have  $T\tilde{A}^{K - \gamma}S = X_sW_sX_s\inv$. Notice that if an overlap of $X_s$ and $\tilde{A}^{K - \gamma}$ contains a whole period $\tilde{A}$, then an overlap of $X_s^{-1}$ and $\tilde{A}^{K - \gamma}$ cannot contain $\tilde{A}$. So if there exists $s \geq \gamma + 3$ such that $T, S$ are contained in $X_s, X_s^{-1}$, respectively, then $W_s$ contains $\tilde{A}^3$. Otherwise $|X_s|\leq \max\{ |S|, |T|\}$ for $s \geq \gamma + 3$  because $|S|, |T|$ do not depend on $s$. In this case $W_s$ contains $\tilde{A}^3$ for all sufficiently large~$s$. Thus $\tilde{A}^3$ is $3\tau$-free of rank~$r$, and so is $\tilde{A}^{N_1}$ for all $N_1 \geqslant 1$.

Now let $t\geq 1$ be minimal such that $W_t$ is not $3\tau$-free of rank~$r$ so   $W_{t - 1}$ is not $\tau$-free of rank~$r$ by Lemma~\ref{accumulation_of_powers_cases}. Therefore by IH~\ref{can_triangle_hyp1}, $W_t = W_{t - 1}'EW_{t - 1}''$, where $E$ is $\tau$-free of rank~$r$, $W_{t - 1}', W_{t - 1}''$ are a non-empty suffix and prefix of $W_{t - 1}$, respectively. Then it follows from Lemma~\ref{lem:power with 3tau+2} that
\begin{equation*}
\can_{r - 1}(\underbrace{A_t\cdot \ldots \cdot A_t}_{N \textit{ times}}) = X_tD_2(MD_3)^{N-1}MD_1X_t,
\end{equation*}
where $W_t = D_2MD_1$, and $(D_1, D_2, D_3)$ is a canonical triangle of rank~$r - 1$. Since $W_{t - 1}$ is $3\tau$-free of rank~$r$, by construction, $MD_3$ is of the form $W'E'W''D_3$, where $E'$ is $\tau$-free and $W', W''$ are $3\tau$-free of rank~$r$ (some parts can be empty). Hence Lemma~\ref{lem:merging of powers} implies that either $(MD_3)^{N-1}$ is $11\tau + 1$-free of rank~$r$, or $MD_3 = a^k$ for some $a^n \in \Rel_r$.

For sufficiently large $N$ and $K = N\cdot 2^t$,  the common part of $(MD_3)^{N-1}$ and $\tilde{A}^{K - \gamma}$ has length $> |MD_3| +|\tilde{A}|$. Hence by Lemma~\ref{common_part_of_powers} $MD_3 = Z^{k_1}$, $\widetilde{A} = Z^{k_2}$ for some word $Z$. So if $(MD_3)^{N-1}$ is $11\tau + 1$-free of rank~$r$, $\tilde{A}^{N_1}$ is also $11\tau + 1$-free of rank~$r$ for all $N_1 \geqslant 1$. If $MD_3 = a^k$, then $\widetilde{A} = a^{k_2}$, since $a$ is primitive. Hence $\widetilde{A}^n \in \Rel_r$, a contradiction.

For sufficiently large $N$ by IH~\ref{can_of_prefix_and_suffix_hyp} we have $\can_{r - 1}(\tilde{A}^N) = D\tilde{A}_1\tilde{A}^{N - \delta}\tilde{A}_2E$, where $D, E$ are $\tau$-free of rank~$r$, $\tilde{A}_1,  \tilde{A}_2$ are a suffix and prefix of $\tilde{A}$, respectively. By IH~\ref{can_periodic_extension_constriction_hyp} $D, E, \tilde{A}_1, \tilde{A}_2, \delta$ do not depend on $N$. By Lemma~\ref{no_accumulation_of_powers_cor} $D\widetilde{A}_1\widetilde{A}^{N - \delta}\widetilde{A}_2E$ is $12\tau + 1$-free of rank~$r$ for big enough $N$. Since $12\tau + 1 \leqslant \fracn - 5\tau - 2$, it follows from Lemma~\ref{too_small_and_too_big_side} that $\can_{r - 1}(\tilde{A}^N) \in \Can_r$, so $\can_r(\tilde{A}^N) = \can_{r -1}(\tilde{A}^N)$.
For sufficiently large $K$ we obtain
\begin{align*}
\can_r(\underbrace{A\cdot \ldots \cdot A}_{K \textit{ times}})\  &=\  \can_r(T)\cdot_r \can_r(\tilde{A}^{K - \gamma})\cdot_r\can_r(S) \\
&=\  \can_r(T)\cdot_r  \left(D\widetilde{A}_1\widetilde{A}^{K - \gamma - \delta}\tilde{A}_2E\right)\cdot_r \can_r(S)\ = \ T'\tilde{A}^{K - \gamma'}S'.
\end{align*}
\end{proof}

\begin{lemma}[$\mu = n - 8\tau_1 - 3$]
\label{lem:greedy powers}
Let $A \in \Canc_{r - 1}$ such that $A^n \notin \Relr$. Then there exists $\widetilde{B}$ conjugate to $A$ in the group $\Fr / \llangle \Rel_0, \ldots, \Rel_r\rrangle$ such that one of the following holds:
\begin{itemize}
\item
$\widetilde{B} \in \Can_{r - 1}$ and its cyclically reduced part is $3\tau$-free of rank~$r$.
\item
$\widetilde{B} \in \Canc_{r - 1}$  and $\widetilde{B}^K$ is $\mu - \tau$-bounded for all $K \geqslant 1$.
\end{itemize}
\end{lemma}
\begin{proof}
Similarly to the proof of Lemma~\ref{lem:scan_greedy_alg}, we do induction on $d'(A) = \max\{d(X)\}$, where $X$ runs through all cyclic shifts of $A$ and $d(X)$ is the sum of all $\Lambda_r$-measures of all maximal occurrences of rank~$r$ in $X$ with $\Lambda_r$-measure $\geq \beta = 3\tau+2$.

Let $Y = \can_{r - 1}(A^{N_1}) = LA^NR$, $N = N_1 - \delta$, and let $u$ be a maximal occurrence of rank~$r$ properly contained in $A^N$ with $\Lambda_r(u) > \mu - \tau$. Since $A^n \notin\Relr$, by Remark~\ref{occurrences_in_periodic_word} $u$ is contained in $A^3$ and there are periodic shifts of $u$ starting in every period of $X$ that are maximal occurrences in $Y$.

We now turn the occurrence of $u$ in $Y$ starting in the second period of $A^N$. Then there are the following configurations in the resulting word.

\medskip\noindent
1. There exists a cyclic shift  $V_2V_1$ of $A = V_1V_2$ of the form $LuR$ such that the result of the turn of $u$ in $Y$ is of the form $(LV_1)(L'Q\tilde{R})V_2A\cdots AR$, where $L'$ is a prefix of $V_2$ that contains $b_1^{\tau}$ and $\tilde{R}$ is a suffix of $uR$ that contains $b_2^{\tau}$ for some $b_1^n, b_2^n \in \Rel_r$. Then Corollary~\ref{can_power_context_hyp_cor1} implies that the result of the turns of all periodic shifts of $u$ in $A$ is equal to $(LV_1)(L'Q\tilde{R})\cdots (L'Q\tilde{R})V_2R$, so $L'Q\tilde{R} \in \Canc_{r - 1}$. Clearly $L'Q\widetilde{R}$ and $A$ are conjugate in the group $\Fr / \llangle \Rel_0, \ldots, \Rel_r\rrangle$.

We have 
\begin{align*}
d(L^{\prime}Q\widetilde{R}) &\leqslant d(V_2V_1) - (\mu - \tau) + (n - \mu + 3\tau) + 2\beta + 2\varepsilon =\\
&= d(V_2V_1) - (\mu - \tau) + (n - \mu + 3\tau) + 2(3\tau + 2) + 2(2\tau + 1) =\\
&= d(V_2V_1) - (n - 30\tau - 12).
\end{align*}
Hence $d'(L^{\prime}Q\widetilde{R}) < d(L^{\prime}Q\widetilde{R}) + 2\beta \leqslant d(V_2V_1) - (n - 30\tau - 12) + 6\tau + 4 = d(V_2V_1) - (n - 36\tau - 16) < d'(A)$. So the claim holds by the induction hypothesis.

\medskip\noindent
2. Assume that $\abs{u} \leqslant \abs{A}$ and we are not in Case~1. Then for every cyclic shift of $A$ of the from $LuR$ we see that $L$, $R$ do not contain words of the form $a^{\tau}Mb^{\tau}$, $a^n, b^n \in \Rel_r$. Consider $\can_{r - 1}(LuR) = L_1\widetilde{u}R_1$, where $\Lambda_r(u) - 2\tau < \Lambda_r(\widetilde{u}) < \Lambda_r(u) + 2\tau$. Let us turn $\widetilde{u}$ and let $B$ be the resulting word. If the cyclically reduced part of $B$ is $3\tau$-free of rank~$r$, then we take it as $\widetilde{B}$ and we are done.

Let $\widetilde{v}$ be the complement of $\widetilde{u}$. Then the turn is of Type~2 and $\Lambda_r(\widehat{\widetilde{v}}, B) < n - (\mu - 3\tau) + 2\tau = 13\tau + 3$. So $B$ is $13\tau + 3$-bounded. Then by Lemma~\ref{lem:power with 3tau+2} the periodic part of $\can_{r - 1}(B\cdot\ldots\cdot B)$ is $27\tau + 7$-bounded. Since $27\tau + 7 < \mu - \tau = n - 9\tau - 3$, we take a period of this periodic part as $\widetilde{B}$.

\medskip\noindent
3. The last case is $|u|>|A|$. Then there exists a cyclic shift of $A$ equal to $u_1$ a prefix of $u$ with $\FracM_r(u_1) > \FracM_r(u) - 1 > \mu - \tau - 1$. Then we take $\can_{r - 1}(u_1) = L_1\widetilde{u}_1R_1$ and argue as above. Using the same notations, we have that $\Lambda_r(\widehat{\widetilde{v}}, B) < n - (\mu - 3\tau - 1) + 2\tau = 13\tau + 4$, so $B$ is $13\tau + 4$-bounded. Hence the periodic part of $\can_{r - 1}(B\cdot\ldots\cdot B)$ is $27\tau + 9$-bounded. Since $27\tau + 9 < \mu - \tau = n - 9\tau - 3$, we take a period of this periodic part as $\widetilde{B}$.
\end{proof}

\begin{proposition}\label{prop:IH 18}
IH~\ref{can_power_hyp} holds for rank~$r$.
\end{proposition}
\begin{proof}
By IH~\ref{IH can of equivalent words} and IH~\ref{can_power_hyp} for rank~$r - 1$ we can assume that $A \in \Canc_{r - 1}$. Let $\tilde{B}$ be the conjugate of $A$ given  by Lemma~\ref{lem:greedy powers}. If $\widetilde{B} \in \Can_{r - 1}$ and its cyclically reduced part is $3\tau$-free of rank~$r$, the result follows from Corollary~\ref{cor: no_accumulation_of_powers2}.

If the second case of Lemma~\ref{lem:greedy powers} holds, then $\can_{r - 1}(\widetilde{B}^K) = Z_1\widetilde{B}_1\widetilde{B}^{K - \delta}\widetilde{B}_1Z_2 \in \Scan_{\mu, r}$, where $\widetilde{B}_1, \widetilde{B_2}$ are a prefix and suffix of $B$, respectively, $Z_1, Z_2$ are $\tau$-free of rank~$r$, and $\delta$ does not depend on $K$. Then Corollary~\ref{periodic_can_subword} implies that $\can_r(\widetilde{B}^K) = \widetilde{L}\widetilde{X}^{K - \gamma}\widetilde{R}$, where $\widetilde{X}$ and $\widetilde{B}$ are conjugate in the group $\Fr / \llangle \Rel_1, \ldots, \Rel_r\rrangle$, and $\widetilde{X}, \gamma$ do not depend on $K$. Since $A \equiv Y\cdot \widetilde{B}\cdot Y^{-1} \mod \llangle \Rel_1, \ldots, \Rel_r\rrangle$, we have 
\[\can_r(\underbrace{A\cdot \ldots \cdot A}_{K \textit{ times}}) = \can_r(Y)\cdot_r \can_r(\widetilde{B}^K)\cdot_r \can_r(Y^{-1}).\] So, $A$ satisfies IH~\ref{can_power_hyp} with $\widetilde{A} = \widetilde{X}$.
\end{proof}

\section{Completion of the proof of Theorem~\ref{main thm}}\label{proof_completion_section}
It is left to show that the canonical form stabilizes and that our relators $\bigcup_{i\in\mathbb{N}} \Rel_i$ yield a quotient group isomorphic to the free Burnside group $B(m,n)$.
We start with the first point:

\begin{lemma}
\label{can_stabilization_aux}
Assume that $A, B \in \cap_{i = 0}^{\infty}\Can_i$. Then there exists $r_0$ such that $\can_{r_0}(A\cdot B) \in \cap_{i = 0}^{\infty}\Can_i$.
\end{lemma}
\begin{proof}
Since $A, B \in \Can_i$ for all $i \geqslant 0$, by IH~\ref{can_triangle_hyp1}  we have
\[
\can_i(A\cdot B) =  A_iD_3^{(i)}B_i,
\ A = A_iD_1^{(i)}X,\ B = X\inv D_2^{(i)}B_i,
\]
where $X\cdot X\inv $ is the maximal cancellation in $A\cdot B$ and $(D_1^{(i)}, D_2^{(i)}, D_3^{(i)})$ is a canonical triangle of rank~$i$,  $A_{i+1}$ is a prefix of $A_i$ and $B_{i+1}$ is a suffix of $B_i$. Let $r_0\geq 0$ be such that for all $i\geq r_0$ we have  $A_{r_0}=A_i, B_{r_0}=B_i$.
 Since the maximal cancellation $X$ does not depend on $i$, this implies $D_1^{(i)} = D_1^{(r_0)}$ and $D_2^{(i)} = D_2^{(r_0)}$ and hence, by IH~\ref{can_triangle_hyp1}, also $D_3^{(i)} = D_3^{(r_0)}$ for all $i\geq r_0$. We obtain
$\can_i(A\cdot B) = 
A_{r_0}D_3^{(r_0)}B_{r_0}=\can_{r_0}(A\cdot B)$ for all $i\geq r_0$ and $r_0$ is as required.
\end{proof}

\begin{proposition}
\label{can_stabilization}
For every word $A\in \Can_{-1}$ there exists $r_0$ such that $\can_{r_0}(A) \in \cap_{i = 0}^{\infty}\Can_i$.
\end{proposition}
\begin{proof}
Clearly we may assume that $A$ is reduced, so $A \in \Can_0$, and do induction on $|A|$. If $|A| = 1$, then it follows from Remark~\ref{final_can_immediately_condition} that $A \in \cap_{i = 0}^{\infty}\Can_i$.

For the induction step assume that $A = A_1x$, where $x$ is a single letter. By our induction assumption there is some $s$ such that $\can_s(A_1), \can_s(x)=x\in  \cap_{i = 0}^{\infty}\Can_i$. By Lemma~\ref{can_stabilization_aux} and Corollary~\ref{can_of_can_product_hyp}, there exists some $r_0\geq s$ such that
\[
\can_{r_0}(A) = \can_{r_0}(A_1x) = \can_{r_0}(\can_{r_0}(A_1)\cdot x)\in \cap_{i = 0}^{\infty}\Can_i.
\]
\end{proof}

By Proposition~\ref{can_stabilization} the sequence $\can_i(A)$, $i \geqslant 0$, stabilizes after a finite number of steps (depending on $A$). Therefore we can now define: 
\begin{definition} 
\label{can_stabilized_def}
For $A\in \Can_{-1}$, \emph{the canonical form}  $\can(A)$ of $A$  is defined as $\can(A)=\can_i(A)$ where $i$ is such that $\can_i(A)\in \bigcap_{i = 0}^{\infty}\Can_i$, and $\Can=\{\ \can(A)\mid \ A\in  \Can_{-1}\} = \bigcap_{i = 0}^{\infty}\Can_i$.
\end{definition}

It follows directly from the definition, IH~\ref{IH idempotent} and Remark~\ref{can_preserves_coset_hyp} that we have
\begin{corollary}
\label{final_can_idempotent}
$\can(\can(A) )= \can(A)\equiv A \mod\llangle\ \Rel_i\mid\ i\geq 0\ \rrangle$ for $A \in \Can_{-1}$.
\end{corollary}

\begin{lemma}
\label{power_subwords_of_can_all_ranks}
Let $a$ be a primitive word and $a^{\tau}$ be an occurrence in $A_i \in \Can_i$ for every $i \geqslant 0$. Then $a^n \in \Rel_r$ for some $r\geq 0$.
\end{lemma}
\begin{proof}
By the definition of $\Canc_i$, $a\in\Canc_i$ for all $i\geq 0$. The proof is by induction on $\abs{a}$. If $\abs{a} = 1$, then by definition $a^n \in \Rel_1$.

If $|a| > 1$ and $a$ cyclically contains $b^{\tau}$, then $|b|<|a|$. By the induction hypothesis $b^n \in \Rel_j$ for some~$j$. Let $r$ be minimal such that $a$ does not cyclically contain any occurrence of the form $b^{\tau}$ with $b^n\in\Rel_r$. Then we have that $a^n\in\Rel_r$ by the definition of $\Rel_r$.
\end{proof}

\begin{proposition}
\label{can_of_nth_powers}
If $A=XWX\inv \in \Can_0$ with $W$ cyclically reduced, then $W^n\in \llangle \Rel_i\mid\ i\in\mathbb{N}\rrangle$.
\end{proposition}
\begin{proof}
Write $\mathcal{H}=\llangle \Rel_i\mid\ i\in\mathbb{N} \rrangle$ and suppose $W^n\notin\mathcal{H}$.  By Corollary~\ref{final_can_idempotent} we may assume $A=\can(A)$. Let $r$ be minimal such that $A$ does not contain  any maximal occurrence of rank~$r$ of $\FracM_r$-measure $\geq 3\tau$. By IH~\ref{can_power_hyp} for all $j \geq 0$ we have
\begin{equation*}
\can_j(\underbrace{A\cdot \ldots \cdot A}_{K \textit{ times}}) = T_j\tilde{A}_j^{K - \gamma_j}S_j \textit{ for } K \geqslant \gamma_j,
\end{equation*}
where $\tilde{A}_j, T_j, S_j, \gamma_j$ depend only on $A$ and $j$, and $A$ and $\tilde{A}_j$ are conjugate in $\Fr/ \mathcal{H}$. Lemma~\ref{cor: no_accumulation_of_powers2} implies that $\tilde{A}_j=\tilde{A}_r$ for all ranks $j \geq r$. Hence $\tilde{A}_r^{\tau}$ is a subword of words from $\Can_i$ for all $i$. Thus by Lemma~\ref{power_subwords_of_can_all_ranks} $\tilde{A}_r^n\in \mathcal{H}$ and hence $W^n \in\mathcal{H}$, a contradiction.
\end{proof}

Since the sets $\Rel_i, i\geq 0,$ consist of $n$-th powers, we now obtain:
\begin{corollary}
\label{rel_covers_all_powers}
The normal subgroup of $\Fr$ generated by $\llangle \ \Rel_i\mid\ i\in\mathbb{N}\rrangle $ coincides with the normal subgroup generated by all $n$-th powers.
\end{corollary}

\begin{theorem}
\label{final_can_of_equivalent_words}
For every $A, B \in \Can_{-1}$ the words $A$ and $B$ represent the same element of the group $B(m, n)$ if and only if $\can(A) = \can(B)$.
\end{theorem}
\begin{proof}
If $\can(A) = \can(B)$, then   $A\equiv B \mod \llangle\ \Rel_i\mid\ i\geq 0\ \rrangle$ by Corollary~\ref{final_can_idempotent}. Thus clearly $A$ and $B$ represent the same element of the group $B(m, n)$.

Converseley, if $A$ and $B$ represent the same element in $B(m, n)$, then, by definition, $A\equiv B \mod \llangle w_1^n, \ldots , w_k^n \rrangle $ for some cyclically reduced words $w_i$. By Corollary~\ref{rel_covers_all_powers} we have $w_i^n \in \llangle  \Rel_0, \Rel_1, \ldots, \Rel_r\rrangle $ for some $r$  and so $A \equiv B \mod \llangle  \Rel_0, \ldots, \Rel_r\rrangle $. Thus  $\can_r(A) = \can_r(B)$ by IH~\ref{IH can of equivalent words} and, by construction, $\can(A) = \can(B)$.
\end{proof}

Finally we are ready to prove Theorem~\ref{main thm}:

\begin{remark}
\label{final_can_immediately_condition}
If $A \in \Can_0$ and all subwords of $A$ of the form $a^ka_1$, $a = a_1a_2$, satisfy  $\Lambda_a(a^ka_1) \leq \frac{n}{2} - 5\tau - 2$, then by iterated application of IH~\ref{IH immediate} we have $A \in \cap_{i = 0}^{\infty}\Can_i$.
\end{remark}

\begin{proof}[The proof of Theorem~\ref{main thm}]
By Theorem~\ref{final_can_of_equivalent_words} words $A, B \in \Can_{-1}$ represent the same element in $B(m, n)$ if and only if $\can(A) = \can(B)$. 
Now consider the set of cube free words in $\Can_0$, which is an infinite set by~\cite{Fife}. Clearly for our choice of the exponent~$n$ we have $3 \leq  \frac{n}{2} - 5\tau - 2$ and hence by Remark~\ref{final_can_immediately_condition} every cube free word is in $\Can$. Thus $\Can$ is infinite and hence so is $B(m, n)$.
\end{proof}

\end{document}